\documentclass{article}

\usepackage{lipsum}
\usepackage{epstopdf}
\usepackage{longtable}
\usepackage{amsmath, amsfonts, bm, amssymb}
\usepackage{amsthm}
\usepackage{graphicx, adjustbox}
\usepackage{bbm}
\usepackage{algorithm, algpseudocode}
\usepackage{hyperref}

\usepackage{cleveref}

\newcommand{\bC}{\mathbb C}

\newcommand{\bR}{\mathbb R}

\newcommand{\tone}{\mathtt 1}

\DeclareMathOperator*{\dgets}{. \! \gets}
\DeclareMathOperator*{\dtimes}{. \! \ast}
\DeclareMathOperator*{\ddiv}{. \! /}

\newcommand{\rd}{\mathrm d}
\newcommand{\OT}{\mathrm{OT}}
\newcommand{\LOT}{\mathrm{LOT}}

\newcommand{\Om}{\Omega}

\newcommand{\cP}{\mathcal P}
\newcommand{\cM}{\mathcal M}
\newcommand{\cC}{\mathcal C}

\DeclareMathOperator*{\argmax}{arg\,max}
\DeclareMathOperator*{\argmin}{arg\,min}

\newtheorem{thm}{Theorem}
\newtheorem{prop}{Proposition}
\newtheorem{defi}{Definition}
\newtheorem{rem}{Remark}
\newtheorem{examp}{Example}

\title{A registration method for reduced basis problems using linear optimal transport\thanks{Submitted to SIAM Journal on Scientific Computing}}

\author{Tobias Blickhan$^{1,2}$ }
\date{   $^1$ Max Planck Institute for Plasma Physics \\ $^2$  Technische Universit\"at M\"unchen }

\usepackage{amsopn}

\usepackage{geometry}
\ifpdf
\hypersetup{
  pdftitle={A registration method for reduced basis problems using linear optimal transport},
  pdfauthor={Tobias Blickhan}
}
\fi

\begin{document}

\maketitle

% REQUIRED
\begin{abstract}
	We present a registration method for model reduction of parametric partial differential equations with dominating advection effects and moving features. Registration refers to the use of a parameter-dependent mapping to make the set of solutions to these equations more amicable for approximation using classical reduced basis methods. The proposed approach utilizes concepts from optimal transport theory, as we utilize Monge embeddings to construct these mappings in a purely data-driven way. The method relies on one interpretable hyper-parameter. We discuss how our approach relates to existing works that combine model order reduction and optimal transport theory. Numerical results are provided to demonstrate the effect of the registration. This includes a model problem where the solution is itself a probability density and one where it is not.
\end{abstract}

\section{Introduction}

The field of reduced complexity modelling is of practical interest for numerous multi-query and real-time applications. Reduced methods can be constructed in an offline (or training) phase, leveraging the strengths of high-performance computing infrastructure. After construction reduced models can reduce the time, money, and energy spent on optimization loops, inverse problems, or routine calculations.
%After construction, they can be deployed to personal computers, handheld devices, or even embedded systems to provide fast approximations at very low computational cost. Using reduced models in optimization loops, inverse problems, or routine calculations can greatly reduce time, energy consumption, and cost of performing these tasks.

Classical methods such as the reduced basis approach \cite{quarteroni_reduced_2016,ohlberger_reduced_2016} provide ways to save orders of magnitude in computational cost while at the same time ensuring rigorous error bounds for the reduced simulation of elliptic and parabolic equations. However, they are notoriously ill-suited when working with hyperbolic systems or solutions with moving features and sharp discontinuities.

In this work, we present a registration method for model order reduction based on methods from optimal transport theory. %For advection-dominated problems, a common and successful approach is to perform a registration step to the data before applying a linear reduction method such as proper orthogonal decomposition \cite{rim_manifold_2023,taddei_registration_2020,nair_transported_2019,iollo_advection_2014}. 
The role of the registration step is to align the features in the solution data to make it more amicable for linear approximation. The method we propose is fully data-driven and therefore requires no previous knowledge of the form of the alignment mapping. Bijectivity of the mapping is connected to its construction from optimal transport theory and not enforced with penalty terms. One hyper-parameter set by the user controls the trade-off between accuracy and smoothness of the mapping. In an effort to keep the presentation self-contained, we recall some results from optimal transport theory and research that connects it to reduced complexity modeling.

\section{Outline}

The article is structured as follows: In \cref{sec:RB}, we briefly review reduced basis and registration methods. \Cref{sec:OT} contains the concepts from optimal transport theory that will be used in our proposed method, which is described in \cref{sec:method}. \Cref{sec:other_works} summarizes similar approaches from the literature and points out similarities and differences to our work. Two numerical examples are given in \cref{sec:numerical_examples}. Specifics of computations related to optimal transport are presented in the supplementary material, \cref{sec:computation}. 
\section{The reduced basis method}\label{sec:RB}

In this work, we will follow the classical reduced basis (RB) approach for parametrized partial differential equations (PPDEs) \cite{quarteroni_reduced_2016}. Given parameters $\mu \in \mathcal A$, we seek to solve many iterations of the following problem: Find a function $u \in V_h$, a discretized function space, such that
\begin{equation}
	\mathcal{L}(u; \mu) = 0 \quad \text{in } \Omega,
\end{equation}
together with suitable initial and boundary conditions which can again depend on $\mu$.  
\begin{rem}
	In this work, we will treat the solution to the discretized problem as the true solution to the PPDE problem. Consequently, we omit the subscript $_h$ from $u$ to avoid notational clutter. Many of the following concepts are also well-defined if $V_h$ is replaced by any Hilbert space.
\end{rem}
The set of all possible solutions to this problem is called the \emph{solution manifold} 
\begin{equation}\label{eq:solution_manifold}
	\cM := \{ u(\mu) \in V_h  : \mathcal{L}(u(\mu); \mu) = 0 \text{ where } \mu \in \mathcal A \}.
\end{equation}
If it is possible to approximate $\cM$ by a linear subspace of small dimension, the PPDE can be solved quickly and accurately for many values of $\mu$ within this subspace. This property of $\mathcal{M}$ is referred to as its \emph{n-width} (\cite{quarteroni_reduced_2016}, Section 5.4):
\begin{defi}[n-width]\label{def:n-width}
	The \emph{(Kolmogorov) n-width} of $\mathcal{M}$ is given by
	\begin{equation}\label{eq:n-width}
		\inf_{\substack{ V_n \subset V_h \\ \dim V_n = n }} \sup_{u \in \mathcal M} \inf_{u_{\mathrm{rb}} \in V_n} \Vert u_{\mathrm{rb}} - u \Vert_{V_h}.
	\end{equation}
\end{defi}
An indicator of the n-width (a worst-case error measure) can be obtained by sampling $\cM$ and using an eigenvalue decomposition to estimate an average error. This method is known \emph{proper orthogonal decomposition} (POD)
\begin{defi}\label{def:POD}
	The $i$th POD basis element of the set $\{ u(\mu_i) \}_{i=1}^{n_s}$ is given by
	\begin{equation}\label{eq:pod_basis}
		\zeta_i = (\lambda^u_i)^{-\frac{1}{2}} \sum_{j=1}^{n_s} u(\mu_j) (\mathrm{v}^u_i)_j,
	\end{equation}
where $\lambda^u_i$ and $\mathrm{v}^u_i$ are the $i$th eigenvalue and eigenvector of the snapshot correlation matrix
\begin{equation}
	\bC^u := \{ \langle u(\mu_i), u(\mu_j) \rangle_{V_h} \}_{1 \leq i,j \leq n_s}.
\end{equation}
\end{defi}
According to the Eckart–Young–Mirsky theorem, the space spanned by the POD basis of size $n$ is optimal in the sense that $\sum_{i=1}^n \Vert u(\mu_i) - \Pi_{\mathrm{POD}} u(\mu_i) \Vert^2_{V_h} = \sum_{i>n} \lambda_i$ is minimal, where $\Pi_{\mathrm{POD}}$ denotes the orthogonal projection onto $\mathrm{span} \, \{ \zeta_i \}_{i=1}^n$. Consequently, the decay of the eigenvalues of $\bC^u$ is a good indicator for how well a problem can be approximated using the classical RB method. While this decay can be shown to be exponential for, e.g., regular elliptic problems (\cite{quarteroni_reduced_2016}, Section 5.5), this decay can be very slow for advection-dominated problems or problems with parameter-dependent geometry. This issue is well-known. One common strategy to overcome it is to find a suitable, parameter dependent mapping $\Phi_\mu$ such that the manifold of mapped solutions
\begin{equation}\label{eq:mapped_manifold}
	\Phi_\mu ( \cM ) :=  \{ u(\mu) \circ \Phi^{-1}_\mu : u(\mu) \in \cM \}
\end{equation}
has a much smaller n-width \cite{rim_manifold_2023,taddei_registration_2020,nair_transported_2019,chetverushkin_model_2019,welper_interpolation_2017,iollo_advection_2014}. This approach we refer to as \emph{registration methods}. It leads to two challenges: First, a suitable mapping must be constructed (offline) and evaluated (online). Second, the mapped problem must be solved (online). The present work will focus on the first challenge, but we will also address how we tackle the second.
\section{Optimal transport}\label{sec:OT}

Optimal transport (OT) theory provides a notion of distance between probability measures $\rho, \sigma \in \mathcal P(\Omega_1) \times \mathcal{P}(\Omega_2) $. In particular, a cost is modelled through a function $c: \Omega_1 \times \Omega_2 \rightarrow \bR$ where $c(x,y)$ gives the cost of moving mass from $x$ to $y$. In this work, we will look at the case where $\Om_1 = \Om_2 = \Om \subset \bR^d$. In all cases, we use the quadratic cost function $c(x,y) = \tfrac{1}{2} \vert x - y \vert^2$. The goal is to move mass that is distributed according to $\rho$ to the configuration given by $\sigma$, while minimizing the total cost. Splitting mass is allowed, which leads to the following minimization problem:

\begin{defi}[OT distance]\label{def:OT_dist}
	The OT distance between $\rho, \sigma \in \cP(\Omega)$ with $\Om \subset \bR^d$ compact is given by
	\begin{equation}\label{eq:OT_primal}
		\OT(\rho,\sigma)^2 := \inf_{\pi \in \Pi(\rho,\sigma)} \int_{\Omega \times \Omega} c(x,y) \rd \pi(x,y),
	\end{equation}
	where $\Pi(\rho,\sigma)$ is the set of admissible \emph{transport plans}, i.e. probability measures on $\Omega \times \Omega$ with fixed marginals:
	\begin{equation}\label{eq:adm_plans}
		\Pi(\rho,\sigma) := \{ \pi \in \cP(\Om \times \Om) : \pi(\cdot,\Om) = \rho \text{ and } \pi(\Om,\cdot) = \sigma \}.
	\end{equation}
\end{defi}

The OT distance is also known as the \emph{Wasserstein distance} or metric. A better name considering its history would be the \emph{Kantorovic metric} \cite{vershik_long_2013}.

A thorough introduction into the topic of optimal transport is beyond the scope of this work. We will repeat only a few select results for convenience. Proofs to the given propositions and theorems can be found in the excellent textbooks on the topic \cite{peyre_computational_2019, villani_topics_2016, santambrogio_optimal_2015}. 

\subsection{Transport potentials and c-transform}

\begin{prop}[OT duality]
	The dual problem of \cref{eq:OT_primal} is given by
	\begin{equation}\label{eq:OT_dual}
		\OT(\rho,\sigma)^2 = \sup_{\psi_\rho, \psi_\sigma \in \cC_b(\Om)} \left \{ \int_\Om \psi_\rho \rd \rho + \int_\Om \psi_\sigma \rd \sigma \; : \; \psi_\rho(x) + \psi_\sigma(y) \leq c(x,y) \right \},
	\end{equation}
	where $\cC_b$ denotes the space of bounded continuous functions. The duality gap is zero, so the two optimization problems are equivalent.
\end{prop}

\begin{defi}[Transport potentials]\label{def:potentials}
	The functions $\psi_\rho, \psi_\sigma \in \cC_b(\Om)$ from \cref{eq:OT_primal} we call \emph{transport potentials}. The name \emph{Kantorovich potentials} can also be found in the literature to refer to those potentials solving the dual problem.
\end{defi}

 Given that $\rho$ and $\sigma$ are probability measures, to maximise the objective in \cref{eq:OT_dual}, we want to make the potentials $\psi_\rho, \psi_\sigma$ as large as possible without violating the constraint $\psi_\rho \oplus \psi_\sigma \leq c$. This leads to the notion of \emph{c-transform}:

\begin{defi}[c-transform]\label{def:c-transform}
	The c-transform of a function $\psi: \Omega \rightarrow \bR \cup \{ +\infty \}$ is given by
	\begin{equation}\label{eq:c-transform}
		\psi^c(y) := \inf_{x \in \Omega} \left ( c(x,y) - \psi(x) \right ).
	\end{equation}
\end{defi}

\begin{prop}
	The potentials solving \cref{eq:OT_dual} are c-conjugate to each other, i.e. $\psi_\rho(x) = \psi_\sigma^c(x)$.
\end{prop}

\begin{rem}
	In the case $c(x,y) = \tfrac{1}{2} \vert x - y \vert^2$, the c-transform is related to the Legendre transform by
	\begin{equation}\label{eq:psi_phi_legendre}
		\frac{\vert y \vert^2}{2} - \psi^c(y) = \frac{\vert y \vert^2}{2} - \inf_{x \in \Omega}  \left ( \frac{1}{2} \vert x-y \vert^2 - \psi(x) \right ) = \left ( \frac{\vert x \vert^2}{2} - \psi(x) \right )^* =: \varphi^*(y).
	\end{equation}
\end{rem}
If the measures are absolutely continuous with respect to the Lebesge measure, there is no mass splitting, i.e. there exists a \textit{transport map} $T$ such that $\sigma = T_\sharp \rho$:
\begin{defi}[Push-forward]\label{def:pushfwd}
	The push-forward of $\rho$ under $T$, denoted by $T_\sharp \rho$, is defined by $(T_\sharp \rho)(\Omega') = \rho(T^{-1}(\Omega'))$ for all measureable $\Omega' \subseteq \Omega$. If both $\rho$ and $\sigma$ are absolutely continuous with densities denoted again $\rho,\sigma$, and $T$ is a $C^1$ diffeomorphism, then $\sigma = T_\sharp \rho$ is equivalent to
	\begin{equation}\label{eq:push-forward}
		\rho = \left ( \sigma \circ T \right) \, | \det DT |.
	\end{equation}
\end{defi}

\begin{thm}[Brenier's theorem]\label{thm:brenier}
	Assume $\Omega \subset \bR^d$ compact. If at least one of $\rho,\sigma$ has a density with respect to the Lebesgue measure, then the unique solution to \cref{eq:OT_primal} is concentrated on the graph $(x, T(x))$ of the \emph{transport map} $T$. In particular,
	\begin{equation}
		\inf_{\pi \in \Pi(\rho,\sigma)} \frac{1}{2} \int_{\Omega \times \Omega} \vert x - y \vert^2 \rd \pi(x,y) = \inf_{T : \sigma = T_\sharp \rho} \frac{1}{2} \int_{\Omega} \vert T(x) - x \vert^2 \, \rd \rho(x),
	\end{equation}
	and $T$ is the gradient of a convex function:
	\begin{equation}\label{eq:Transport_map}
		T(x) = \nabla \varphi = x - \nabla \psi(x).
	\end{equation}
	The function $\psi$ one of the two optimal transport potentials solving \cref{eq:OT_dual}, i.e. $\psi_\rho$ corresponds to $T_{\rho \rightarrow \sigma}$ and vice versa.
\end{thm} 

The map $T$ is also called the \textit{Monge map}.

\begin{rem}
	The push-forward density is also sometimes denoted $\rho \circ T^{-1}$. In this notation, $\sigma = \rho \circ T_{\rho \rightarrow \sigma}^{-1} = \rho \circ ( \mathrm{id} - \nabla \psi_{\rho \rightarrow \sigma} )^{-1} = \rho \circ ( \mathrm{id} - \nabla \psi_{\rho \rightarrow \sigma}^c )= \rho \circ ( \mathrm{id} - \nabla \psi_{\sigma \rightarrow \rho} )$. We denote by $\psi_{\rho \rightarrow \sigma}$ the potential such that $(\mathrm{id} - \nabla \psi_{\rho \rightarrow \sigma} )_\sharp \rho = \sigma$.
\end{rem}

\subsection{Displacement interpolation and OT barycenters}

\begin{figure}[h!]
	\centering
		\includegraphics[width=0.25\textwidth]{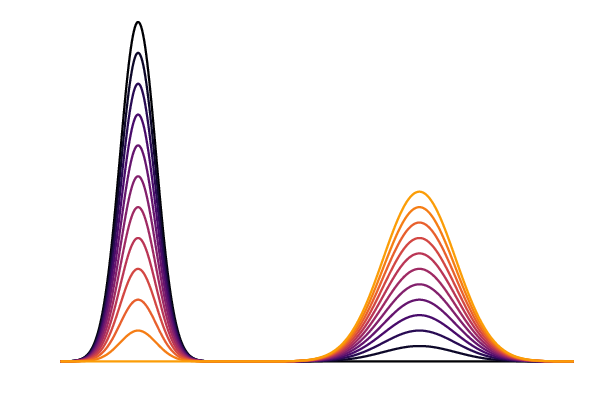}
		\includegraphics[width=0.25\textwidth]{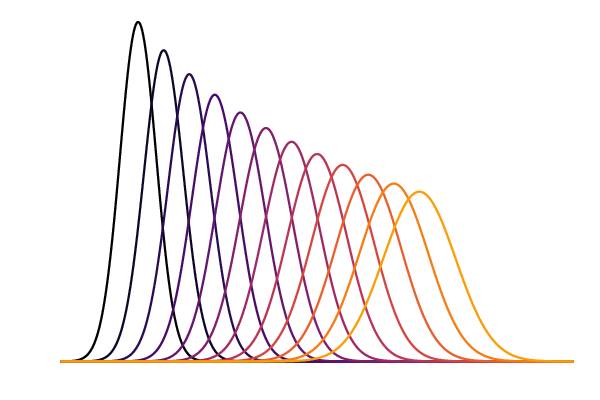}
	\caption{$L^2$ interpolation and displacement interpolation between two Gaussian distributions.}
	\label{fig:displacement_interpolation}
\end{figure}

Given the map $T$, one can visualize the optimal transport by introducing a time-like variable $t$ to define $T_t(x) := (1-t)x + t T(x)$ and consider $\rho_t := (T_t)_\sharp \rho$. The resulting path $t \mapsto \rho_t$ through the space of probability measures is referred to as \emph{displacement interpolation} \cite{mccann_convexity_1997} and illustrated in \cref{fig:displacement_interpolation}.

The generalization to convex combinations of several distributions is called an \textit{OT barycenter}:

\begin{defi}[OT barycenter]\label{def:OT_barycenter}
	Given a set of $\{\rho_i\}_{i=1}^n \subset \cP(\Omega)$ and non-negative, normalized weights $\{\omega_i\}_{i=1}^n$, the OT barycenter of $\{\rho_i\}_{i=1}^n$ is given by
	\begin{equation}\label{eq:OT_barycenter}
		\mathrm{OTBar}( \{ \omega_i, \rho_i \}_{i=1}^n ) := \argmin_{\sigma \in \cP(\Omega)} \sum_{i=1}^n \omega_i \, \OT(\rho_i, \sigma)^2.
	\end{equation}
	When $\omega_i \equiv n_s^{-1} \; \forall i$, we omit the weights and write $\mathrm{OTBar}( \{ \rho_i \}_{i=1}^n )$.
\end{defi}

\subsection{Linear optimal transport}

OT theory provides a metric on the space of probability measures with many favourable properties (e.g. the OT distance metrizes weak convergence for compact $\Omega \subset \bR^d$). Moreover, it allows us to (formally) treat $\cP(\Omega)$ as a Riemannian manifold, where the elements of the tangent space of $\rho \in \cP(\Omega)$ are elements of $L^2(\rho, \bR^d)$. This allows a geometric and physical interpretation for the transport potentials $\psi$, which are the central objects of method. For more information on this geometric interpretation, we refer to \cite{otto_geometry_2001} and \cite{santambrogio__2016, ambrosio_gradient_2005}.

Given two distributions $\rho$ and $\sigma$, the optimal transport map $T_{\rho \rightarrow \sigma}$ from $\rho$ to $\sigma$ is an element of the tangent space $L^2(\rho, \bR^d)$. Linear optimal transport (LOT) \cite{wang_linear_2013, merigot_quantitative_2020, moosmuller_linear_2023} works with these vector spaces. We recall the following definitions:

\begin{defi}[Monge embedding]\label{def:monge_embedding}
	Given a fixed reference distribution $\bar \rho \in \cP(\Omega)$, let $T_{\bar \rho \rightarrow \sigma}$ denote the Monge map from $\bar \rho$ to $\sigma$. The mapping $\sigma \mapsto T_{\bar \rho \rightarrow \sigma}$ is called \emph{Monge embedding}.
\end{defi}

\begin{defi}[LOT distance]\label{def:lot_distance}
	\begin{equation}\label{eq:LOT_distance}
		\LOT^{\bar \rho}(\rho_1,\rho_2)^2 := \int_\Omega \vert T_{\bar \rho \rightarrow \rho_1} - T_{\bar \rho \rightarrow \rho_2} \vert^2 \, \rd \bar \rho.
	\end{equation} 
\end{defi}

\begin{defi}[LOT barycenter]\label{def:lot_barycenter}
	\begin{equation}\label{eq:LOT_barycenter}
		\mathrm{LOTBar}^{\bar \rho}( \{ \omega_i; \rho_i \}_{i=1}^n) := \argmin_{\sigma \in \cP(\Omega)} \sum_{i=1}^n \omega_i \, \LOT^{\bar \rho}(\rho_i, \sigma)^2 = \left ( \sum_{i=1}^n \omega_i \, T_{\bar \rho \rightarrow \rho_i} \right )_\sharp \bar \rho.
		% = \argmin_{\sigma \in \cP(\Omega)} \sum_{i=1}^n \omega_i \, \LOT^{\bar \rho}(\rho_i, \sigma)^2
	\end{equation}
\end{defi}

A natural question to ask is how well the LOT distance approximates the true OT distance. For general measures $\rho, \sigma, \bar \rho$, it holds that $\OT(\rho,\sigma) \leq \LOT(\rho,\sigma) \leq \mathrm{const.}(\Omega, d) \OT(\rho,\sigma)^{\tfrac{2}{15}}$ (\cite{merigot_quantitative_2020}, Theorem 3.1), so the Monge embedding is not an isometry. With regularity assumptions on the measures, the exponent can be improved to $\tfrac{1}{2}$ for measures along a Lipschitz curve through 2-Wasserstein space \cite{gigli_holder_2011}.

There are, however, special cases. For example, if $\rho$ and $\sigma$ are connected by shifts and scalings, i.e. a transformation of the form $x \mapsto \bar b x +\bar a$ where $0 < \bar b \in \bR$ and $\bar a \in \bR^d$, then $\OT = \LOT$. Furthermore, if the transformation connecting $\rho$ to $\sigma$ is $\varepsilon$-close to the form $\bar b x + \bar a$, then $\LOT(\rho,\sigma) \leq \OT(\rho,\sigma) + C_1 \varepsilon + C_2 \sqrt \varepsilon$. We refer to \cite{moosmuller_linear_2023} for more details.

\begin{rem}
	The examples shown in \cref{sec:numerical_examples} are fundamentally scaling and translation mappings. In this context, OT and linear OT can be expected to perform well. For more complex transformations, and for transported measures that are not smooth densities (e.g.: curves), the fact that OT mappings do not preserve topology can become an issue. This difficulty is discussed in the context of shape analysis and medical imaging in \cite{feydy_geometric_2020} (see especially figures 3.31 and 3.32).
\end{rem}

\subsection{Entropic optimal transport}

The classical OT problem can be modified by adding a regularizing term to the cost function \cref{eq:OT_primal} as popularized by \cite{cuturi_sinkhorn_2013} (see \cite{leonard_survey_2014} for a historical perspective going back to statistical physics considerations from the early 20th century). The resulting regularized problem can be solved numerically in a very efficient way (see \cref{sec:computation}).
\begin{defi}[Entropic OT]\label{def:entropic_ot}
	Let $\rho,\sigma, c,\Pi(\rho,\sigma)$ as in \cref{def:OT_dist} and $\varepsilon > 0$. The OT problem with entropic regularization reads
	\begin{multline}\label{eq:OT_entropic_primal}
		\OT_\varepsilon(\rho,\sigma)^2 := \min_{\pi^\varepsilon \in \Pi(\rho,\sigma)} \int_{\Omega \times \Omega} c(x,y) \rd \pi^\varepsilon(x,y) + \varepsilon \int_{\Omega \times \Omega}  \log \left ( \frac{\rd \pi^\varepsilon(x,y)}{\rd \rho(x) \, \rd \sigma(y)} \right ) \rd \pi^\varepsilon(x,y) \\ - \varepsilon \int_{\Omega \times \Omega} \rd \pi^\varepsilon(x,y) + \varepsilon \int_{\Omega \times \Omega} \rd \rho(x) \, \rd \sigma(y).
	\end{multline}
	The corresponding dual problem has the form
	\begin{multline}\label{eq:OT_entropic_dual}
		\OT_\varepsilon(\rho,\sigma)^2 = \max_{\psi^\varepsilon_\rho, \psi^\varepsilon_\sigma \in \cC_b(\Omega)} \int_\Om \psi^\varepsilon_\rho(x) \rd \rho(x) + \int_\Om \psi^\varepsilon_\sigma(y) \rd \sigma(y)  \\ - \varepsilon \int_{\Omega \times \Omega} \exp \left ( \frac{\psi_\rho^\varepsilon(x) + \psi^\varepsilon_\sigma(y) -c(x,y)}{\varepsilon} \right ) \rd \rho(x)  \rd \sigma(y) + \varepsilon.
	\end{multline}
\end{defi}
Note that we recover the constraint $\psi_\rho \oplus \psi_\sigma \leq c$ as $\varepsilon \rightarrow 0$ in \cref{eq:OT_entropic_dual}. A formal computation shows that the stationarity conditions for this problem read
\begin{align}\label{eq:schroedinger_eq}
	\exp \left (\frac{- \psi^\varepsilon_\sigma(y)}{\varepsilon} \right ) &= \int_\Omega \exp \left (\frac{\psi^\varepsilon_\rho(x) - c(x,y)}{\varepsilon} \right ) \rd \rho(x) \quad \sigma - \text{a. e.} \quad \text{and} \\
	\exp \left (\frac{- \psi^\varepsilon_\rho(x)}{\varepsilon} \right ) &= \int_\Omega \exp \left (\frac{\psi^\varepsilon_\sigma(y) - c(x,y)}{\varepsilon} \right ) \rd \sigma(y) \quad \rho - \text{a. e.}.
\end{align}
\begin{rem}
	At optimality, the double integral in \cref{eq:OT_entropic_dual} evaluates to one due to \cref{eq:schroedinger_eq} and $\OT_\varepsilon(\rho,\sigma)^2$ is given by the sum of two weighted integrals of the potential functions.
\end{rem}
%\begin{defi}[Scaling functions]\label{def:scaling_functions}
%	The functions $a: x \mapsto \exp \left ( \tfrac{1}{\varepsilon} %\psi^\varepsilon_\rho(x) \right )$ and $b: y \mapsto \exp \left ( \tfrac{1}{\varepsilon} %\psi^\varepsilon_\sigma(y) \right )$ are called \emph{scaling functions}.
%\end{defi}
Taking the logarithm of \cref{eq:schroedinger_eq} defines the \emph{softmin}, which replaces the \emph{c-transform}:
\begin{defi}[softmin]\label{def:softmin}
	\begin{align}\label{eq:softmin}
			\psi^{c,\varepsilon}(y) &= - \varepsilon \log \int_\Omega \exp \left (\frac{\psi^\varepsilon(x) - c(x,y)}{\varepsilon} \right ) \rd \rho(x) \nonumber \\
			&=:  \min^\varepsilon_{x \sim \rho} \left \{ c(x,y) - \psi^\varepsilon(x) \right \}
	\end{align}
\end{defi}
\begin{rem}
	The fundamental difference to the classical OT problem is that for any $\varepsilon > 0$, the solution to \cref{eq:OT_entropic_primal} necessarily has full support with respect to the product measure $\rho \otimes \sigma$. In contrast, it was concentrated on the graph of the Monge map for $\varepsilon = 0$ in cases covered by \cref{thm:brenier}. The optimal solution lies in the interior of the set of admissible solutions and is characterized by optimality conditions of vanishing gradients \cite{feydy_geometric_2020}.
\end{rem}
%\begin{rem}
%	The primal $OT_\varepsilon$ problem can be (formally) reformulated as
%		\begin{equation}
%		\OT_\varepsilon(\rho,\sigma)^2 = \min_{\pi^\varepsilon \in \Pi(\rho,\sigma)} \varepsilon \int_{\Omega \times \Omega} \log \left ( \frac{1}{k^\varepsilon(x,y)}  \frac{\rd \pi^\varepsilon(x,y)}{\rd \rho(x) \, \rd \sigma(y)} \right ) \rd \pi^\varepsilon(x,y),
%	\end{equation}
%	where $k^\varepsilon(x,y) := \exp \left ( - \tfrac{1}{\varepsilon} c(x,y) \right )$ is called the \emph{Gibbs kernel}. This formulation is known as the \emph{Schr\"odinger bridge problem} and, as the name suggests, 
%\end{rem}
The parameter $\varepsilon$ determines the strength of the regularization. As $\varepsilon \rightarrow 0$, $\min^\varepsilon_{x \sim \rho} \rightarrow \min_{x \in \mathrm{supp} \, \rho}$ \cite{feydy_geometric_2020}. Furthermore, it holds that $\Vert \nabla^k \psi \Vert_\infty = \mathcal{O}(1 + \varepsilon^{1-k})$ \cite{genevay_sample_2019}. The properties of $\OT_\varepsilon$, including its convergence as $\varepsilon$ goes to zero has been studied extensively. We refer to \cite{peyre_computational_2019, feydy_geometric_2020} for details. A useful practical interpretation is this: as the softmin operation is a Gaussian convolution, the entropic transport plan practically ignores features below the scale of $\sqrt{\varepsilon}$.

Just as in \cref{eq:psi_phi_legendre}, the case with no regularization, the potentials are linked to convex functions, defined through an approximate maximum:
\begin{prop}\label{prop:convexity_softmax}
	The function
	\begin{align}
		y \mapsto \frac{1}{2} |y|^2 - \psi^{c,\varepsilon}(y) &= \varepsilon \log \int_\Omega \exp \left ( \frac{1}{\varepsilon} \left ( x \cdot y - \frac{1}{2} | x |^2 + \psi^\varepsilon(x) \right ) \right ) \rd \rho(x) \nonumber \\
		&=: \max_{x \sim \rho}^\varepsilon \left \{ x \cdot y - \left ( \frac{1}{2} |x|^2 - \psi^\varepsilon(x) \right ) \right \}
	\end{align} is convex.
\end{prop}
\begin{proof}
	Evaluate the function at $y_t := ty_1 + (1-t)y_2 : 0 < t < 1$ and apply H\"older's inequality with exponents $\tfrac{1}{t}, \tfrac{1}{1-t}$.
\end{proof}
\begin{rem}\label{rem:convexity_softmax}
	The particular value of $\varepsilon$ plays no role for convexity. The limit cases $\varepsilon \rightarrow 0: \max^\varepsilon_{x \sim u} \rightarrow \max_{x \in \mathrm{support} \, \rho}$ and $\varepsilon \rightarrow \infty: \max^\varepsilon_{x \sim \rho} \rightarrow \int_\Omega \, \rd \rho$ are also convex.
\end{rem}
In practice, the $\OT_\varepsilon$ dual problem is solved by iteratively applying the softmin operation until the potentials fulfill \cref{eq:schroedinger_eq} and $\psi_\rho^\varepsilon = \psi_\sigma^{c,\varepsilon}$. The solution to the primal problem is then given by $\pi^\varepsilon = \rho \otimes \sigma \, \exp \left ( \frac{1}{\varepsilon} \left ( \psi_\rho^\varepsilon \oplus \psi_\sigma^\varepsilon - c \right ) \right )$. \\
\begin{rem}\label{rem:sinkhorn_cost}
	The number of iterations needed to solve the entropic OT problem in practice go up dramatically as $\varepsilon \rightarrow 0$. In particular, in simple cases where the solution for $\varepsilon = 0$ is a smooth map, the error after the $l$th iteration is of order $(1-\varepsilon)^l$ (\cite{peyre_computational_2019}, Remark 4.15). For the moderately small values of $\varepsilon$ (when compared to the characteristic scale of the cost function) used in this work, this is not yet very restrictive. Beyond that, however, simulated annealing and multi-scale methods become necessary (c.f. \cite{feydy_geometric_2020}, Section 3.3.3 and \cite{peyre_computational_2019}, Section 4.2).
\end{rem}
The entropic OT problem does not admit a transport map as a solution, as the transport plan is necessarily supported on the entirety of $\rho \otimes \sigma$. It is a natural question what the map $x \mapsto x - \nabla \psi_\rho^\varepsilon(x)$ corresponds to. From the stationarity condition \cref{eq:schroedinger_eq}, we find
\begin{align}
	\nabla \psi_\rho^\varepsilon(x) &= \frac{\int_\Omega (x - y) \exp \left ( \frac{1}{\varepsilon}  \left ( \psi_\sigma^\varepsilon(y) - c(x,y) \right ) \right ) \rd \sigma(y)  }{ \int_\Omega \exp \left ( \frac{1}{\varepsilon}  \left ( \psi_\sigma^\varepsilon(y) - c(x,y) \right ) \right ) \rd \sigma(y)  } \nonumber \\
	&= x -  \frac{ \int_\Omega y \exp \left ( \frac{1}{\varepsilon}  \left ( \psi_\sigma^\varepsilon(y) - c(x,y) \right ) \right ) \rd \sigma(y)   }{\int_\Omega \exp \left ( \frac{1}{\varepsilon}  \left ( \psi_\sigma^\varepsilon(y) - c(x,y) \right ) \right ) \rd \sigma(y) } \nonumber \\
	&=: x - T^\varepsilon_{\rho \rightarrow \sigma}(x)
\end{align}
\begin{defi}[Entropic Monge map]
	We call $T_{\rho \rightarrow \sigma}^\varepsilon = \mathrm{id} - \nabla \psi_\rho^\varepsilon$ the \emph{entropic Monge map} between $\rho$ and $\sigma$.
\end{defi}
It must be stressed that $T^\varepsilon_\sharp \rho \neq \sigma$ in general. Nevertheless, the map has appealing properties: It is defined for all $y \in \Omega$ (not only $\rho$ - almost everywhere) and converges to the Monge map in $L^2(\rho)$ as $\varepsilon \rightarrow 0$ \cite{pooladian_entropic_2022}. The entropic Monge map can also be interpreted as an extension of the expected value $x \mapsto \mathbb{E}_{\pi^\varepsilon} [ Y | X = x ]$ from $\mathrm{support}(\rho)$ to the entire domain. It is also referred to as the \emph{barycentric mapping} \cite{feydy_geometric_2020} as it has the form of a weighted mean with normalized weights.

\subsection{Sinkhorn divergences}\label{sec:sinkhorn_dvg}

While the OT problem without regularization defines a distance on $\mathcal{P}$, this does not hold for the entropic version $\OT_\varepsilon$. In general, $\OT_\varepsilon(\rho,\rho) \neq 0$ and there exist measures $\sigma$ such that $\OT_\varepsilon(\rho,\sigma) < \OT_\varepsilon(\rho, \rho)$ \cite{feydy_geometric_2020}. Indeed, as stated in \cite{janati_debiased_2020}, Proposition 1, we have
\begin{equation}
	\OT_\varepsilon(\rho,\rho)^2 =  \max_{\psi \in \cC_b(\Omega)} \left ( 2 \int_{\Omega} \psi \, \rd \rho - \varepsilon \int_{\Omega \times \Omega} \exp \left ( \frac{\psi \oplus \psi - c }{\varepsilon} \right ) \rd \rho \, \rd \rho + \varepsilon \right ) = 2 \int_{\Omega} \psi_{\rho \rightarrow \rho}^\varepsilon \rd \rho
\end{equation}
where $\psi_{\rho \rightarrow \rho}^\varepsilon$ solves $e^{- \tfrac{\psi_{\rho \rightarrow \rho}^\varepsilon(x)}{\varepsilon}} = \int_{\Omega} e^{ \tfrac{\psi_{\rho \rightarrow \rho}^\varepsilon(y)  - c(x,y) }{\varepsilon} } \rd \rho(y) \; \rho - a.e.$.
\begin{rem}
	Consequently, fitting problems of the form $\min_\sigma \OT_\varepsilon(\rho,\sigma)$ usually do not converge to $\rho$.
\end{rem}
As discussed in \cite{ramdas_wasserstein_2017, genevay_learning_2018}, the \emph{Sinkhorn divergence}
\begin{equation}\label{eq:sinkhorn_dvg}
	S_\varepsilon(\rho,\sigma)^2 := \OT_\varepsilon(\rho,\sigma)^2 - \frac{1}{2} \OT_\varepsilon(\rho,\rho)^2 - \frac{1}{2} \OT_\varepsilon(\sigma,\sigma)^2
\end{equation}
remedies this issue and defines a distance for any value of $\varepsilon$. As shown in \cite{ramdas_wasserstein_2017, genevay_learning_2018}, $S_\varepsilon(\rho, \sigma)$ converges to $\OT(\rho, \sigma)$ as $\varepsilon$ goes to zero under mild assumptions on the densities. For $\varepsilon \rightarrow +\infty$, $S_\varepsilon$ approaches a \emph{maximum mean discrepancy} distance with the cost function as its kernel. We refer to \cite{feydy_geometric_2020}, Section 3.2.3 for details.

The solution to \cref{eq:sinkhorn_dvg} involves four potentials, two for the transport between $\rho$ and $\sigma$ stemming from the $\OT_\varepsilon(\rho, \sigma)$ term which we will denote $\psi_{\rho \rightarrow \sigma}, \psi_{\sigma \rightarrow \rho}$ and two potentials $\psi_{\sigma \rightarrow \sigma}, \psi_{\sigma \rightarrow \sigma}$ from the added terms:
\begin{equation}
	S_\varepsilon(\rho,\sigma)^2 = \int_\Omega ( \psi^\varepsilon_{\rho \rightarrow \sigma} - \psi^\varepsilon_{\rho \rightarrow \rho}) \rd \rho + \int_\Omega ( \psi^\varepsilon_{\sigma \rightarrow \rho} -  \psi^\varepsilon_{\sigma \rightarrow \sigma}) \rd \sigma.
\end{equation}
We can define the entropic Monge map analogously as 
\begin{equation}\label{eq:debiased_monge_map}
	T_{\rho \rightarrow \sigma}^{\varepsilon, \mathrm{debiased}} := \mathrm{id} - \nabla \psi_{\rho \rightarrow \sigma} + \nabla \psi_{\rho \rightarrow \rho}.
\end{equation}

\subsection{Entropic bias}\label{sec:entropic_bias}

Let $\pi, \varsigma \in \mathcal{P}(\Omega \times \Omega)$ and 
\begin{equation}
	\mathrm{KL}(\pi | \varsigma) := \int_{\Omega \times \Omega}  \left ( \log \left ( \frac{\rd \pi}{\rd \varsigma} \right ) \rd \pi - \rd \pi + \rd \varsigma  \right ).
\end{equation}
We observe that $\OT_\varepsilon(\rho,\sigma)^2 = \min_{\pi^\varepsilon \in \Pi(\rho,\sigma)} ( \int_{\Omega \times \Omega} c(x,y) \rd \pi^\varepsilon(x,y) + \varepsilon \mathrm{KL}(\pi^\varepsilon | \rho \otimes \sigma) )$ and $\mathrm{KL}(\pi | \varsigma_1) = \mathrm{KL}(\pi | \varsigma_2) + \mathrm{KL}(\varsigma_2 | \varsigma_1) \; \forall \varsigma_1, \varsigma_2 \in \mathcal{P}(\Omega \times \Omega)$. 

Let $\OT_\varepsilon^\varsigma(\rho,\sigma)^2 := \min_{\pi^\varepsilon \in \Pi(\rho,\sigma)} ( \int_{\Omega \times \Omega} c(x,y) \rd \pi^\varepsilon(x,y) + \varepsilon \mathrm{KL}(\pi^\varepsilon | \varsigma) )$. The choice $\varsigma = \rho \otimes \sigma$ in the definition of $\OT_\varepsilon$ (no superscript) is natural, as $\pi^\varepsilon$ is guaranteed to be absolutely continuous with respect to this product measure. However, many other choices (e.g. one that is constant on $\mathrm{supp}(\rho \otimes \sigma)$) are also possible. Choosing two different $\varsigma_1, \varsigma_2$ changes the value of $\OT_\varepsilon^{\varsigma_1}$ by the constant $ \mathrm{KL}(\varsigma_2 | \varsigma_1)$.

When computing OT barycenters, a different choice of $\varsigma$ leads to 
\begin{equation}
	\mathrm{OTBar}_\varepsilon^\varsigma( \{ \omega_i;  \rho_i \} _{i = 1}^n ) = \argmin_{\sigma \in \cP(\Omega)} \sum_{k=1}^n \omega_k \left \{ \OT_\varepsilon(\rho_k, \sigma)^2 + \varepsilon \mathrm{KL}(\varsigma | \rho_k \otimes \sigma) \right \}.
\end{equation}
Selecting the Lebesgue measure as $\varsigma$ leads to entropic smoothing, blurring the barycenter. This effect is discussed in \cite{cuturi_fast_2014, solomon_convolutional_2015}. In \cite{janati_debiased_2020}, it is shown that when all input measures are Gaussian $\rho_k = \mathcal{N}(\mu_k,  \mathrm{var})$, the entropic barycenter using the Lebesgue measure in $\mathrm{KL}$ will be a Gaussian $\mathcal{N}(\sum_k \omega_k \mu_k,  \mathrm{var} + \varepsilon)$. The choice of the product measure leads to $\mathcal{N}(\sum_k \omega_k \mu_k, \mathrm{var} - \varepsilon)$ for $\mathrm{var} > \varepsilon$ and $\delta_{\sum_k \omega_k \mu_k}$ otherwise.

This smoothing (resp. shrinking) bias can be seen as a feature of the method, as in \cite{rigollet_entropic_2018} it is shown that the entropic shrinking corresponds to a maximum-likelihood deconvolution technique. Alternatively, one can remove the effect of the choice in $\mathrm{KL}$ by replacing $\OT_\varepsilon$ with $S_\varepsilon$: It is straightforward to show that the value of $S_\varepsilon$ no longer depends on $\varsigma$, but only on $\mathrm{KL}(\pi^\varepsilon_{\rho, \sigma} | \pi^\varepsilon_{\rho, \rho})$ and $\mathrm{KL}(\pi^\varepsilon_{\rho, \sigma} | \pi^\varepsilon_{\sigma, \sigma})$. For the example of Gaussians $\rho_k = \mathcal{N}(\mu_k,  \mathrm{var})$, the $S_\varepsilon$ barycenter is a $\mathcal{N}(\sum_k \omega_k \mu_k,  \mathrm{var}) \; \forall 0<\varepsilon<+\infty$ (\cite{janati_debiased_2020}, Theorem 3).

In the following, replacing $\OT_\varepsilon$ by $S_\varepsilon$ will be referred to as \emph{de-biased} OT and de-biased barycenter will be used to refer to $ \argmin_{\sigma \in \cP(\Omega)} \sum_{k=1}^n \omega_k S_\varepsilon(\rho_k, \sigma)^2$.

\section{Transport mappings for reduced bases}\label{sec:method}

\begin{rem}
	In this section, all quantities related to OT ($\psi, \psi^c, T, \dots$) are denoted by their un-regularized form. The proposed method is applicable in this setting under sufficient regularity assumptions discussed in \cref{sec:regularity}. For computational feasibility, we use entropic regularization in practice. To apply the method in this case, one has to make the appropriate replacements, i.e. the c-transform becomes an application of the softmin, $\OT$ becomes $\OT_\varepsilon$ or $S_\varepsilon$, $\psi_{\rho \rightarrow \sigma}$ becomes either $\psi_{\rho \rightarrow \sigma}^\varepsilon$ or $\psi_{\rho \rightarrow \sigma}^\varepsilon - \psi_{\rho \rightarrow \rho}^\varepsilon$, and so forth. When specific assumptions or steps change depending on what regularisation is used, it is explicitly stated.
\end{rem}

As a motivating example, consider the pure advection equation.
\begin{examp}\label{ex:advection_eq}
	The PPDE $\partial_t \rho + \bar a \partial_x \rho = 0$ with initial condition $\rho_0 = \mathbbm{1}_{(-1,0]}$, a fixed advection $\bar a > 0$, and parameter $t \in \mathcal A = [0,T]$ on the domain $\Omega = [-1,1]$ is a prototypical example of a problem with very slow n-width decay of the solution manifold $\cM = \{ x \mapsto \rho_0(x - \bar a t) : t \in [0,1] \}$. The n-width in this case decays as slow as $\sim n^{1/2}$ \cite{ehrlacher_nonlinear_2020}.
\end{examp}
In contrast, if we set $\bar \rho = \rho_0$, the Monge embeddings of the set of solutions are of the form $\{ T_{\bar \rho \rightarrow \rho(t) }(y) = y + \bar a t  : t \in [0,1] \}$ with the corresponding potentials $\{ \psi^c_{t}(y) = - \bar a t y  : t \in [0,1] \}$. Clearly, this is a one-dimensional linear space. 

\subsection{Dimension reduction on the Monge embeddings}\label{sec:embedding_reduction}

Consider a PPDE problem with solution $u(\mu)$ and related densities $\rho(u)(\mu) =: \rho(\mu)$. In some cases, $\rho(u) = u$ is a possible choice, as in \cref{ex:advection_eq}. The requirements for $\rho$ are that it returns probability densities that coincide with the features of the solution that have to be registered. In \cite{iollo_mapping_2022}, a scalar testing function $\mathcal T$ is chosen to determine the distribution of features (which are sets of points $\{x \in \Omega :\mathcal T(x; u) > 0 \}$). Examples for $\mathcal T$ considered therein include $\Vert \nabla \times u \Vert, \Vert \nabla u \Vert$, the derivative of the Mach number of a flow, and a shock discontinuity indicator (equation (34) therein). For the proposed method, $\rho$ should return continuous densities supported on the entire domain, see \cref{sec:regularity}. In the numerical examples, we use $\rho(u) = \tfrac{\vert u(\mu) \vert^2}{\int \vert u(\mu) \vert^2}$ as one example.

The discretization of $u$ and $\rho$ also need not agree. It can be beneficial to discretise the latter on a regular tensor grid, see \cref{sec:computation}.

Given a set of snapshots $\{ u(\mu_i) \}_{i=1}^{n_s} \subset V_h$, compute $\{ \rho(u)(\mu_i) \}_{i=1}^{n_s}$ and denote by $\bar \rho$ a suited reference density, e.g. $\rho(\bar \mu)$ for a certain parameter value $\bar \mu$, or a weighted OT barycenter of $\{ \rho(\mu_i) \}_{i=1}^{n_s}$. Next, calculate the Monge embeddings $\{ T_{\bar \rho \rightarrow \rho(\mu_i)} \}_{i=1}^{n_s}$. We denote by $\psi^c_i$ the transport potential such that $T_{\bar \rho \rightarrow \rho(\mu_i)}(y) = y - \nabla \psi^c_i(y)$.
\begin{defi}\label{def:transport_modes}
	The \emph{transport modes} of a set of probability measures $\{ \rho(\mu_i) \}_{i=1}^{n_s}$ and a reference $\bar \rho$ are given by $y \mapsto \xi^c_j(y) = (\lambda^\psi_j)^{-1/2} \sum_{i=1}^{n_s}  (\mathrm{v}^\psi_j)_i \psi^c_i(y)$, where $\lambda^\psi_j$ and $\mathrm{v}^\psi_j$ are $j$th non-zero eigenvalue and eigenvector of the Monge embedding correlation matrix 
	\begin{equation}
		\bC^\psi := \{ \langle \nabla \psi^c_i, \nabla \psi^c_j \rangle_{L^2(\bar \rho)} \}_{1 \leq i,j \leq n_s}.
	\end{equation}
\end{defi}

Note the similarities with \cref{def:POD} for the POD modes.

\begin{examp}
	For snapshots of the pure advection equation at different times $t_i \in [0,T]$: $\{ \rho_0(x - \bar a t_i) \}_{i=1}^{n_s}$, we find $( \bC^\psi )_{ij} = \bar  a^2 t_i t_j$ with one non-zero eigenvalue $\lambda^\psi_1 = \bar  a^2 \sum_{i=1}^{n_s} t_i^2$ and eigenvector $(\mathrm{v}^\psi_1)_i = t_i ( \sum_{j=1}^{n_s} t_j^2 )^{-1/2}$. The corresponding transport mode is given by $\xi_1^c(y) = -y$. 
\end{examp}

If the eigenvalues of $\bC^\psi$ decay fast enough, all transport potentials $\psi^c(\mu)$ can be accurately (in the sense of a $\bar \rho$-weighted $L^2$ norm of their derivatives) approximated by a linear combination of the form $\psi^c(\mu) \approx \sum_{j=1}^m w_j(\mu) \xi^c_j$ where $m \ll n_s$. 

\begin{rem}
	We assume here that the map $\mu \rightarrow \psi^c(\mu)$ is very regular. This can be the case even if the mapping $\mu \mapsto u(\mu)$ is not (which can cause problems for reduced order modelling in the first place, see \cite{quarteroni_reduced_2016}, Section 5.3) and even if $\rho(u) = u$. Consider for example $\cM = \{ u(\mu) = \mathcal{N}(\mu, \mathrm{var}) : \mu \in \bR \}$. The map $\mu \mapsto u(\mu) \in L^1$ has a Lipschitz constant of $\tfrac{1}{\mathrm{var}}$, while the map $\mu \mapsto T_\mu = \mathrm{id} - \mu +\bar \mu$ is 1-Lipschitz for $\bar \rho = \rho(u)(\bar \mu)$.
\end{rem}
%\begin{rem}\label{ref:convexity_of_maps}
%	In order to guarantee that the approximate transport potential $\sum_{j=1}^m w_j(\mu) \xi^c_j$ is in fact a transport potential (i.e.: convex) it would be sufficient to take convex combinations of the snapshot potentials: $\psi^c(\mu) \approx \sum_{i=1}^{n_s} \omega_i(\mu) \psi^c_i$ with non-negative weights $\omega_{1,\dots,m}$ that sum to one. In that case, the function $\tfrac{1}{2} \vert y \vert^2 - \sum_{i=1}^{n_s} \omega_i(\mu) \psi^c_i$ is again a convex function. This is the case in the definition of the LOT barycenter \eqref{eq:LOT_barycenter}. Using a linear combination of transport modes, the resulting function $\tfrac{1}{2} \vert y \vert^2 - \sum_{j=1}^m w_j(\mu) \xi^c_j(y)$ will still be convex in practice, but this is not guaranteed by construction but a consequence of the quality of approximation through the modes $\{ \xi_i^c \}_{i=1}^m$.
%\end{rem}

\subsection{Reference reduced basis}

Evaluating 
\begin{equation}
	u(\mu_i) \circ \left ( \mathrm{id} - \nabla \sum_{j=1}^m w_j(\mu) \xi^c_j \right ) =: u(\mu_i) \circ \Phi_\mu^{-1} 
\end{equation}
applies the approximated transport mapping to the $i$th snapshot and yields elements of the mapped snapshot manifold $\Phi_\mu ( \cM )$. By construction, we expect this set to be more amicable to linear approximation. Returning to the simplest example, In the pure advection case, there is only one transport mode of the form $\xi_1^c(y) = - y$. The approximation of the snapshot potentials $\psi^c_i = - \bar a t_i y \in \mathrm{span} \{ \xi_1^c \} \; \forall i$ is exact and the mapped snapshot manifold consists of one single element.

More generally, we proceed by building a reduced basis in the reference space using the correlation matrix of transported snapshots 
\begin{equation}
	\bC^{\Phi_* u} := \{ \langle u(\mu_i) \circ \Phi_{\mu_i}^{-1},  u(\mu_j) \circ \Phi_{\mu_j}^{-1} \rangle_{V_h} \}_{1 \leq i,j \leq n_s}.
\end{equation}
Just as in the classical RB method described in \cref{sec:RB}, we obtain a set of reduced basis functions which we will denote by $\phi_{1,\dots,n_m}$. Now, any element of $\cM$ can be approximated via
\begin{equation}\label{eq:OTRB_approx}
	u_{\mathrm{trb}}(\mu) := \sum_{i=1}^{n_m} \tilde u(\mu)_i \, \phi_i \circ \left( \mathrm{id} - \nabla  \left [ \sum_{j=1}^m w_j(\mu) \xi^c_j  \right ]^c \right ) = \sum_{i=1}^{n_m} \tilde u(\mu)_i \, \phi_i \circ \Phi_\mu.
\end{equation}
\begin{rem}
	We use the properties of the c-transform to invert the mapping. This trick is possible since the gradients of Legendre transforms are inverses of each other (\cite{hiriart-urruty_fundamentals_2001}, Remark 0.1) and the c-transform and Legendre transform are related through \cref{eq:psi_phi_legendre}, as long as $\sum_{j=1}^m w_j(\mu) \xi^c_j$ is in fact a transport potential, i.e. a convex function (see \cref{sec:boundary_cond}).
\end{rem}
We conclude this subsection with an example for a one-dimensional PPDE that forms boundary layers from \cite{taddei_registration_2020}:

\begin{prop}\label{prop:boundary_layer}
	The solutions to the equation
	\begin{eqnarray}
		- \partial_{xx}^2 u_\mu + \mu^2 u_\mu = 0
	\end{eqnarray}
	on the domain $\Omega = (0,1)$ with boundary conditions $u_\mu(0) = 1$, $u_\mu(1) = 0$ and $\mu, \bar \mu \in [ \mu_{\mathrm{min}}, \mu_{\mathrm{max}}] =: \mathcal{A}$, $\mu_{\mathrm{max}} = \epsilon^{-2} \mu_{\mathrm{min}}$, $\mu_{\mathrm{min}} > 1$, $\epsilon \in (0,1)$ satisfy
	\begin{equation}\label{eq:boundary_layer_nm}
		\inf_{\xi^c_{1} \in \mathrm{span} \{ \psi^c_\mu : \mu \in \mathcal{A} \} } \; \sup_{\mu \in \mathcal A} \; \inf_{\substack{ w_1(\mu) \in \bR \\ \Phi^{-1}(y) = y - w_1(\mu) \partial_y \xi^c_1(y) \\ \Phi_\mu^{-1}: \Omega \rightarrow \Omega \text{ is a bijection}}} \Vert u_{\bar \mu} - u_\mu \circ \Phi_{\mu}^{-1} \Vert_{L^2(\Omega)} \leq e^{-\mu_{\mathrm{min}}}( 4 + \epsilon).
	\end{equation}
	where $\rho(u) = \tfrac{u}{\int u}$, $\bar \rho = \mathrm{OTBar}( \{ \rho_\mu : \mu \in \mathcal{A} \})$, and $\psi^c_\mu$ denotes the function such that $T_{\bar \rho \rightarrow \rho_{\mu}}(y) = y - \partial_y \psi^c_\mu(y)$ and $\int \psi^c_\mu = 0$.
\end{prop} 

In other words, we can show a bound on the Kolmogorov n-m-width (in the limit of tolerance $\rightarrow 0$, c.f. \cite{taddei_registration_2020}, section 3.2) of $\mathcal{M}$ for $n=m=1$. The proof of this proposition can be found in \cref{sec:n-m-supp}.

\subsection{Online phase}

To solve the PPDE problem for a new parameter value $\mu$, we determine the mapping $\Phi_\mu$, which is determined by the values of $w_{1,\dots,m}(\mu)$. Next, we plug the expression \eqref{eq:OTRB_approx} into the discretized PDE, and solve for $\tilde u_{1,\dots,n_m}(\mu)$. In this work, we learn the mapping $\mu \mapsto w_{1,\dots,m}(\mu)$ using a Gaussian process \cite{rasmussen_gaussian_2006} and the data from the snapshot set $\{ \mu_i, w_{1,\dots,m}(\mu_i) \}_{i=1}^{n_s}$. The functions could also be described by interpolation or any related method. We use the Gaussian process as it is computationally cheap also for high-dimensional data. 

The system of equations for $\tilde u_{1,\dots,n_m}(\mu)$ is then obtained by Galerkin projection using the reference reduced basis $\phi_{1,\dots,n_m}$. For example, a bilinear form corresponding to a Laplace operator reads
\begin{equation}\label{eq:laplace_reference}
	\int_\Omega \nabla (\phi_j \circ \Phi_\mu) \cdot  \nabla (\phi_j \circ \Phi_\mu) \, \rd x = \int_{\Phi_\mu(\Omega)} \nabla \phi_j \cdot [ D\Phi_\mu^{-1} ]^{-1}  [ D\Phi_\mu^{-1} ]^{-T} \nabla \phi_j \vert \det D\Phi_\mu^{-1} \vert \, \rd y
\end{equation}
where $\Phi_\mu(\Omega) = \Omega$ and $D\Phi_\mu^{-1} = \mathrm{Id} - \sum_{j=1}^m w_j(\mu) \, D^2 \xi^c_j$. 
\begin{rem}\label{rem:online_assembly}
	The drawback is that these forms have to be assembled for every new parameter value, and the computational cost for this depends on the dimension of the full-order problem. This is a challenge to any projection-based model order reduction method that utilizes a parameter-dependent mapping and requires hyper-reduction techniques to solve, see \cref{sec:hyperreduction}. %The are a number of strong hyper-reduction techniques available, for example the empirical quadrature \cite{yano_lp_2019}, which apply to this case (see also \cref{sec:numerical_examples} and \cref{sec:summary}).
\end{rem}
\begin{rem}
	If the parameters are time-dependent, or time is itself a parameter, the online phase will also feature an additional advection-like term:
	\begin{equation}
		\frac{\rd}{\rd t} u_{\mathrm{trb}}(\mu) =  \sum_{i=1}^{n_m} \left (  \frac{\rd \tilde u(\mu)_i}{\rd t}  \, \phi_i \circ \Phi_\mu + \tilde u(\mu)_i \, \frac{\rd \Phi_{\mu}}{\rd t} \cdot \left ( \nabla \phi_i \circ \Phi_\mu \right ) \right )
	\end{equation}
	In the reference domain, this requires the evaluation of $\tfrac{\rd \Phi_{\mu}}{\rd t} \circ  \Phi_{\mu}^{-1} = - [ D\Phi_\mu^{-1} ]^{-T} \, \tfrac{\rd \Phi^{-1}_{\mu}}{\rd t} $. Evaluating the latter expression is done using $\tfrac{\rd \Phi^{-1}_{\mu}}{\rd t} = - \sum_{j=1}^m \tfrac{\rd w_j(\mu)}{\rd t} \nabla \xi^c_j$.
\end{rem}
In summary, our proposed approach relies on snapshot remapping. The difference to other existing methods of this form is how the mappings are obtained. Our approach is data-driven and based on a POD of Monge embeddings. Other choices for parameter dependent mappings in the literature include problem-dependent parametrizations \cite{chetverushkin_model_2019}, polynomial expansion \cite{nair_transported_2019,welper_interpolation_2017}, and high-fidelity piece-wise polynomial mappings \cite{taddei_registration_2020} (see also \cref{sec:other_works} for a comparison of this approach and our method).

\subsection{Invertibility and boundary conditions of the mapping}\label{sec:boundary_cond}

For now, assume $\Omega = [0,1]^{d \in \{2,3\}}$, the unit square or cube. Proposition 2.3 in \cite{taddei_registration_2020} proves two sufficient conditions in order for a mapping of the form $\Phi^{-1}(y) = y - \sum_{j=1}^{m} w_j \nabla \xi^c_j$ to be a bijection in this case: Firstly, $\nabla \xi^c_j \cdot \hat e_i = 0 : 1 \leq i \leq d$ on all edges (faces), where $\hat e_{1,\dots,d}$ are normal vectors, and secondly $\det D\Phi^{-1} > 0$ in $\Omega \cup \partial \Omega$. In the case of more general mappings from $\Omega_2$ to $\Omega_1$, the first condition reads $\mathrm{dist}(\Phi^{-1}(y), \partial \Omega_1) = 0 \; \forall y \in \partial \Omega_2$ (Proposition 2.4 therein).

A natural question is if these conditions are met by the mappings defined in this section. We can give some answers in the non-regularized case. Suppose that $\psi^c$ is the Kantorovich potential for the transport from $\bar \rho \in \cP(\Omega_2)$ to $\rho \in \cP(\Omega_1)$, and both $\bar \rho$ and $\rho$ admit densities. Then, from \cref{thm:brenier}, we know that the transport map is given by $\nabla \varphi^* : y \mapsto y - \nabla \psi^c(y)$. Substituting this into the push-forward condition \eqref{eq:push-forward} gives, assuming enough regularity of $\varphi^*$ (see \cref{sec:regularity}),
\begin{equation}\label{eq:monge_ampere}
	\bar \rho = (\rho \circ \nabla \varphi^*) \, \det D^2 \varphi^*.
\end{equation}
We omitted the absolute value around the Jacobian determinant since we know $\varphi^*$ is convex. This PDE is known as the \emph{Monge-Amp\`{e}re equation}. It is subject to the transport condition $\nabla \varphi^*: \Omega_2 \rightarrow \Omega_1$. To fulfil this, it is enough to require that the boundary maps into the boundary, $\nabla \varphi^*: \partial \Omega_2 \rightarrow \partial \Omega_1$ \cite{urbas_mass_1998}, which equivalent to the first condition stated above.

As for the second condition, if both $\rho$ and $\bar \rho$ are bounded from above and below, we can write $\det D^2 \varphi^* = \tfrac{\bar \rho}{\rho \circ \nabla \varphi^*}$ and indeed one can show that $\varphi$ is strictly convex in this case, as the left hand side is strictly positive (see \cref{sec:regularity}). 

A sufficient condition to enforce positivity of the Jacobian determinant $\det D\Phi^{-1}$ would be that $\Phi^{-1}(y) = y - \sum_{j=1}^{n_s} \omega_j \nabla \psi^c_j$, $\omega_{1,\dots,n_s}$ are normalized, non-negative weights, and $\det (\mathrm{Id} - D^2 \psi^c_{j^*}) > 0, \omega_{j^*} > 0$ for at least one $j^*$. This is the setting of (L)OT barycenters. In this work, we opt to go for linear combinations over convex ones in order to make use of the POD compression on the tangent space at the cost of guaranteed bijectivity. A similar approach is taken in \cite{taddei_registration_2020}, see \cref{sec:registration_methods}. %It is still possible to recover the individual contributions of the snapshot potentials through the formula
%\begin{equation}\label{eq:bary_weights_from_lin_weights}
%	\sum_{j=1}^m w_j \xi^c_j = \sum_{i=1}^{n_s} \left ( \sum_{j=1}^m w_j ( \lambda^\psi_j)^{-1} (\mathrm{v}^\psi_j)_i \right ) \psi^c_i.
%\end{equation}

In the entropic case, the transport maps are also gradients of convex functions (\cref{rem:convexity_softmax}). However, the entropic Monge map does not exactly fulfil the OT boundary condition. This is also the case if the support of $\bar \rho, \rho$ does not coincide with $\Omega$. 

The boundary condition can be enforced in a post-processing step by $H^1$-projecting the $\psi^{c}$. To be precise, we solve the system
\begin{multline}
	\int_\Omega ( \kappa^2 \nabla \psi^c \cdot \nabla v + \psi^c v ) + \delta^{-1} \int_{\partial \Omega} ( \nabla \psi^c \cdot \hat n) v \\ = \int_\Omega ( \kappa^2 \nabla \psi^{c}_{\mathrm{pre-proj.}}(\mu_i) \cdot \nabla v + \psi^{c}_{\mathrm{pre-proj.}}(\mu_i) v ) \quad \forall v \in V_h
\end{multline}
for every $i = 1,\dots,n_s$. Here, $\psi^{c}_{\mathrm{pre-proj.}}$ denotes the output of the entropic OT calculations and $\psi^c$ is the projected potential used for the mapping. There are two parameters to set: $\delta$ is a penalty term to enforce the boundary condition and is set to $10^{-9}$ in our numerical experiments. The value of $\kappa$ determines the scale on which the function changes shape to fit the boundary condition. Since we expect the error introduced by the entropic smoothing to be of the scale $\sqrt \varepsilon$ and we want to limit the number of free parameters in our method, we set $\kappa^2 = {\varepsilon}^{-1}$.
\begin{rem} 
	If $\psi^{c}_{\mathrm{pre-proj.}}$ is very far from fulfilling the boundary conditions, this step can deform the potential to the point that $y \mapsto \tfrac{|y|^2}{2} -  \psi^c(y)$ is no longer convex and the mapping no longer invertible.
\end{rem}

\subsection{Regularity of the mapping}\label{sec:regularity}

The theory on the regularity of OT potentials is an involved topic and a proper treatment is beyond the scope of this work. The interested reader can refer to \cite{villani_topics_2016,loeper_regularity_2009,urbas_mass_1998,caffarelli_regularity_1992}.

In the un-regularized case, the transport potential $\psi$ from $\rho$ to $\sigma$ introduced in Brenier's \cref{thm:brenier} is differentiable $\rho$-almost everywhere as long as $\rho$ is absolutely continuous with respect to the Lebesgue measure. The brief explanation of this fact is that through the c-transform, $\psi$ shares the modulus of continuity with the cost function (\cite{santambrogio_optimal_2015}, Section 1.2). Clearly, in our method, we implicitly assume higher regularity of $\psi$. This can be argued for with the following result, which is a special case of the regularity theory developed by Caffarelli as presented in \cite{villani_topics_2016}, Theorem 4.14: When $\Omega$ is open, bounded, uniformly convex, and of class $C^{2,\alpha}$, $\rho$ and $\sigma$ are H\"older continuous functions $\in C^{0,\alpha}(\bar \Omega)$, bounded above and below by positive constants, then the solution $\varphi$ to the Monge-Amp\`{e}re \cref{eq:monge_ampere} is of class $C^{2,\alpha}(\Omega) \cap C^{1,\alpha}(\bar \Omega)$. In particular, through the classical bootstrapping arguments from elliptic regularity, $\varphi$ is $C^\infty$ in the interior if $\rho$ and $\sigma$ are.

Even in this case, the constants appearing in these estimates can be prohibitively large. Consider the example from \cref{prop:boundary_layer}. As shown in \cref{sec:n-m-supp}, the transport map reads $T_{\rho_{\bar \mu} \rightarrow \rho_{\mu}}(y) = 1 - \frac{1}{\mu} \sinh^{-1} \left ( \frac{\sinh \mu}{\sinh \bar \mu} \sinh(\bar \mu(1-y)) \right )$ and therefore $\partial_yT_{\rho_{\bar \mu} \rightarrow \rho_{\mu}}(1) \approx \tfrac{\bar \mu}{\mu} e^{\mu - \bar \mu}$, which can take extreme values for typical $\mu$ and $\bar \mu$. 

For the regularized case, the transport potentials are smooth as they are defined by Gaussian convolutions in \cref{eq:schroedinger_eq} and we recall that $\Vert \nabla^k \psi \Vert_\infty = \mathcal{O}(1 + \varepsilon^{1-k})$ \cite{genevay_sample_2019}. However, the same caveats apply. We ultimately require our regularized transport map to resemble the un-regularized one. Increasing $\varepsilon$ too much will degrade the quality of the mapping.
The projection step as defined in \cref{sec:boundary_cond} will preserve the regularity of the potentials (with a constant depending on $\kappa$) as it is an elliptic problem.

\begin{figure}[h!]
	\centering
	\includegraphics[width=0.45\textwidth]{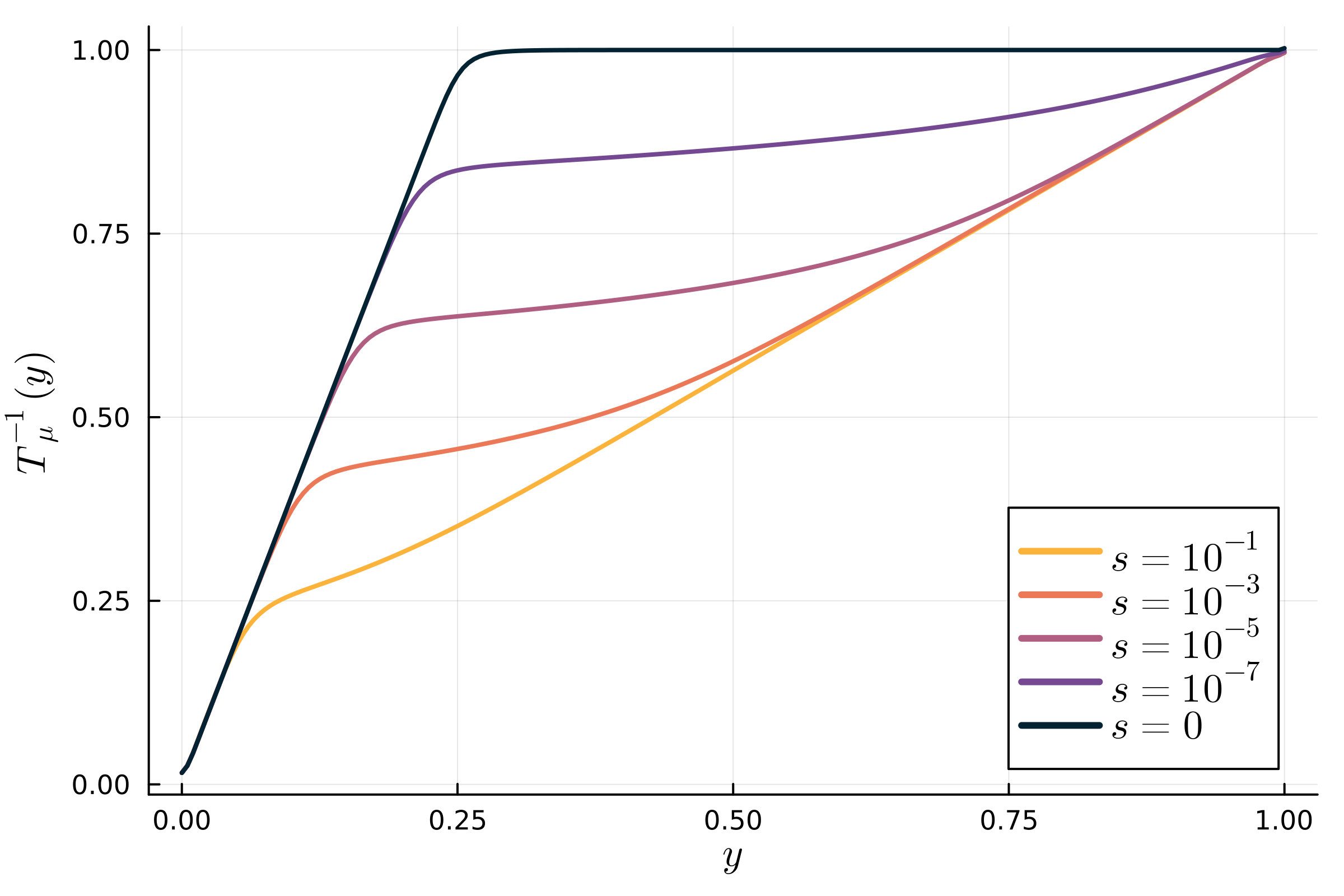}
	\includegraphics[width=0.45\textwidth]{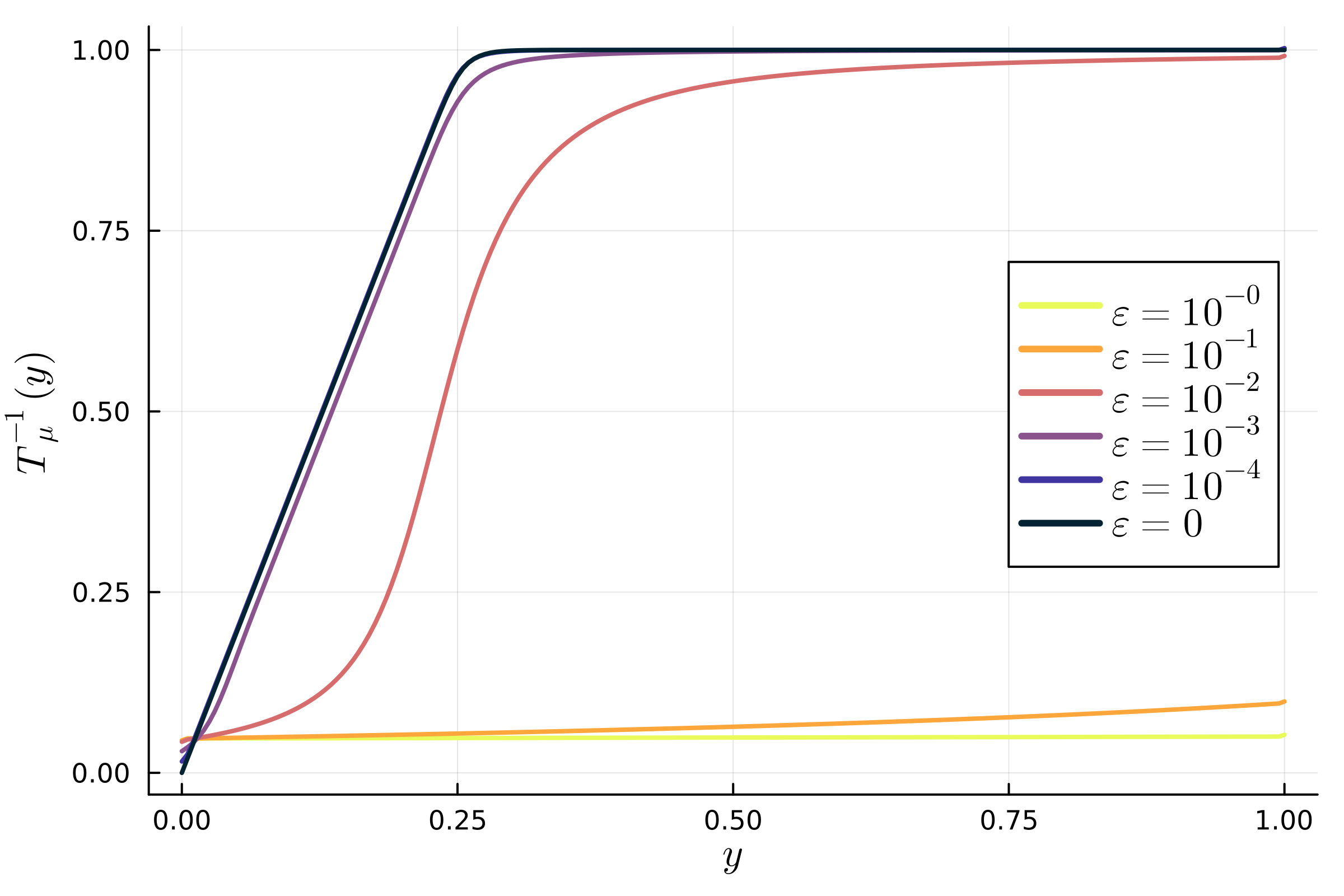}
	\caption{
		Transport maps for the problem from \cref{prop:boundary_layer} with $\epsilon^2 = 10^{-1}$ and $\mu_{\mathrm{min}} = 20$. Left: transport map $T_{\mu_{\mathrm{min}}}^{-1}$ for different values of $s$ with $\varepsilon$ constant at $10^{-4}$. Right: different values of $\varepsilon$ with $s \equiv 0$. We see that the parameter $s$ that sets a lower bound to the densities can control the derivatives of the mapping in this example. Note that modification of $T$ beyond the point $x^\star \approx \tfrac{\mu}{\bar \mu} \approx 0.26$ does not impact the approximation result in \cref{prop:boundary_layer} (\cref{sec:n-m-supp}) so that the choice $s = 10^{-7}$ provides the same error bounds while keeping the derivatives of $T$ and $T^{-1}$ under control.}
	\label{fig:T_reg}
\end{figure}

In practice, it seems to be effective to increase the lower bound on the densities. For example, when using $\rho(u) = (1-s) \tfrac{u}{\int u} + s$ with $s > 0$ in \cref{prop:boundary_layer}, we see that the derivative of the transport map is controlled (\cref{fig:T_reg}). The entropic smoothing for OT barycenters discussed in \cref{sec:entropic_bias} can provide this effect.

Lastly, we have to remark that we cannot directly apply these regularity results to the numerical examples of \cref{sec:numerical_examples}: our domain is only of class $C^0$ and the potentials $\psi$ are approximated using $H^1$ conforming finite element functions that are only piece-wise $C^1$.

\subsection{Hyper-reduction}\label{sec:hyperreduction}

When the online phase can be made fully independent of the size $N$ of the high-fidelity problem, its computational cost can be reduced dramatically. For now, the assembly of linear and bilinear forms has to be done online and depends on $N$ (\cref{rem:online_assembly}). We can remedy this shortcoming using the \textit{Empirical Interpolation Method} (EIM). In particular, we utilize a version of EIM based on a POD of the parameter-dependent forms \cite{chaturantabut_nonlinear_2010, taddei_registration_2020}. We briefly recall the method using the example of the mapped Laplace operator from \cref{eq:laplace_reference}.

Based on the data from the training set, a collection of snapshots
\begin{equation}
	\{ [ D\Phi_{\mu_i}^{-1} ]^{-1}  [ D\Phi_{\mu_i}^{-1} ]^{-T} \vert \det D\Phi_{\mu_i}^{-1} \vert \}_{i = 1}^{n_s} =: \{ K_{\mu_i}  \}_{i = 1}^{n_s}
\end{equation}
is used to obtain a POD basis from the correlation matrix $\bC^K : \bC^K_{ij} = \int_\Omega \mathrm{tr} ( K_{\mu_i}^T K_{\mu_j} ) \, \rd y, 1 \leq i,j \leq n_s$. Using an energy criterion $\tau_{\mathrm{eim}}$, the eigenvectors $\Xi_q : 1 \leq q \leq Q$ corresponding to the $Q$ largest eigenvalues are selected to span an approximation space. Coefficients $\theta_q(\mu)$ and functions $X_q : 1 \leq q \leq Q$ are determined such that $K_\mu \approx \sum_{q = 1}^{Q} \theta_q(\mu) X_q$ for all $\mu$. The way the interpolation points and functions are selected guarantees that the matrix $B \in \bR^{Q \times Q} : B_{q' q} = X_{q'}(y_q^{\mathrm{eim}})$ is lower-triangular with unit diagonal, so the interpolation problem is well-defined and quickly (i.e. in $\mathcal{O}(Q^2)$ time) solved. The procedure is outlined in \cref{sec:eim_appendix}, \cref{alg:eim}.

Online, $K_\mu$ is evaluated at select points $\{ y_q^{\mathrm{eim}} \}_{q=1}^Q$ and the interpolation problem $K_\mu(y_q^{\mathrm{eim}}) = \sum_{q' = 1}^{Q} \theta_{q'}(\mu) X_{q'}(y_q^{\mathrm{eim}}) : 1 \leq q \leq Q$ is solved to obtain $\{ \theta_q(\mu) \}_{q=1}^Q$. The full form $\int_{\Omega} \nabla \phi_j \cdot K_\mu \nabla \phi_j \, \rd y$ is approximated using $\sum_{q=1}^Q \theta_q(\mu) \int_{\Omega} \nabla \phi_j \cdot X_{q,ij} \nabla \phi_j \, \rd y$ where the integral defines a $Q \times n_m \times n_m$ tensor that can be pre-computed offline. As a result, no integration in the high-fidelity space has to be done online. The online cost of the EIM procedure consists of $Q$ evaluations of $K_\mu$, the $\mathcal{O}(Q^2)$ interpolation problem, and a $\mathcal{O}(Qn_m^2)$ tensor contraction. Importantly, it does not depend on $N$. Oversampling techniques which have been shown to improve the stability of the empirical interpolation method in the presence of noisy data \cite{peherstorfer_stability_2020} are not used in this work.
\section{Comparison to other works}\label{sec:other_works}

There have been a number of works that link model order reduction with OT techniques. In this section, we will briefly present some of them and discuss how they relate to the present work.

\subsection{Optimal transport barycenter coordinates}

Recall the definition of an OT barycenter \eqref{eq:OT_barycenter}:
\begin{equation}
	\mathrm{OTBar}( \{ \omega_i; \rho_i \}_{i=1}^n) := \argmin_{\sigma \in \cP(\Omega)} \sum_{i=1}^n \omega_i \, \OT(\rho_i, \sigma)^2. \nonumber
\end{equation}
Given a suitable set of probability densities $\{\rho_1, \dots, \rho_n\}$ denoted \emph{atoms} in the spirit of dictionary learning, for any $\rho \in \cP(\Omega)$, one can define optimal weights as
\begin{equation}\label{eq:optimal_weights}
	\{ w^{\mathrm{opt}}_i(\rho) \}_{i=1}^n := \argmin_{\substack{\omega_1, \dots, \omega_n \geq 0 \\ \sum_{i=1}^n \omega_i = 1}} \mathrm{Loss} \left ( \rho,  \mathrm{OTBar}(\omega_1, \dots, \omega_n; \rho_1, \dots, \rho_n) \right ),
\end{equation}
given a suitable loss function $\cP(\Omega) \times  \cP(\Omega) \rightarrow \bR$. This method yields a parametrization of elements of $\cP(\Omega)$ through a small number of $n$ non-negative, normalized weights.

 In \cite{bonneel_wasserstein_2016}, the authors apply this procedure to several application cases from computer graphics to medical imaging. The optimal weights are called \emph{barycenter coordinates} and different loss functions are used. The optimization problem \eqref{eq:optimal_weights} is solved using a standard L-BFGS quasi-Newton solver. 
 
 A different approach is proposed in \cite{werenski_measure_2022}. Therein, the authors use the fact that at the minimum defining the OT barycenter \ref{eq:OT_barycenter}, the Fr\'{e}chet derivative of $\sigma \mapsto  \sum_{i=1}^n \omega_i \OT(\rho_i, \sigma)^2$ has to be zero. Under mild assumptions outlined in their work, it holds that $\nabla_\sigma \sum_{i=1}^n \omega_i \OT(\rho_i, \sigma)^2 = - \sum_{i=1}^n \omega_i ( T_{\sigma \rightarrow \rho_i} - \mathrm{id} ) =  \sum_{i=1}^n \omega_i \nabla \psi_{\sigma \rightarrow \rho_i}$. One obtains a quadratic optimization problem on the simplex of normalized weights featuring a correlation matrix as introduced in \cref{def:transport_modes}:
 \begin{equation}
 	\left \Vert \nabla_\rho \sum_{i=1}^n \omega_i \OT(\rho_i, \rho)^2 \right \Vert^2_{ L^2(\rho) } = \sum_{i,j = 1}^n \omega_i \omega_j \langle \nabla \psi_{\rho \rightarrow \rho_i}, \nabla \psi_{\rho \rightarrow \rho_j} \rangle_{L^2(\rho)} \overset{!}{=} 0,
 \end{equation}
 subject to $\omega_1, \dots, \omega_n \geq 0, \sum_{i=1}^n \omega_i = 1$. 
 
 In the later works \cite{mueller_geometric_2022,schmitz_wasserstein_2018}, the authors do not work with a given set of atoms, but instead determine them by optimization: given a training set $\{ \rho_i \}_{i=1}^{n_s}$, the optimization problem reads
\begin{equation}\label{eq:schmitz_barycenter_coordinates}
	\min_{\substack{\omega_j(\rho_i) \geq 0 \; \forall i,j \\ \sum_{i=j}^n \omega_j(\rho_i) = 1 \; \forall i \\ \sigma_1, \dots, \sigma_n \in \cP(\Omega)}} \sum_{i=1}^{n_s} \mathrm{Loss} \left ( \rho_i,  \mathrm{OTBar}(\{ \omega_j(\rho(i)); \sigma_j \}_{j=1}^n ) \right ).
\end{equation}

Equations \eqref{eq:optimal_weights} and especially \eqref{eq:schmitz_barycenter_coordinates} are complicated multilevel optimization problems, which however can be tackled using a tailored warm-start technique as outlined in Section 4.2 of \cite{schmitz_wasserstein_2018}. The gradients involved can either be computed from closed formulas provided therein, or using automatic differentiation techniques. However, even with those adaptations, solving for truly optimal barycentric weights online remains computationally restrictive.

Another option to construct the set of atoms $\{ \sigma_i \}_{i=1}^n$ is through greedy algorithms, see \cite{battisti_waserstein_2022,ehrlacher_nonlinear_2020}. Promoting sparsity in the weight vectors $\{ \omega_i \}_{i=1}^n$ is beneficial to memory footprint and run-time cost on the one hand, and mitigates issues of possible redundancy in the set of atoms on the other hand \cite{do_approximation_2023,mueller_geometric_2022}.

Most of the works mentioned utilize entropic regularisation to make the OT computations feasible, except for \cite{ehrlacher_nonlinear_2020,battisti_waserstein_2022}, where the authors apply barycentric approximations to one-dimensional PPDE problems. In one spatial dimension, there is a closed formula available for the transport from $\rho$ to $\sigma$, given by their cumulative distribution functions (cdf):
\begin{equation}
	T_{\rho \rightarrow \sigma} = \mathrm{cdf}(\sigma)^{[-1]} \circ \mathrm{cdf}(\rho).
\end{equation}
The superscript $^{[-1]}$ denotes the pseudo-inverse as defined in chapter 2 of \cite{santambrogio_optimal_2015}. Note that for the case $\Omega = [0,1]$ and $\bar \rho = \mathbbm{1}_{[0,1]}$, the inverse cdf operation coincides with the Monge embedding: $T_{\bar \rho \rightarrow \sigma} = \mathrm{cdf}(\sigma)^{[-1]}$. The method of \emph{tangent principal component analysis} (tPCA) \cite{ehrlacher_nonlinear_2020} on the set of $\{ \mathrm{cdf}(\rho(\mu_i))^{[-1]} \}_{i=1}^{n_s}$, which proved numerically unstable, is a special case of performing a POD on the Monge embeddings.

Convex displacement interpolation (CDI), introduced in \cite{iollo_mapping_2022}, is based on a linear approximation of the displacement interpolation between two measures. The CDI method is non-intrusive and similar in spirit to the barycentric coordinates. In order to solve the OT problems needed to build the displacement interpolations, the authors rely on the closed form available for the OT between multivariate Gaussian densities \cite{peyre_computational_2019}. See also \cite{rim_displacement_2018}, where displacement interpolation is used to interpolate between solutions of PPDEs at different parameter values.

OT barycentric coordinates have proven effective at parametrizing sets of probability densities in numerous applications. One drawback of these methods, however, is that it is very hard to go beyond interpolation-based methods in the online phase. The reconstruction $\omega_1(\mu), \dots, \omega_n(\mu) \mapsto \mathrm{OTBar}(\mu)$ is costly and non-linear, complicating the online evaluation of a PDE residual such as $\Vert \mathcal{L}(\mathrm{OTBar}(\mu); \mu) \Vert_{L^2}$.

Methods based on OT barycenters are inherently limited when it comes to extrapolation tasks (see \cite{battisti_waserstein_2022}, section 5.3). We avoid this limitation in our work by moving to the tangent space of $\mathcal{P}(\Omega)$ and replacing the convex combinations with linear ones.

\subsection{Registration methods}\label{sec:registration_methods}

As discussed in \cref{sec:RB}, the principle of recasting PPDE problems using parametric mappings $\Phi_\mu$ in order to obtain a mapped solution manifold $\Phi_\mu ( \cM ) = \{ u(\mu) \circ \Phi^{-1}_\mu : u(\mu) \in \cM \}$ that is better suited for linear reduction techniques, has been studied extensively. As an example, we will discuss the method presented in \cite{taddei_registration_2020,taddei_space-time_2021}. Therein, mappings of the form
\begin{equation}
	y \mapsto y + \sum_{j=1}^{m^{\mathrm{hf}}} w^{\mathrm{hf}}(\mu)_j \chi_j^{\mathrm{hf}}(y)
\end{equation}
are proposed. As in \cref{sec:method}, $w^{\mathrm{hf}}(\mu)_{1,\dots,m^{\mathrm{hf}}} \in \bR$ are parameter-dependent coefficients, while $\chi^{\mathrm{hf}}_{1,\dots,m^{\mathrm{hf}}}$ are elements of a general approximation space such as Legendre polynomials or Fourier expansions. The mappings are constructed such that
\begin{equation}
	\left \Vert \, u(\mu) \circ \left ( y + \sum_{j=1}^{m^{\mathrm{hf}}} w^{\mathrm{hf}}(\mu)_j \chi_j^{\mathrm{hf}}(y) \right ) - \bar u \, \right \Vert
\end{equation}
is small, given a reference $\bar u$. Besides this \emph{proximity measure}, the optimization penalizes the $H^2$ seminorm of the mappings and enforces certain constraints to keep the Jacobian of the mappings strictly positive. To guarantee a sufficiently rich set of mappings to optimize over, $m^{\mathrm{hf}}$ has to be rather large, which can be restrictive when evaluating mappings in the online phase. Consequently, the authors also opt for a POD approach, reducing the number of mapping terms to $m$ based on an eigenvalue decomposition of the matrix with elements $\bC^w= \{ w^{\mathrm{hf}}(\mu)_i \cdot w^{\mathrm{hf}}(\mu)_j \}_{1 \leq i,j \leq m^{\mathrm{hf}}}$.
%\begin{equation}\label{eq:mapping_taddei}
%	\sum_{j=1}^{m^{\mathrm{hf}}} w^{\mathrm{hf}}(\mu)_j \chi_j^{\mathrm{hf}} \approx \sum_{j=1}^m w(\mu)_j \chi_j.
%\end{equation}
This method is similar to the one we propose in the present work. Note that also in this case, it cannot be guaranteed that this approximate mapping is invertible, even if the high fidelity one is. However, the method proved stable in the numerical test cases considered.

%Comparing this approach to the one from \cref{sec:method}, the main difference lies in the calculation of the high-fidelity mapping. The appeal of using the transport maps from LOT is that there is only one hyperparameter $\varepsilon$ which controls both the proximity measure (as $\varepsilon \rightarrow 0$, $T^\varepsilon(\mu)_\sharp \bar \rho \rightarrow \rho(\mu)$ in the $L^2(\bar \rho)$ norm \cite{pooladian_entropic_2022}) and the derivatives of the mapping ($\Vert \nabla^{(k)} T^\varepsilon \Vert_\infty = \mathcal{O}(1 + \varepsilon^{-k})$ \cite{genevay_sample_2019}). 

In \cite{taddei_optimization-based_2022}, figure 5, the authors show the performance of a registration method based on a non-linear optimization problem when using a general polynomial space $\{\chi^\mathrm{hf}_j\}_{j=1}^{m^\mathrm{hf}}$ and a space comprised only of gradient functions, i.e. $\chi^\mathrm{hf}_j = \nabla \psi^\mathrm{hf}_j : 1 \leq j \leq m^{\mathrm{hf}}$ of similar dimension. Methods that rely on transport maps to approximate $\Phi$ are by construction confined to the latter space. In this example, it is shown that while the performance is very similar for small and medium deformations, the optimization method using gradient functions deteriorates for large deformations, leading to $\sim 10 \times$ higher iteration numbers and even inverted elements (i.e. non-bijective $\Phi$). The method proposed in this work cannot be directly applied to this test case, as it is a point registration problem and un-regularized OT gives no guarantees of invertible mappings when the transported measures give mass to small sets.

In \cite{iollo_advection_2014}, the authors propose a snapshot registration method based on \emph{advection modes}, which are obtained using the matrix $\{ \OT(\rho(\mu_i),\rho(\mu_j)) \}_{1 \leq i,j \leq n_s}$ and applying techniques of Euclidean distance matrix methods \cite{dokmanic_euclidean_2015}. The mapping of the snapshots is done in the sense of push-forward measures, that is
\begin{equation}
	\rho(\mu) \approx \left ( \bar \rho \circ \Phi_\mu + r(\mu) \circ \Phi_\mu  \right ) \det D\Phi_\mu.
\end{equation}
Here, $r(\mu)$ denotes a residual term that needs to be determined in the online phase. Especially in the case where $u$ itself is used to construct the transport mappings, this would be a natural choice for our method as well, since the transport mappings are constructed to fulfil this push-forward relation. Up to numerical inaccuracies, we would expect $n_m = 1$.

We chose the form \cref{eq:OTRB_approx}, which corresponds to the push-forward of a function, for two reasons: First, the relation $\bar \rho = \rho(\mu_i) \circ \left ( \mathrm{id} - \nabla \psi^c(\mu_i) \right ) \det \left ( \mathrm{Id} - D^2 \psi^c(\mu_i) \right ) \; \forall i$ only holds if the transport mappings are constructed directly from $u$ itself and not another derived quantity $\rho(u)$. 
%Second, mappings of the form $\phi \, \circ \, \Phi_\mu$ are used in standard finite element codes to go from a reference domain to the physical domain. These routines can then be utilized in the implementation. 
Second, substituting the push-forward relation for a density into a PDE residual requires evaluating derivatives of $\det \left ( \mathrm{Id} - D^2 \psi^c \right )$, which requires even more regularity of $\psi^c$, and higher order basis functions in the discretization.

To summarize, our work aims to hit a compromise between interpolating methods based on optimal transport and registration approaches. By using linear OT and transforming our snapshots as functions, not densities, we accept the need for more reference basis functions $\phi_{1,\dots,n_m}$ in return for a lower-order residual that has to be evaluated in the online phase.

\section{Numerical experiments}\label{sec:numerical_examples}

We will demonstrate the proposed method and the impact of some of the hyperparameters on two test cases. For the Finite Element calculations, we rely on the \texttt{Gridap.jl} library\footnote{https://github.com/gridap/Gridap.jl} \cite{verdugo_software_2022,badia_gridap_2020}, while Gaussian processes are calculated using \texttt{GaussianProcesses.jl}\footnote{https://github.com/STOR-i/GaussianProcesses.jl}. The computational OT routines used are available in the package \texttt{WassersteinDictionaries.jl}\footnote{https://github.com/JuliaRCM/WassersteinDictionaries.jl}.

\subsection{Poisson's equation with moving source}

Let $\Omega = [0,1]^2$, discretized by a $64\times64$ grid of quadrilateral cells. For $V_h$, we chose $H^1_0$-conforming Legendre basis functions of order $p=3$, which leads to $N=36481$ degrees of freedom. The grid size is denoted with $h$. The equation we solve is
\begin{equation}
	\Delta u(x; \mu) = f(x; \mu) : x \in \Omega, u(x; \mu) = 0 : x \in \partial \Omega,
\end{equation}
where $f$ is a narrow Gaussian with variance $\mathrm{var} = 10^{-3}$ and mean $( \tfrac{1}{2}, \tfrac{1}{2} ) + \mu$, $\mu \in ( -\tfrac{7}{20}, \tfrac{7}{20} )^2 = \mathcal{A}$.

To construct the mappings and reduced basis we draw  $n_s = 100$ samples of $\mu$ uniformly from $\mathcal{A}$. The solutions to this equation are not probability densities, so we use $\rho(u)(\mu) = \tfrac{u(\mu)^2}{ \int u(\mu)^2}$ and $\rho(u)(\mu) = f(\mu)$ to calculate the transport mappings. These computations, which rely on collocation, are performed on a finer $192 \times 192$ grid of quadrilaterals. The iterative OT solver is set to stop when the $l_1$ error of the marginal constraint reaches $10^{-3}$ or at the maximum number of iterations of $\lceil \tfrac{10}{\varepsilon} \rceil$.
%As reference density $\bar \rho$ we chose $\rho(u)(\bar \mu)$ where $\bar \mu$ is the mean value of $\{ \mu_i \}_{i=1}^{n_s}$.
As reference density $\bar \rho$ we chose the OT barycenter of the training snapshot densities. Note that it is not crucial for the proposed method to use the OT barycenter. Any reference density, even $\bar \rho \equiv |\Omega|^{-1}$, can be used. However, when $\bar \rho$ is an average of the data $\{ \rho(\mu_i) \}_{i=1}^{n_s}$ in a meaningful sense (which the OT barycenter is in these cases), we expect the transport maps to have a simpler structure and require fewer modes $m$ to approximate.
We employ an annealing strategy as described in \cite{feydy_geometric_2020}, section 3.3.3, initializing the regularization parameter to one and then halving it at every iteration until we reach $\varepsilon$. To invert the mapping at the last step, we use the c-transform with $\varepsilon_{\mathrm{fine}} = \tfrac{h^2}{10}$. 

Given a threshold $\tau(\mathcal{E})$, $m$ is chosen such that $1 - \mathcal{E}(m; \lambda) < \tau$, where
\begin{equation}\label{eq:eigenvalue_energy}
	\mathcal{E}(m; \lambda) := \frac{\sum_{i=1}^{m} \lambda_i}{ \sum_{j = 1}^{\mathrm{rank} \, \bC} \lambda_j }.
\end{equation}
We set $\tau_{\mathrm{eim}} = \tfrac{\tau}{10}$ to account for the difficulty in approximating the moving source term. 

The Gaussian process we use is taken as-is from the reference package. In particular, we select a zero mean function, a squared exponential kernel with characteristic length and standard deviation set to one (the default settings). The log standard deviation of observation noise is set to $-6$. These parameters have not been optimized.

\subsubsection{The case $\rho(u)(\mu) \propto u(\mu)^2$}

Since $\rho(\mu)$ is positive on all of $\Omega$ in this case, we use de-biased calculations, replacing $\OT_\varepsilon$ with $S_\varepsilon$ when needed.The transport maps are now given by the de-biased potentials as defined in \cref{eq:debiased_monge_map}. We see that de-biasing improves the performance of the method, both by increasing the accuracy and by reducing the number of approximation functions $n_m$ and $Q$, in \cref{tab:errors_debias}. 

\begin{figure}[h!]
	\centering
	\includegraphics[width=0.45\textwidth]{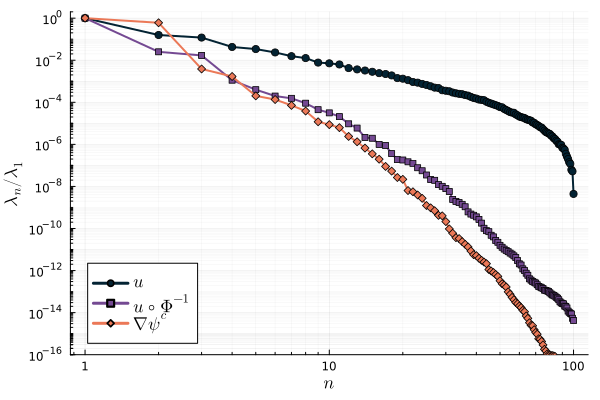}
	\includegraphics[width=0.45\textwidth]{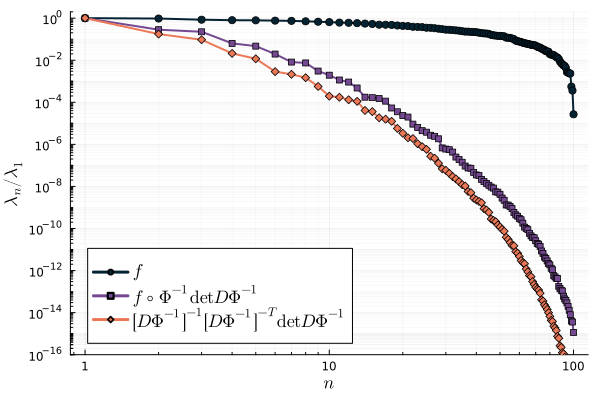}
	\caption{Left: Eigenvalues of the correlation matrices $\bC^u$, $\bC^{\Phi_\ast u}$, and $\bC^\psi$. Right: Eigenvalues of the correlation matrices  $\bC^f$, $\bC^{K}$, and $\bC^{\Phi_\ast f}$ used in the hyper-reduction. $\varepsilon = 10^{-2}, \tau = 10^{-4}, \rho(u)(\mu) \propto u(\mu)^2,$ and de-biasing are used. }
	\label{fig:evds_1}
\end{figure}

\Cref{fig:evds_1} displays the eigenvalue decay for the correlation matrices of snapshots $\bC^u$, transported snapshots $\bC^{\Phi_\ast u}$, and Monge embeddings $\bC^\psi$. As expected, the eigenvalues of the mapped snapshots are indicative of a much faster n-width decay of $\Phi_\mu ( \cM )$ compared to that of $\cM$.

\begin{figure}[h!]
	\centering
	\includegraphics[width=0.23\textwidth]{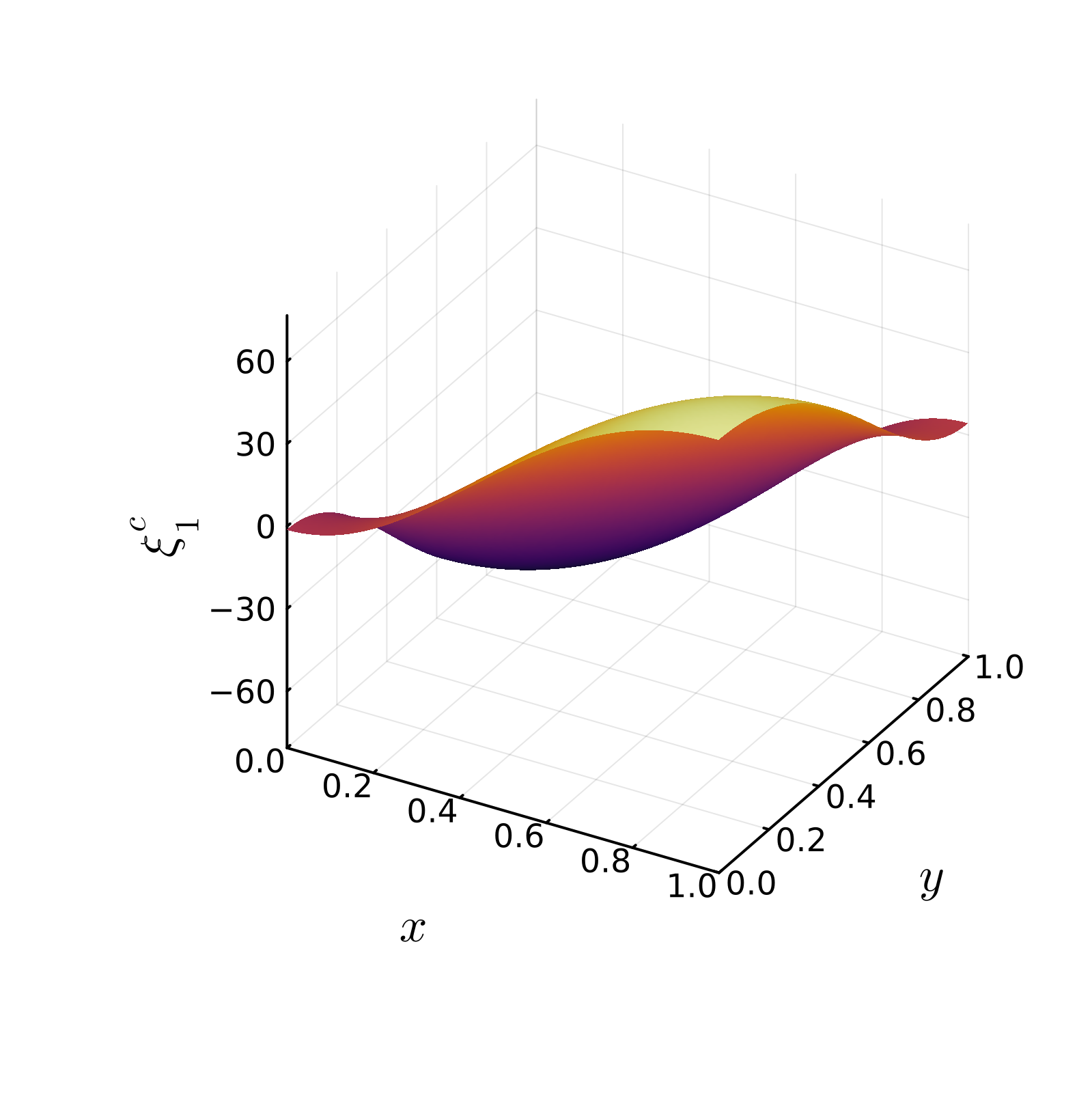}
	\includegraphics[width=0.23\textwidth]{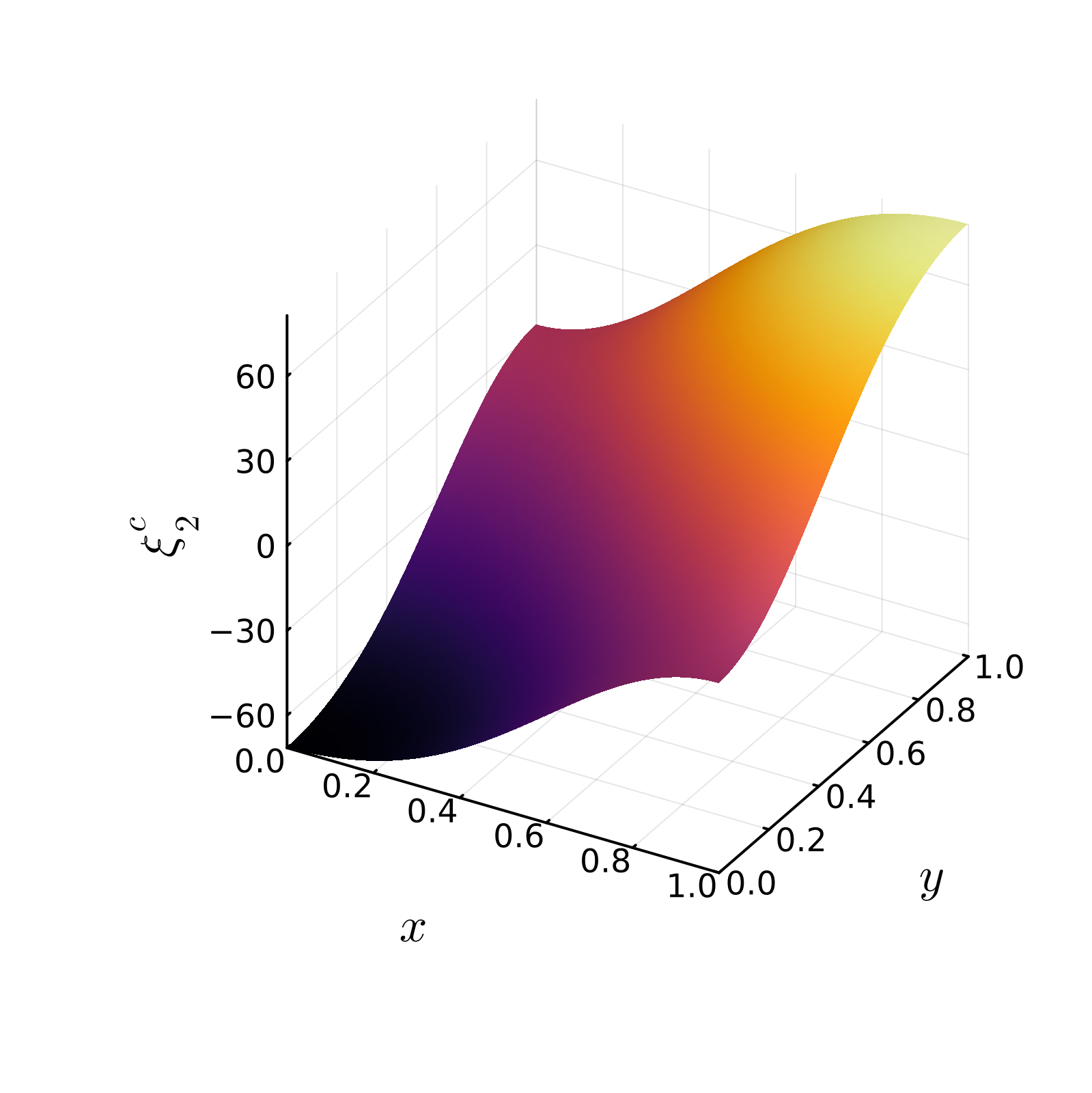}
	\includegraphics[width=0.23\textwidth]{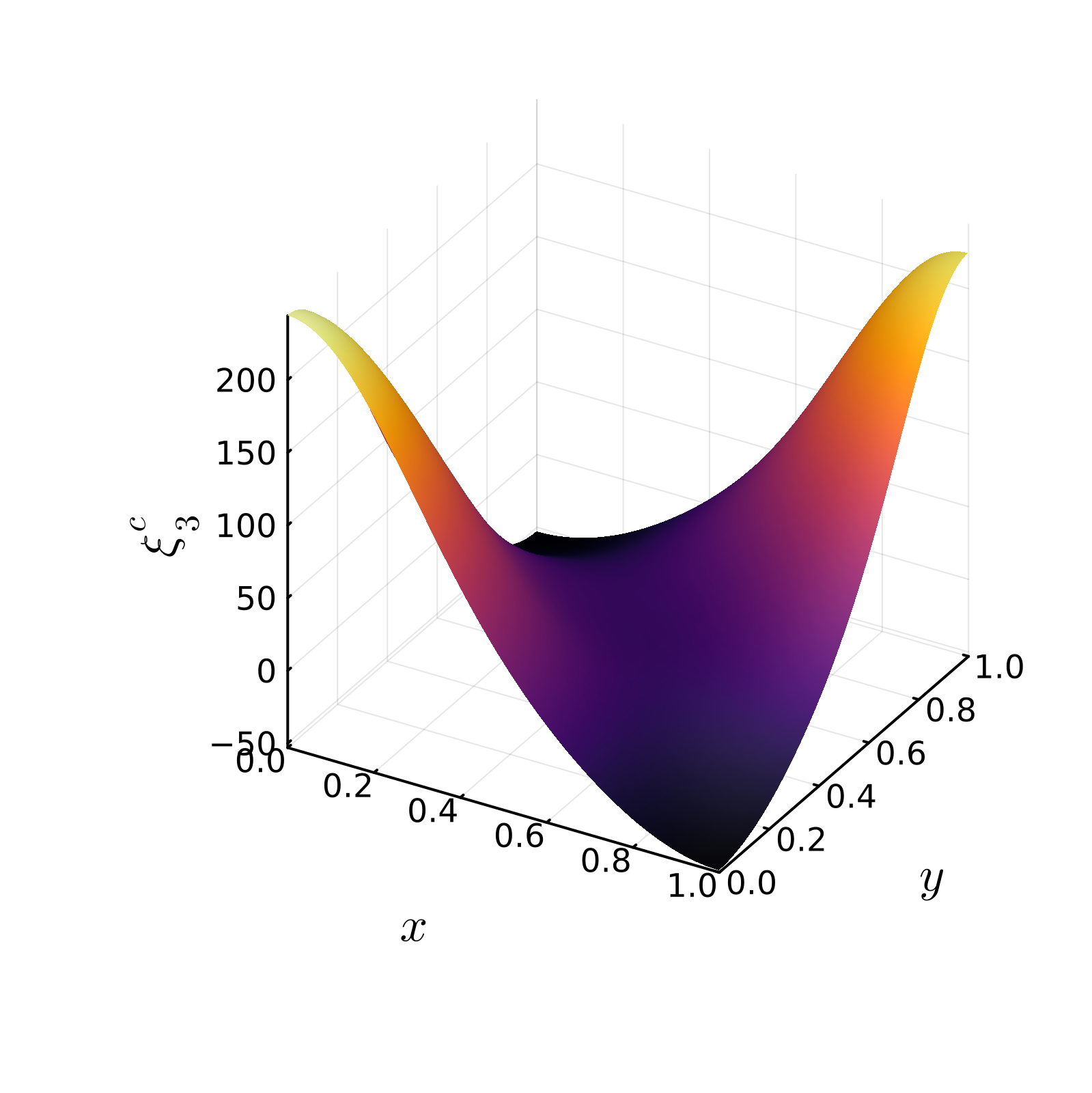}
	\includegraphics[width=0.23\textwidth]{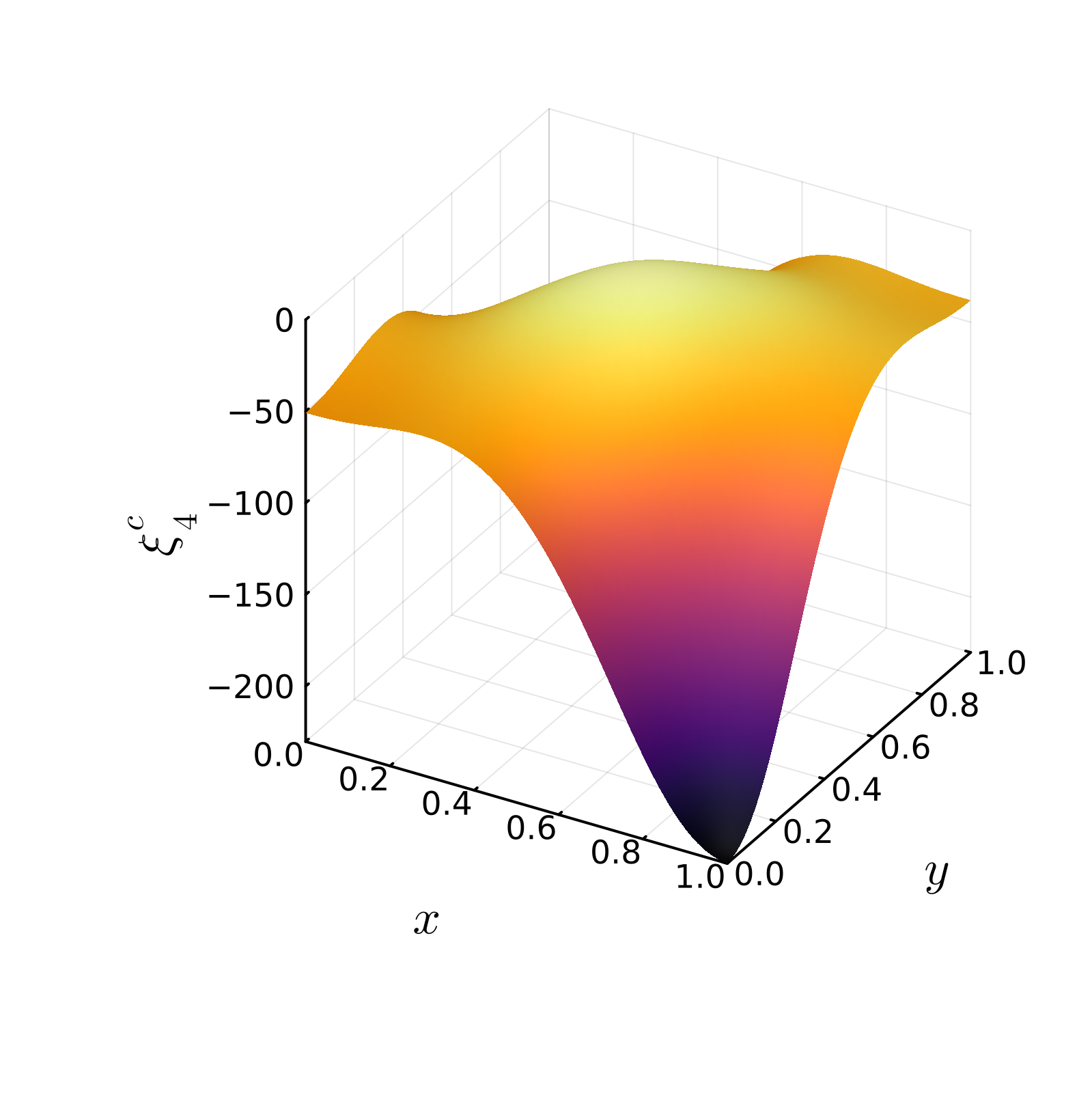} \\
	\includegraphics[width=0.23\textwidth]{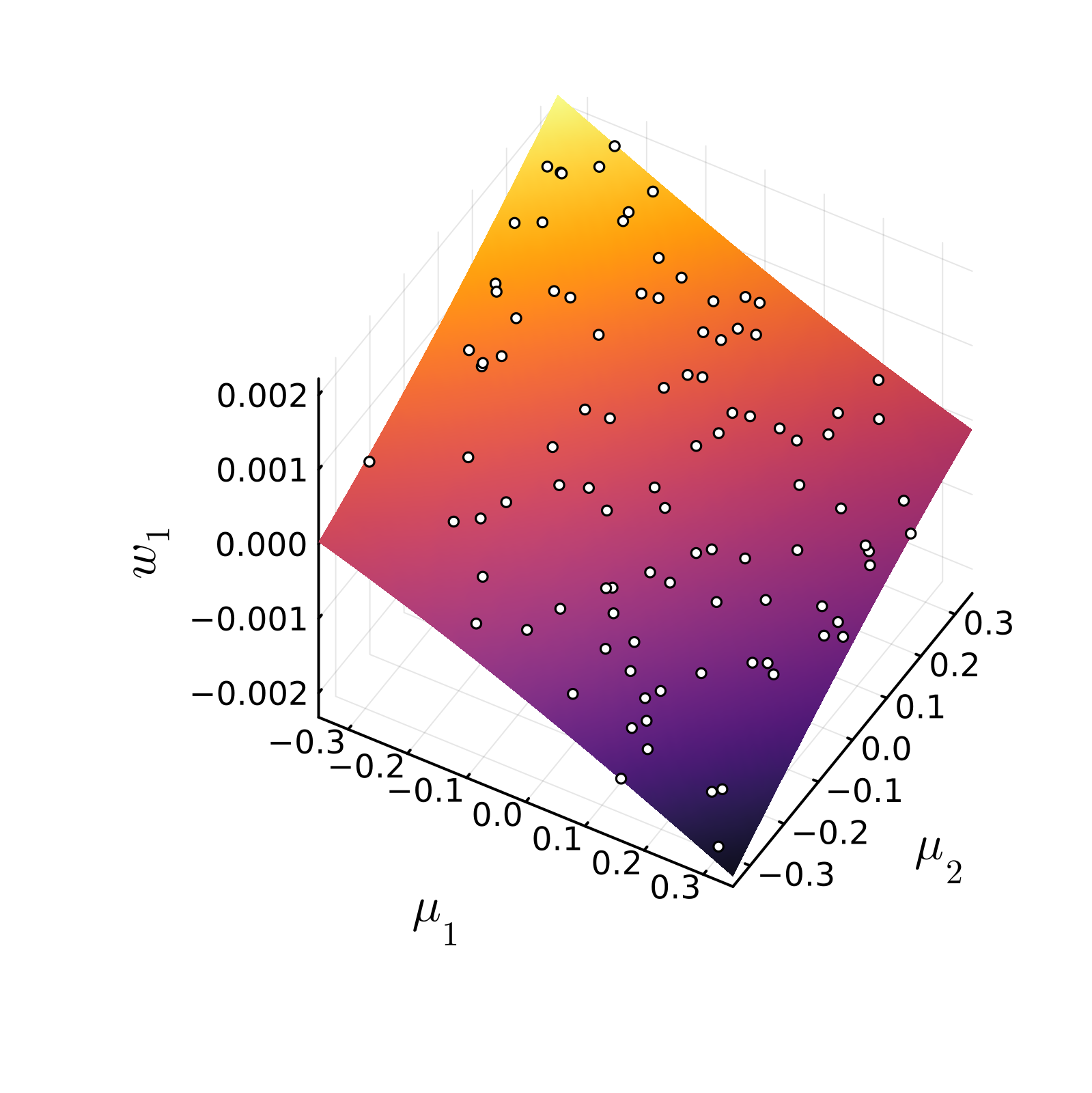}
	\includegraphics[width=0.23\textwidth]{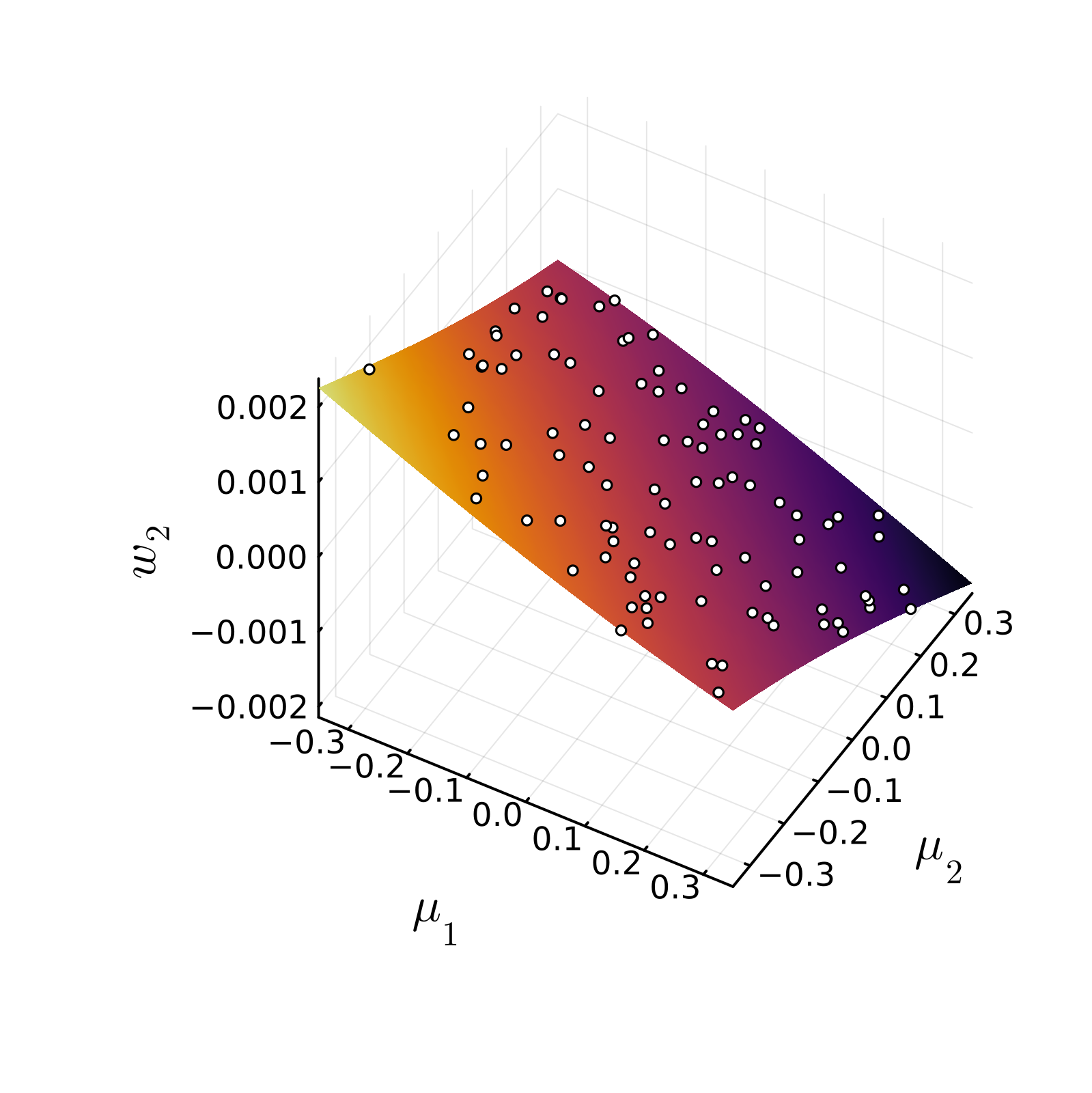}
	\includegraphics[width=0.23\textwidth]{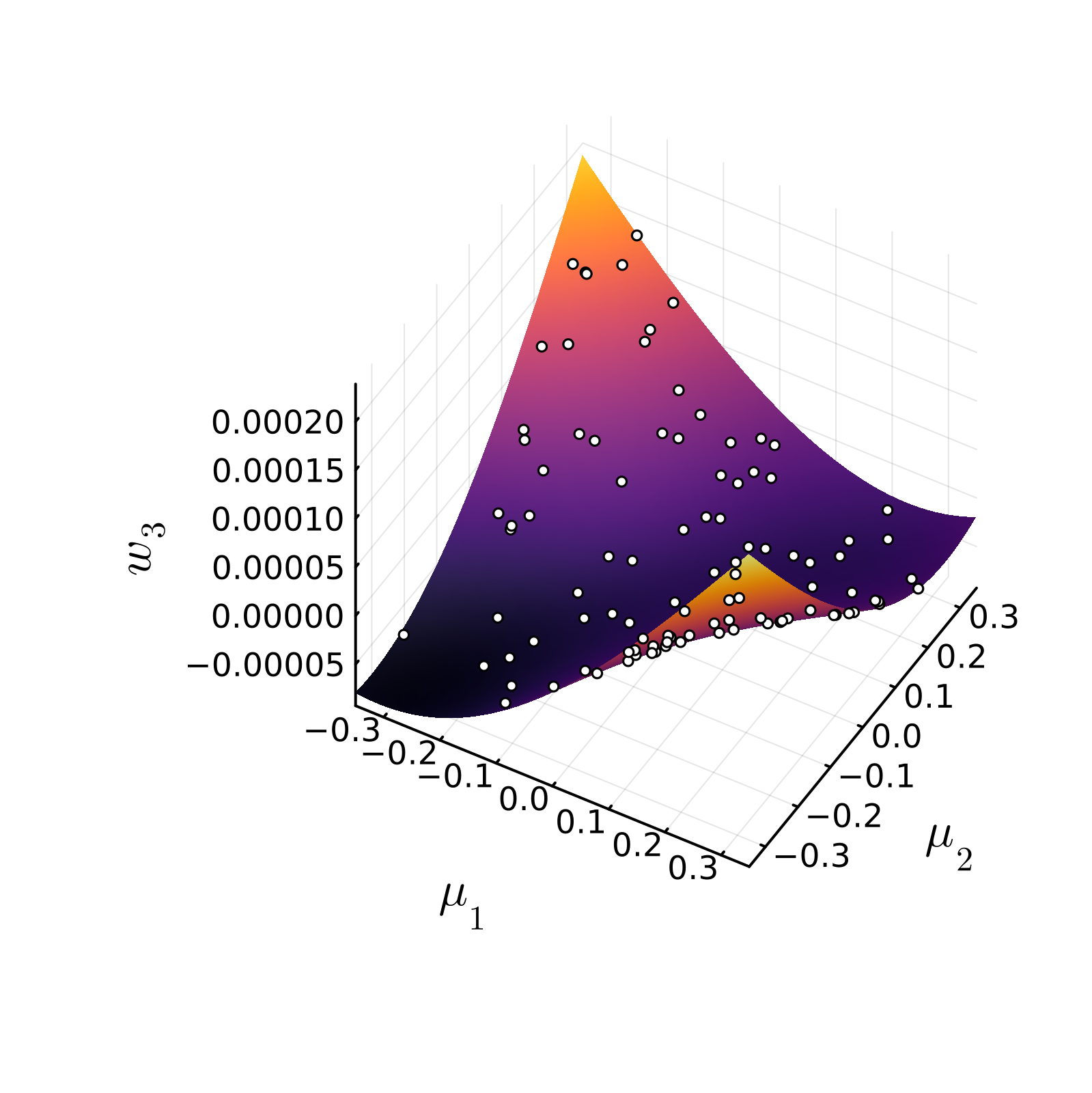}
	\includegraphics[width=0.23\textwidth]{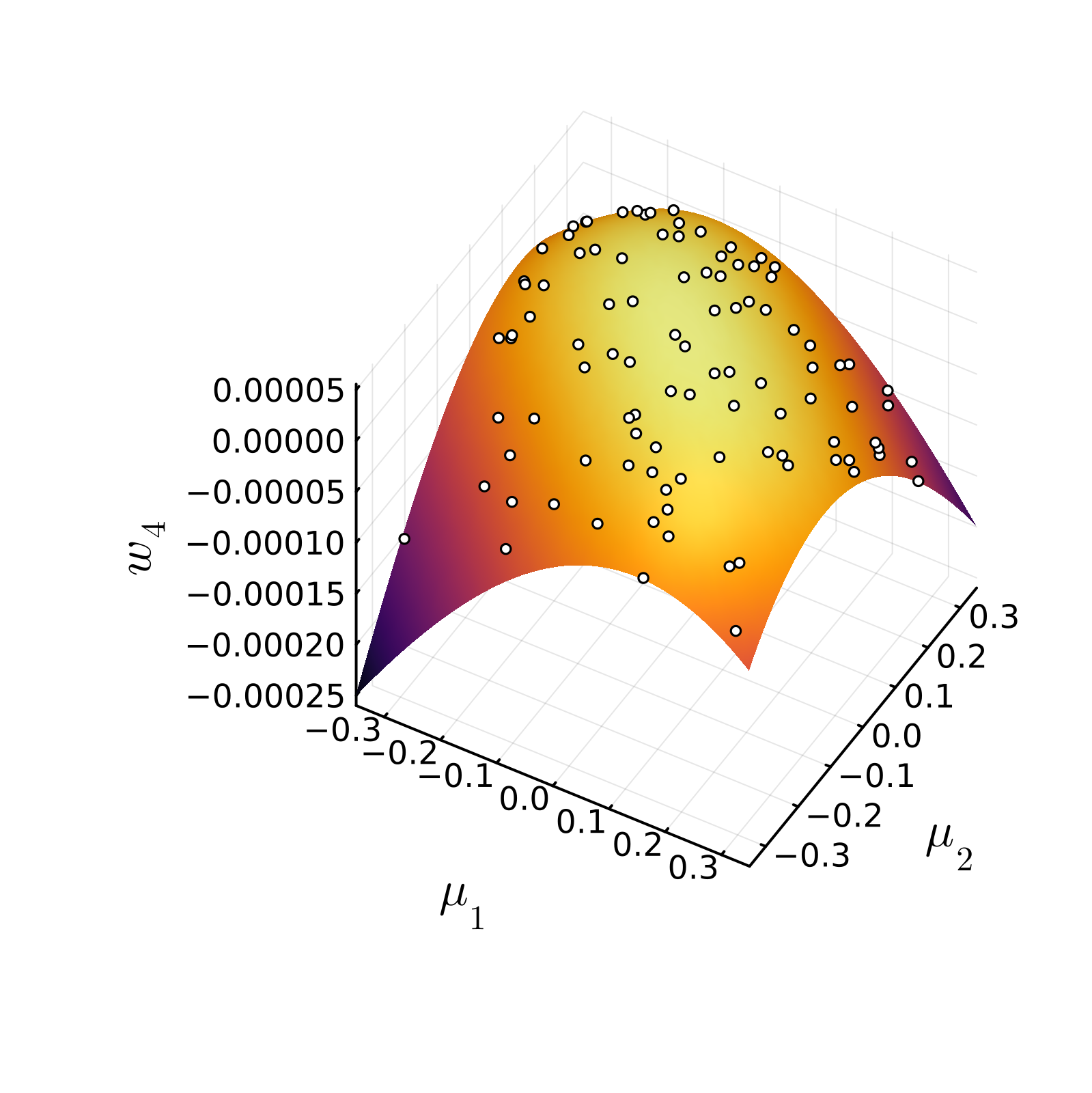}
	\caption{First row: first four transport modes with an added constant such that $\xi^c_{1,\dots,4}(\tfrac{1}{2},\tfrac{1}{2}) = 0$. Second row: Gaussian process approximation of the functions $w_{1,\dots,4}(\mu)$. Values used to construct the basis are marked in white. The parameters chosen as in \cref{fig:evds_1}.}
	\label{fig:xis_1}
\end{figure}

\Cref{fig:xis_1} shows the fist four transport modes $\xi^c_{1,\dots,4}$ and the Gaussian process approximations of $\mu \mapsto w_j(\mu)$, the transport mode coefficients used in the mapping $\Phi_\mu^{-1} = \mathrm{id} - \nabla \sum_{j=1}^m w_j(\mu) \xi^c_j$. The first two modes are essentially translations, the third mode is a contraction (or expansion, depending on the sign of its coefficient), and the fourth mode is a contraction along one diagonal and an expansion along the other.

As indicated by the very fast eigenvalue decay of the correlation matrix $\bC^\psi$, transport mappings can be approximated accurately as a linear combination of only a few transport modes. 
%To verify this, we calculate mismatch for the potentials $\{ \psi^c(\mu_j) \}_{1 \leq j \leq n_s}$ from the training set. The results are shown in \cref{tab:transport_errors} for different values of $m$. In particular, $m$ is chosen to match a certain retained eigenvalue energy 
%\begin{table}[h!]
%	\centering
%	\begin{tabular}{c c | ccc }
%		$\tau(\mathcal{E})$, & $m$ & \multicolumn{3}{c}{relative $L^2(\bar u)$ error of $\nabla \psi^c(\mu)$ } \\ \hline
%		& & avg & sd & max \\
%		$10^{-2}$ & 2 & $2.73 \times 10^{-2}$ & $5.54 \times 10^{-2}$ & $3.97 \times 10^{-1}$ \\
%		$10^{-3}$ & 4 & $4.76 \times 10^{-4}$ & $3.50 \times 10^{-4}$ & $1.93 \times 10^{-3}$ \\
%		$10^{-4}$ & 7 & $1.02 \times 10^{-4}$ & $1.04 \times 10^{-4}$ & $6.43 \times 10^{-4}$ \\
%	\end{tabular}
%	\caption{Average, standard deviation, and maximum of the approximation error of transport potentials in the training set for varying retained eigenvalue energy. Choosing $\tau(\mathcal{E}) = 10^{-1}$ also yields $m=2$.}
%	\label{tab:transport_errors}
%\end{table}
%
Next, we compare the error in the solution of the PPDE for the entire online phase. This includes approximating the mapping with transport modes, obtaining the coefficients of the transport modes with a Gaussian process, solving the PPDE in the reference domain as in \eqref{eq:laplace_reference}, and mapping the solution back to the physical domain as in \eqref{eq:OTRB_approx}. Average and maximum errors are calculated for a test set using $n_t = 50$ samples from $\mathcal{A}$. The results are compared to the classical POD method with no registration step, i.e. the $m=0$ case. The hyper-reduction uses EIM to evaluate the mapped Laplacian as described in \cref{sec:hyperreduction} as well as the right-hand-side term $\int_{\Omega} \phi_i f_\mu \det D\Phi_\mu^{-1} \, \rd y \; \forall i = 1,\dots,n_m$. The number of interpolation functions are denoted $Q_K$ and $Q_f$, respectively. Values for the cases of $m=0$ using hyper-reduction are not given, since the set $\{ f(\mu)\}_{\mu \in \mathcal A}$ shows extremely slow n-width decay. In contrast, the n-width decay of $\{ f(\mu) \circ \Phi_\mu^{-1} \, \det D\Phi_\mu^{-1}  \}_{\mu \in \mathcal A}$ allows the use of EIM albeit with large values of $Q_f$.

We also report the relative error of the $H^1$ semi-norm of $u(\mu)$, i.e. the error in the energy of the solution for an electrostatic problem. We choose this as an example of a quantity of interest that can be computed in the reference domain.

\begin{table}[h!]
	\centering
	\begin{adjustbox}{width=1\textwidth}
	\begin{tabular}{ccccc | cc | cc | cc}
		$\tau(\mathcal{E})$ & $n$ & $m$ & $Q_K$ & $Q_f$ & \multicolumn{2}{c | }{relative $L^2$ error of $u(\mu)$ }  & \multicolumn{2}{c | }{relative $H^1$ error of $u(\mu)$ } & \multicolumn{2}{c}{relative error of $\Vert u(\mu) \Vert_{\dot H^1}$ }  \\ \hline
		& & & & & avg 	& max & avg & max & avg & max \\
		%$10^{-3}$ 	& 41 & 0 & -& -   & $7.87 \times 10^{-2}$ & $3.32 \times 10^{-1}$ & $2.99 \times 10^{-1}$ & $7.08 \times 10^{-1}$ & $1.09 \times 10^{-1}$ & $5.07 \times 10^{-1}$ \\
		%					& 6 & 4 & - & -   & $ 4.46\times 10^{-2}$ & $1.59 \times 10^{-1}$ & $1.10 \times 10^{-1}$ & $2.76 \times 10^{-1}$ & $1.43 \times 10^{-2}$ & $7.77 \times 10^{-2}$ \\
		%					& 6 & 4 & 11&28 & $4.79 \times 10^{-2}$ & $1.68 \times 10^{-1}$ & $1.15 \times 10^{-1}$ & $3.24 \times 10^{-1}$ & $2.16 \times 10^{-2}$ & $8.76 \times 10^{-2}$ \\  \hline
		$10^{-3}$ 	& 41 & 0 & -& -   & $7.87 \times 10^{-2}$ & $3.32 \times 10^{-1}$ & $2.99 \times 10^{-1}$ & $7.08 \times 10^{-1}$ & $1.09 \times 10^{-1}$ & $5.07 \times 10^{-1}$ \\
							& 5 & 4 & - & -    & $4.48\times 10^{-2}$  & $1.27 \times 10^{-1}$ & $9.80 \times 10^{-2}$ & $2.16 \times 10^{-1}$ & $1.11 \times 10^{-2}$ & $4.67 \times 10^{-2}$ \\
							& 5 & 4 & 12&19  & $4.94 \times 10^{-2}$ & $1.25 \times 10^{-1}$ & $1.02 \times 10^{-1}$ & $2.20 \times 10^{-1}$ & $2.12 \times 10^{-2}$ & $7.13 \times 10^{-2}$ \\  \hline
		%$10^{-4}$  & 64 & 0 & -& -   & $5.85 \times 10^{-2}$ & $3.19 \times 10^{-1}$ & $2.28 \times 10^{-1}$ & $6.93 \times 10^{-1}$ & $7.75 \times 10^{-2}$ & $4.85 \times 10^{-1}$ \\
		%					& 10 & 6 & - & -   & $1.65 \times 10^{-2}$ & $5.46 \times 10^{-2}$ & $5.21 \times 10^{-2}$ & $1.22 \times 10^{-1}$ & $3.08 \times 10^{-3}$ & $1.45 \times 10^{-2}$ \\
		%					& 10 & 6 & 15&36& $2.02 \times 10^{-2}$ & $7.55 \times 10^{-2}$ & $5.47 \times 10^{-2}$ & $1.32 \times 10^{-1}$ & $1.17 \times 10^{-2}$ & $8.65\times 10^{-2}$ \\  \hline
		$10^{-4}$  & 64 & 0 & -& -   & $5.85 \times 10^{-2}$ & $3.19 \times 10^{-1}$ & $2.28 \times 10^{-1}$ & $6.93 \times 10^{-1}$ & $7.75 \times 10^{-2}$ & $4.85 \times 10^{-1}$ \\
							& 9 & 6 & - & -   & $1.53 \times 10^{-2}$ & $4.04 \times 10^{-2}$ & $4.27 \times 10^{-2}$ & $1.01 \times 10^{-1}$ & $1.71 \times 10^{-3}$ & $8.20 \times 10^{-3}$ \\
							& 9 & 6 & 19&24& $1.67 \times 10^{-2}$ & $5.10 \times 10^{-2}$ & $4.42 \times 10^{-2}$ & $1.08 \times 10^{-1}$ & $8.17 \times 10^{-3}$ & $4.70\times 10^{-2}$ \\  \hline
		%
		%$10^{-5}$  & 82 & 0 & -& -     & $5.10 \times 10^{-2}$ & $2.99 \times 10^{-1}$  & $2.02 \times 10^{-1}$ & $6.71 \times 10^{-1}$  & $6.68 \times 10^{-2}$ & $4.55 \times 10^{-1}$ \\
		%					& 15 & 9 & - & -    & $6.01 \times 10^{-3}$ & $1.78 \times 10^{-2}$ & $3.23 \times 10^{-2}$ & $5.96 \times 10^{-2}$ & $9.07 \times 10^{-4}$ & $2.71 \times 10^{-3}$ \\
		%					& 15 & 9 & 20&46& $6.83 \times 10^{-3}$ & $2.68 \times 10^{-2}$ & $3.27 \times 10^{-2}$ & $6.20 \times 10^{-2}$ & $1.98 \times 10^{-3}$ & $1.36 \times 10^{-2}$ \\
		$10^{-5}$  & 82 & 0 & -& -      & $5.10 \times 10^{-2}$ & $2.99 \times 10^{-1}$  & $2.02 \times 10^{-1}$ & $6.71 \times 10^{-1}$  & $6.68 \times 10^{-2}$ & $4.55 \times 10^{-1}$ \\
							& 11 & 10 & - & -    & $7.23 \times 10^{-3}$ & $2.36 \times 10^{-2}$ & $2.79 \times 10^{-2}$ & $6.20 \times 10^{-2}$ & $5.43 \times 10^{-4}$ & $2.32 \times 10^{-3}$ \\
							& 11 & 10 & 28&30& $9.88 \times 10^{-3}$ & $5.23 \times 10^{-2}$ & $3.02 \times 10^{-2}$ & $9.15 \times 10^{-2}$ & $6.26 \times 10^{-3}$ & $5.31 \times 10^{-2}$ \\
	\end{tabular}
	\end{adjustbox}
	\caption{PPDE solution errors in the test set as a function of the retained eigenvalue energy for $\varepsilon = 10^{-2}, \rho(\mu) \propto u(\mu)^2$, and de-biased calculations. }
		\label{tab:errors_solution}
\end{table}

\begin{table}[h!]
	\centering
	\begin{adjustbox}{width=1\textwidth}
		\begin{tabular}{ccccc | cc | cc | cc}
			de-biasing & $n$ & $m$ & $Q_K$ & $Q_f$ & \multicolumn{2}{c | }{relative $L^2$ error of $u(\mu)$ }  & \multicolumn{2}{c | }{relative $H^1$ error of $u(\mu)$ } & \multicolumn{2}{c}{relative error of $\Vert u(\mu) \Vert_{\dot H^1}$ }  \\ \hline
			& & & & & avg 	& max & avg & max & avg & max \\
			yes   & 9 & 6 & - & -   & $1.53 \times 10^{-2}$ & $4.04 \times 10^{-2}$ & $4.27 \times 10^{-2}$ & $1.01 \times 10^{-1}$ & $1.71 \times 10^{-3}$ & $8.20 \times 10^{-3}$ \\
					 & 9 & 6 & 19&24& $1.67 \times 10^{-2}$ & $5.10 \times 10^{-2}$ & $4.42 \times 10^{-2}$ & $1.08 \times 10^{-1}$ & $8.17 \times 10^{-3}$ & $4.70\times 10^{-2}$ \\  \hline
			no    & 10 & 7 & - & -   & $1.52 \times 10^{-2}$ & $5.42 \times 10^{-2}$ & $4.89 \times 10^{-2}$ & $1.02 \times 10^{-1}$ & $2.61 \times 10^{-3}$ & $9.97 \times 10^{-3}$ \\
					& 10 & 7 & 15&35& $1.97 \times 10^{-2}$ & $8.00 \times 10^{-2}$ & $5.19 \times 10^{-2}$ & $1.18 \times 10^{-1}$ & $1.22 \times 10^{-2}$ & $5.27\times 10^{-2}$ \\
		\end{tabular}
	\end{adjustbox}
	\caption{PPDE solution errors in the test set with and without de-biasing, using $\varepsilon = 10^{-2}, \tau = 10^{-4}, \rho(\mu) \propto u(\mu)^2$. }
	\label{tab:errors_debias}
\end{table}

In order to see how much of the overall error of the method is made when inverting $\Phi_{\mu}^{-1}$, we compare $u_{\mathrm{trb}} \circ \Phi_{\mu}^{-1}$ to $u \circ \Phi_{\mu}^{-1}$ in the case $\tau = 10^{-5}$ and using hyper-reduction. Calculating the approximation error in the reference domain in this way, we find average and maximum $L^2$ errors of $9.31 \times 10^{-3}$ and $5.62 \times 10^{-2}$. The average and maximum $H^1$ errors are $2.13 \times 10^{-2}$ and $8.61 \times 10^{-2}$. 
%Using no hyper-reduction, these errors increase to $5.04 \times 10^{-3}$, $2.57 \times 10^{-2}$, $1.47 \times 10^{-2}$, and $4.68 \times 10^{-2}$, respectively. 
We conclude that the approximation error in the reference domain dominates in the overall error of the method. For some parameter values, the two can even compensate each other.

%\begin{figure}[h!]
%	\centering
%	%\includegraphics[width=0.33\textwidth]{figs/ex1/crosssec_x}
%	%\includegraphics[width=0.33\textwidth]{figs/ex1/crosssec_y}
%	\includegraphics[width=0.5\textwidth]{figs/ex1/crosssecs}
%	\caption{Cross-sections along $[x, 0.5 + \mu_2]$ and $[0.5 + \mu_1, y]$ for the $\mu$-value with the largest relative $H^1$-error in the test set for the proposed method. $u$ denotes the high fidelity solution, $u_{\mathrm{rb}}$ the POD solution, $u_{\mathrm{trb}}$ the solution using the proposed method, and, in the case of $u_{\mathrm{trb,eim}}$, hyper-reduction using $\varepsilon = \tau = 10^{-3}$.}
%	\label{fig:crossecs}
%\end{figure}
%The solution with the largest $H^1$ error for the proposed method in the case $\tau(\mathcal{E}) = 10^{-4}$ is depicted in \cref{fig:crossecs}. %Upon close inspection, one can see small oscillations in the solution obtained with our method. These oscillations are an artefact of the mapping process back to the physical domain, i.e. the evaluation of $\phi_{1,\dots,n} \circ \Phi_\mu$. While we use $p=3$ Lagrange elements and can control the derivatives of $\Phi_\mu$ through the parameter $\varepsilon$, the finite element approximation of $\Phi_\mu$ is obtained as the gradient of an $H^1$ function and therefore discontinuous. We have verified that we can reduce the size of these oscillations using a finer mesh in numerical experiments. Another option might be to use a more regular approximation using e.g. $H^2$ conforming spline finite elements.

\subsubsection{The case $\rho(u)(\mu)= f(\mu)$}

This case is added here to see how well the effect of entropic smoothing can be used to apply the method even when $\rho$ takes very small values in $\Omega$. In this case, we do not use de-biased potentials and barycenters. Indeed, when using de-biased quantities, the performance of the method is heavily degraded with this choice of $\rho$.

\begin{figure}[h!]
	\centering
	\includegraphics[width=0.45\textwidth]{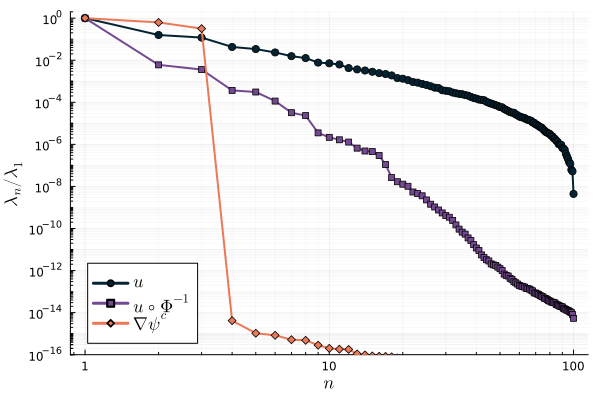}
	\includegraphics[width=0.45\textwidth]{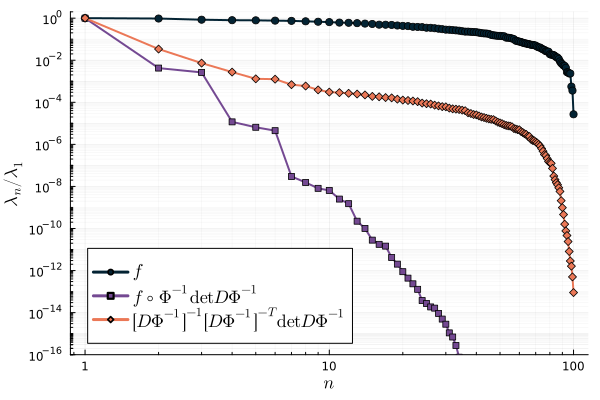}
	\caption{Left: Eigenvalues of the correlation matrices $\bC^u$, $\bC^{\Phi_\ast u}$, and $\bC^\psi$ for the case $\rho(\mu) = f(\mu)$. Right: Eigenvalues of the correlation matrices  $\bC^f$, $\bC^{K}$, and $\bC^{\Phi_\ast f}$ used in the hyper-reduction. No de-biasing, $\varepsilon = 10^{-2}, \tau = 10^{-4}, \rho(\mu) = f(\mu)$. }
	\label{fig:evds_f}
\end{figure}

As expected, there are only three eigenvalues of $\bC^\psi$ different from machine zero in this case - corresponding to two translations and one scaling mode (i.e. $y \mapsto \mathrm{const.} \cdot y^2$) which is a consequence of the entropic smoothing. The eigenvalue decay of $\bC^K$ and $\bC^{\Phi_* f}$ is also improved (\cref{fig:evds_f}). Approximation errors for this case are shown in \cref{tab:errors_solution_f} and are comparable to the case $\rho(u) \propto u^2$ even at $m = 3$. However, note that the online cost is mostly dependent on $n$, not $m$, see \cref{sec:runtimes}.

\begin{table}[h!]
	\centering
	\begin{adjustbox}{width=1\textwidth}
		\begin{tabular}{ccccc | cc | cc | cc}
			$\tau(\mathcal{E})$ & $n$ & $m$ & $Q_K$ & $Q_f$ & \multicolumn{2}{c | }{relative $L^2$ error of $u(\mu)$ }  & \multicolumn{2}{c | }{relative $H^1$ error of $u(\mu)$ } & \multicolumn{2}{c}{relative error of $\Vert u(\mu) \Vert_{\dot H^1}$ }  \\ \hline
			& & & & & avg 	& max & avg & max & avg & max \\
			%
			%$10^{-3}$ 	&41 & 0 & -& -   & $7.87 \times 10^{-2}$ & $3.32 \times 10^{-1}$ & $2.99 \times 10^{-1}$ & $7.08 \times 10^{-1}$ & $1.09 \times 10^{-1}$ & $5.07 \times 10^{-1}$ \\
			%					& 3 & 3 & - & -    & $7.45\times 10^{-2}$ & $2.03 \times 10^{-1}$ & $1.33 \times 10^{-1}$ & $2.69 \times 10^{-1}$ & $1.82 \times 10^{-2}$ & $6.83 \times 10^{-2}$ \\
			%					& 3 & 3 & 9&12   & $7.79 \times 10^{-2}$ & $2.19 \times 10^{-1}$ & $1.35 \times 10^{-1}$ & $2.71 \times 10^{-1}$ & $1.73 \times 10^{-2}$ & $9.86 \times 10^{-2}$ \\  \hline
			$10^{-3}$ 	&41 & 0 & -& -   & $7.87 \times 10^{-2}$ & $3.32 \times 10^{-1}$ & $2.99 \times 10^{-1}$ & $7.08 \times 10^{-1}$ & $1.09 \times 10^{-1}$ & $5.07 \times 10^{-1}$ \\
								& 3 & 3 & - & -    & $7.46\times 10^{-2}$ & $2.03 \times 10^{-1}$ & $1.33 \times 10^{-1}$ & $2.68 \times 10^{-1}$ & $1.81 \times 10^{-2}$ & $6.76 \times 10^{-2}$ \\
								& 3 & 3 & 9&12   & $7.84 \times 10^{-2}$ & $2.17 \times 10^{-1}$ & $1.35 \times 10^{-1}$ & $2.69 \times 10^{-1}$ & $1.58 \times 10^{-2}$ & $9.24 \times 10^{-2}$ \\  \hline
			%
			%$10^{-4}$  & 64 & 0 & -& -   & $5.85 \times 10^{-2}$ & $3.19 \times 10^{-1}$ & $2.28 \times 10^{-1}$ & $6.93 \times 10^{-1}$ & $7.75 \times 10^{-2}$ & $4.85 \times 10^{-1}$ \\
			%					& 6 & 3 & - & -    & $3.01 \times 10^{-2}$ & $8.88 \times 10^{-2}$ & $8.05 \times 10^{-2}$ & $1.23 \times 10^{-1}$ & $4.01 \times 10^{-3}$ & $1.29 \times 10^{-2}$ \\
			%					& 6 & 3 & 12&19 & $3.01 \times 10^{-2}$ & $8.12 \times 10^{-2}$ & $8.09 \times 10^{-2}$ & $1.27 \times 10^{-1}$ & $8.69 \times 10^{-3}$ & $4.77 \times 10^{-2}$ \\ \hline
			$10^{-4}$  & 64 & 0 & -& -   & $5.85 \times 10^{-2}$ & $3.19 \times 10^{-1}$ & $2.28 \times 10^{-1}$ & $6.93 \times 10^{-1}$ & $7.75 \times 10^{-2}$ & $4.85 \times 10^{-1}$ \\
								& 6 & 3 & - & -    & $3.00 \times 10^{-2}$ & $8.91 \times 10^{-2}$ & $8.03 \times 10^{-2}$ & $1.24 \times 10^{-1}$ & $4.00 \times 10^{-3}$ & $1.32 \times 10^{-2}$ \\
								& 6 & 3 & 12&19 & $3.00 \times 10^{-2}$ & $8.10 \times 10^{-2}$ & $8.08 \times 10^{-2}$ & $1.28 \times 10^{-1}$ & $8.68 \times 10^{-3}$ & $4.49 \times 10^{-2}$ \\ \hline
			%
			%$10^{-5}$  & 82 & 0 & -& -    & $5.10 \times 10^{-2}$ & $2.99 \times 10^{-1}$  & $2.02 \times 10^{-1}$ & $6.71 \times 10^{-1}$  & $6.68 \times 10^{-2}$ & $4.55 \times 10^{-1}$ \\
			%					& 9 & 3 & - & -     & $1.10 \times 10^{-2}$ & $2.72 \times 10^{-2}$ & $6.70 \times 10^{-2}$ & $8.90 \times 10^{-2}$ & $1.72 \times 10^{-3}$ & $5.50 \times 10^{-3}$ \\
			%					& 9 & 3 & 15&26 & $1.34 \times 10^{-2}$ & $4.30 \times 10^{-2}$   & $6.75\times 10^{-2}$ & $9.33 \times 10^{-2}$ & $6.59 \times 10^{-3}$ & $4.03 \times 10^{-2}$ \\
			$10^{-5}$  & 82 & 0 & -& -    & $5.10 \times 10^{-2}$ & $2.99 \times 10^{-1}$  & $2.02 \times 10^{-1}$ & $6.71 \times 10^{-1}$  & $6.68 \times 10^{-2}$ & $4.55 \times 10^{-1}$ \\
								& 9 & 3 & - & -     & $1.09 \times 10^{-2}$ & $2.75 \times 10^{-2}$ & $6.69 \times 10^{-2}$ & $9.01 \times 10^{-2}$ & $1.73 \times 10^{-3}$ & $5.74 \times 10^{-3}$ \\
								& 9 & 3 & 15&26 & $1.27 \times 10^{-2}$ & $3.85 \times 10^{-2}$  & $6.73\times 10^{-2}$  & $9.31 \times 10^{-2}$ & $5.38 \times 10^{-3}$ & $4.09 \times 10^{-2}$ \\
		\end{tabular}
	\end{adjustbox}
	\caption{PPDE solution errors in the test set as a function of the retained eigenvalue energy. No de-biasing, $\varepsilon = 10^{-2}, \tau = 10^{-4}, \rho(\mu) = f(\mu)$.}
	\label{tab:errors_solution_f}
\end{table}

Again, we compute the error in the reference domain to estimate the error induced by inverting the mapping. Using hyper-reduction, the average $L^2$ error is $6.71 \times 10^{-3}$ with a maximum of $2.76 \times 10^{-2}$. For the $H^1$ error, we find $1.60 \times 10^{-2}$ and $4.91 \times 10^{-2}$. We conclude that the error contributions are of the same order of magnitude in this case.

\subsubsection{Influence of the smoothing parameter}

The parameter $\varepsilon$ influences the method in two ways. Firstly, it acts as a hyper-parameter to balance fidelity and regularity of the mapping $\Phi$ itself. Secondly, when using $\bar \rho = \mathrm{OTBar_\varepsilon}\{ \rho_i \}_{i=1}^{n_s}$, $\varepsilon$ influences the shape of $\bar \rho$. Especially with no de-biasing, the entropic bias (see \cref{sec:entropic_bias}) leads to a smoothing of $\bar \rho$.

%When $\rho(\mu) = f(\mu)$, we know from \cref{sec:entropic_bias} that $\bar \rho \approx \mathcal{N}( \bar \mu := \tfrac{1}{n_s} \sum_i \mu_i, \mathrm{var} + \varepsilon )$. Therefore, $\Phi^{-1}_\mu(y) \approx \sqrt{1 + \tfrac{\varepsilon}{\mathrm{var}} }^{-1} (y + \mu - \bar \mu )$ and $\det D\Phi_{\mu}^{-1} \approx (1 + \tfrac{\varepsilon}{\mathrm{var}} )^{-1}$. Clearly, this degrades for $\varepsilon \gg \mathrm{var}$, which has been verified numerically. Therefore, when studying the dependence of the method on $\varepsilon$, one has to fix the value used to compute the reference density to a value such as $\varepsilon_{\mathrm{Bar}} = 10^{-2}$. This value is not optimized extensively and values in $[\tfrac{1}{2},2] \times \varepsilon_{\mathrm{Bar}}$ are also stable. These results are not shown here in the interest of brevity.

\begin{figure}[h!]
	\centering
	\includegraphics[width=0.3\textwidth]{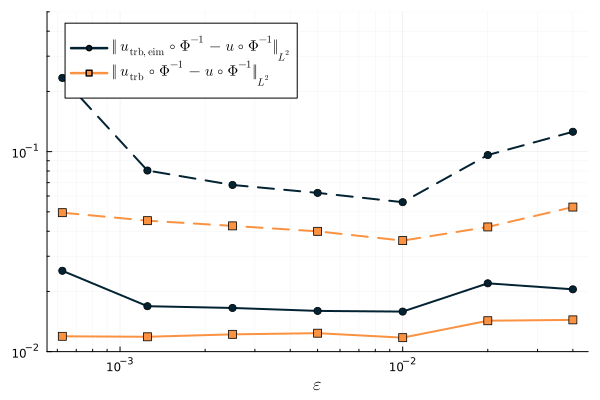}
	\includegraphics[width=0.3\textwidth]{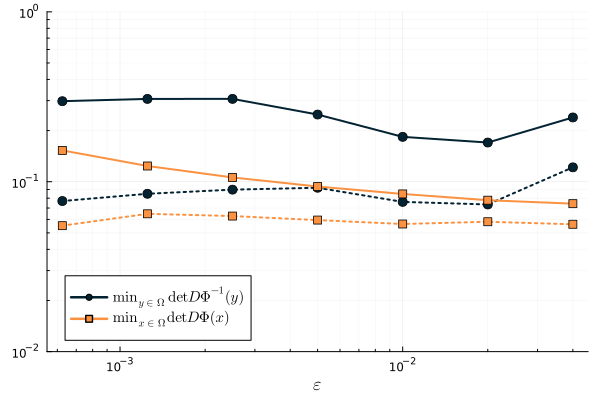}
	\includegraphics[width=0.3\textwidth]{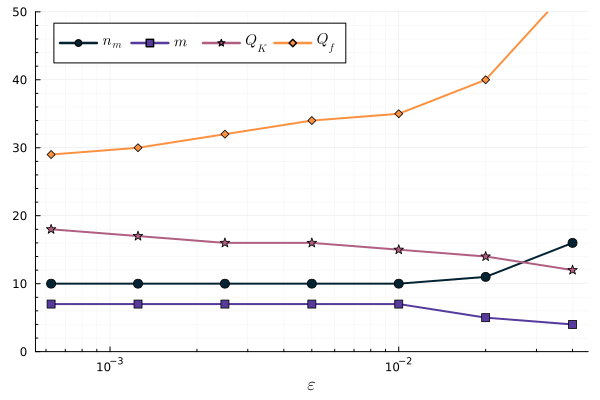} \\
	\caption{Influence of the regularisation parameter $\varepsilon$. In all cases, $\tau = 10^{-4} = 10\tau_\mathrm{eim}$, $\rho(\mu) \propto u(\mu)^2$, and no de-biasing was used. Left: relative average (solid) and maximum (dashed) $L^2$ error over the test set. Middle: minimum value of the determinant of the mapping and its inverse: average (solid) and minimum (dotted) over $\mu \in \mathcal{A}_{\mathrm{test}}$. The inverse mapping is calculated using the entropic c-transform with $\varepsilon_{\mathrm{fine}}$. Right: number of approximation modes.}
	\label{fig:eps_scaling}
\end{figure}

When $\rho(\mu) \propto u(\mu)^2$, the densities are supported on the full domain and we vary $\varepsilon$ for all calculations. We report results in \cref{fig:eps_scaling} for the case with no de-biasing. The results for the de-biased case are very similar. We give the error in the reference domain here to separate the effect of inaccuracies when inverting $\Phi_{\mu}^{-1}$.

Firstly and importantly, we observe that the approximation quality does not strongly depend on $\varepsilon \in [10^{-3}, 10^{-2}]$. Reducing $\varepsilon$ significantly below $10^{-3}$ would require all OT calculations to be moved to the log-domain, as $\min_{x,y} \exp \tfrac{ - c(x,y) }{ \varepsilon}$ becomes numerically zero in double precision.

As $\varepsilon$ becomes too large, the empirical interpolation method does not work as well any more, as the source terms $\{f(\mu_i)\}_i$ are no longer well aligned, leading to an increase in $Q_f$. As discussed in \cref{sec:boundary_cond}, the mapping can even be non-invertible in these cases. Note that the value $\varepsilon = 4 \times 10^{-2}$ corresponds to a characteristic scale of the transport of $\sqrt{\varepsilon} = 0.2$.

\subsubsection{Run times}\label{sec:runtimes}

%\begin{table}[h!]
%	\centering
%	\begin{adjustbox}{width=1\textwidth}
%		\begin{tabular}{cccc | ccccccccc | cc | ccc}
%			$n$ & $m$ & $Q_K$ & $Q_f$ & \multicolumn{9}{c | }{offline (all $\mu \in \mathcal A_{\mathrm{train}}$)} & \multicolumn{2}{c |}{online (per $\mu \in \mathcal A_{\mathrm{test}}$)} & \multicolumn{3}{c}{post-processing  (per $\mu \in \mathcal A_{\mathrm{test}}$)}  \\ \hline
%			& & & & \multicolumn{3}{c}{OT calculations} & gaussian process & mapping & reduced basis & assembly & \multicolumn{2}{c | }{hyper-reduction} & \multicolumn{2}{c |}{${u_\mathrm{rb}}$ or ${u_\mathrm{trb} \circ \Phi_\mu^{-1}} $} &  \multicolumn{3}{c }{re-mapping} \\
%			& & & & $\bar \rho$ and $\{\psi^{c}_{\mathrm{pre-proj.}}(\mu_i)\}_{i=1}^{n_s}$ & $\{\psi^{c}(\mu_i)\}_{i=1}^{n_s}$ & $\{ \xi^c_j \}_{j=1}^m$ & $\{ w(\mu) \}_{i=1}^{m}$ & $\{u(\mu_i) \circ \Phi_{\mu_i}^{-1} \}_{i=1}^{n_s}$ & $\{ \zeta_i \}_{i=1}^n$ or  $\{ \phi_i \}_{i=1}^{n_m}$ & $\int \nabla \zeta_i \cdot \nabla \zeta_j$ & $\mathrm{EIM}_K$ & $\mathrm{EIM}_f$ &  no EIM & EIM & $[ \sum_{j=1}^m w_j(\mu) \xi^c_j ]^c_{\mathrm{pre-proj.}}$ & $\Phi(\mu)$ & $u_{\mathrm{trb}}(\mu)$ \\
%			64 & - & - & - & - & - & - & - & - & 19s & 9.4s & - & - & 140ms & - & - & - & - \\
%			10 & 7 & 19 & 24 & 28s & 19s & 34s & 2.3s & 43s & 18s & - & 87s & 43s & 710ms & 6.1ms & 240ms & 190ms & 410ms \\
%		\end{tabular}
%	\end{adjustbox}
%	\caption{Run times}
%	\label{tab:runtimes}
%\end{table}

The codebase used to generate the results presented in this work is still being developed and parts of it have not yet been optimized. Thus, a thorough treatment of runtimes is beyond the scope of this work. Comparisons between the high-fidelity solver, the POD method, and the presented method depend on the size of high fidelity simulation $N$, the size of training- and test-set $n_s$ and $n_t$, and several other factors. Depending on the choice of these parameters, one or another method might seem favourable. Regardless, we show some examples in \cref{tab:runtimes}. The parameters chosen in these runs are a subset of those in \cref{tab:errors_solution} and \cref{tab:errors_solution_f}, where the corresponding errors can be found.

\begin{table}[h!]
	\centering
	\begin{adjustbox}{width=1\textwidth}
		\begin{tabular}{c | cccc}
			& \multicolumn{4}{c}{offline phase I: transport calculations (all $\mu \in \mathcal A_{\mathrm{train}}$)} \\ \hline
									& OT barycenter and potentials& boundary projection & transport modes & mapped snapshots \\
			$\rho(\mu)$  & $\bar \rho$ and $\{\psi^{c}_{\mathrm{pre-proj.}}(\mu_i)\}_{i=1}^{n_s}$ & $\{\psi^{c}(\mu_i)\}_{i=1}^{n_s}$ & $\{ \xi^c_j \}_{j=1}^m$ & $\{u(\mu_i) \circ \Phi_{\mu_i}^{-1} \}_{i=1}^{n_s}$\\ \hline
			$\propto u(\mu)^2$ & 27s 	& 19s 	& 33s 	& 43s \\
			$f(\mu)$ 					& 6.6s 	& 19s 	&33s 	& 43s \\
		\end{tabular}
	\end{adjustbox} \\

	\vspace{1em}

	\begin{adjustbox}{width=1\textwidth}
			\begin{tabular}{c cccc | ccccc }
					& & & & & \multicolumn{5}{c }{offline phase II (all $\mu \in \mathcal A_{\mathrm{train}}$)} \\ \hline
					& & & & & gaussian process & reduced basis & assembly & \multicolumn{2}{c }{hyper-reduction} \\
					$\rho(\mu)$ 			& $n$ & $m$ & $Q_K$ & $Q_f$& $\{ w(\mu) \}_{i=1}^{m}$ & $\{ \zeta_i \}_{i=1}^n$ or  $\{ \phi_i \}_{i=1}^{n_m}$ & $\int \nabla \zeta_i \cdot \nabla \zeta_j$ & $\mathrm{EIM}_K$ & $\mathrm{EIM}_f$ \\ \hline
					none 							 & 64 & -  & -    & -    & -        &  19s & 9.4s & -       & -	  \\
														 & 82  & -  & -    & - 	  & -        & 18s & 15s   & - 		& -     \\ \hline
					$\propto u(\mu)^2$ 	 & 10  & 7  & 19  & 24  & 2.3s  & 18s & -       & 86  & 42s \\
														 & 11  & 10 & 28 & 30  &3.6    & 17s & - 	  	& 93s  & 50s \\ \hline
					$f(\mu)$ 					  & 6   & 3  & 12  & 19   & 1.3s  & 18s & -       & 83s  & 43s \\
														 & 9   & 3  & 12  & 19   & 1.2s  & 18s & -       & 86s  & 43s \\
				\end{tabular}
		\end{adjustbox} \\

	\vspace{1em}

		\begin{adjustbox}{width=1\textwidth}
		\begin{tabular}{c | cccc | cc | ccc}
			& & & & &\multicolumn{2}{c |}{online (per $\mu \in \mathcal A_{\mathrm{test}}$)} & \multicolumn{3}{c}{post-processing  (per $\mu \in \mathcal A_{\mathrm{test}}$)}  \\ \hline
			& & & & & no EIM & EIM&  \multicolumn{3}{c }{re-mapping} \\
			$\rho(\mu)$ & $n$ & $m$ & $Q_K$ & $Q_f$ & \multicolumn{2}{c |}{${u_\mathrm{rb}}$ or ${u_\mathrm{trb} \circ \Phi_\mu^{-1}} $} & $[ \sum_{j=1}^m w_j(\mu) \xi^c_j ]^c_{\mathrm{pre-proj.}}$ & $\Phi(\mu)$ & $u_{\mathrm{trb}}(\mu)$ \\ \hline
			none & 64 & - & - 	& - 	& 0.14s & - & - & - & - \\
					&	82 & - & - 	& - 	& 0.17s & - & - & - & - \\ \hline
			$\propto u(\mu)^2$ 	 & 10 & 7 & 19& 24	 & 0.71s & 4.9ms  & 0.24s & 0.19s & 0.41s \\
												& 11 & 10& 28& 30	& 0.95s & 5.5ms & 0.24s & 0.20s & 0.41s \\ \hline
			$f(\mu)$					& 6 &	3	&12	& 19  & 0.43s & 4.5ms  & 0.24s & 0.19s	& 0.40s	\\
												& 9 &	3	&15	& 26  & 0.71s & 5.1ms  & 0.24s & 0.19s	& 0.41s	\\
		\end{tabular}
		\end{adjustbox}
		\caption{Run times for different choices of $\rho$. Parameters are as in \cref{tab:errors_solution} and \cref{tab:errors_solution_f}.}
	\label{tab:runtimes}
\end{table}

We clearly see the additional cost induced by the transport and registration. Hyper-reduction leads to large computational speed-up ($\approx 100$ to $200$) in the online phase by removing the dependence on $N$. We also see that $n_m$ has a larger impact on the cost of the method than $m$. The post-processing step to map the solutions back to the physical domain is costly, as it again depends on the size of the full order model.

In several applications, this last step is not needed. Quantities of interest can be obtainable in the reference domain, so the inversion of $\Phi_\mu^{-1}$ is not necessary. For example, the energy $\tfrac{1}{2} \int \nabla u \cdot \nabla u$ that we report here, or other linear functionals of the form ${l}(u) = \int u f_{{l}}$ can be calculated in the reference domain with no cost depending on the dimension of the full problem. For time-dependent problems, mapping back to the physical domain is only needed for diagnostics and plotting of the solution and thus usually not done at every time-step.

\subsection{Non-linear advection equation}\label{sec:nonlinear_advection}

As a second test case, we consider the equation
\begin{equation}
	\partial_t u(x,t) + \bar a (\theta) \cdot \nabla \left ( u(x,t) + \gamma u(x,t)^2  \right ) = \beta \Delta u(x,t),
\end{equation}
where $x \in \Omega = [0,1]^2$ and $t \in [0,T]$. The advecting velocity is given by $\bar a(\alpha) = \tfrac{1}{5} ( \cos \alpha, \sin \alpha )$, depending on the parameter $\alpha \in [0,2\pi]$. The strength of the non-linearity is set to $\gamma = 10^{-2}$ and $\beta$ is set to $10^{-3}$. The parameter space is therefore $\mathcal{A} := [0,1] \times [0,2\pi] \ni (t,\alpha) = \mu$. As an initial condition, we choose a Gaussian centred at $(\tfrac{1}{2},\tfrac{1}{2})$ with a variance of $5 \times 10^{-3}$. The solution is discretized using the same basis as in the previous example on a coarser $32 \times 32$ grid of quadrilaterals, using $N = 9025$ degrees of freedom. Time-integration is performed by an implicit midpoint method with time-step $\Delta t = 5 \times 10^{-2}$.

The reduced basis and transport modes are computed using the solutions at every time-step in $[0, T^{\mathrm{train}} = \tfrac{4}{5} ]$ for ten different values of $\alpha$ on a uniform grid between $0$ and $2\pi$, thus $n_s = 170$. The solutions $u$ themselves are used as the densities $\rho(u) = u$ for the OT computations, which are performed on a $96 \times 96$ grid. We set $\varepsilon = 10^{-2}$ and use the OT barycenter of the training set as a reference density. No de-biasing is used. We let $\tau_{\mathrm{eim}} = \tau = 10^{-3}$. All other parameters are identical to the previous example.

%
%\subsubsection{Transport mappings}
%
%Figure \ref{fig:xis_2} shows the first three transport modes and Figure \ref{fig:ws_2} depicts the corresponding weight functions $(t,\theta) \mapsto w_{1,2,3}(t,\theta)$. 
%\begin{figure}[h!]
%	\centering
%	\includegraphics[width=0.3\textwidth]{figs/ex2/xi1}
%	\includegraphics[width=0.3\textwidth]{figs/ex2/xi2}
%	\includegraphics[width=0.3\textwidth]{figs/ex2/xi3}
%	\caption{First three transport modes with an added constant such that $\xi^c_{1,\dots,3}(0,0) = 0$.}
%	\label{fig:xis_2}
%\end{figure}
%\begin{figure}[h!]
%	\centering
%	\includegraphics[width=0.3\textwidth]{figs/ex2/w1}
%	\includegraphics[width=0.3\textwidth]{figs/ex2/w2}
%	\includegraphics[width=0.3\textwidth]{figs/ex2/w3}
%	\caption{The Gaussian process approximation of the functions $w_{1,\dots,3}(t,\theta)$.}
%	\label{fig:ws_2}
%\end{figure}
%
%\subsubsection{Approxmiation of the solution}

\begin{figure}[h!]
	\centering
	\includegraphics[width=0.45\textwidth]{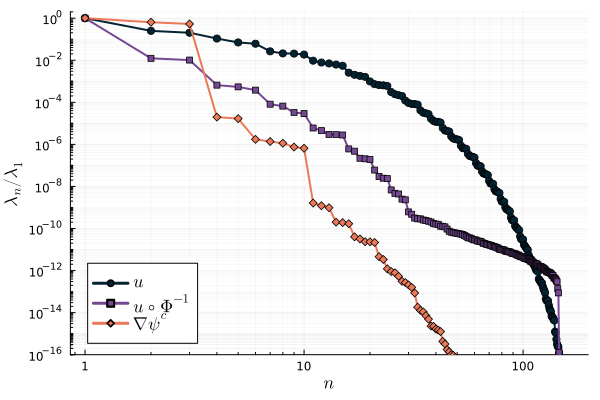}
	\includegraphics[width=0.45\textwidth]{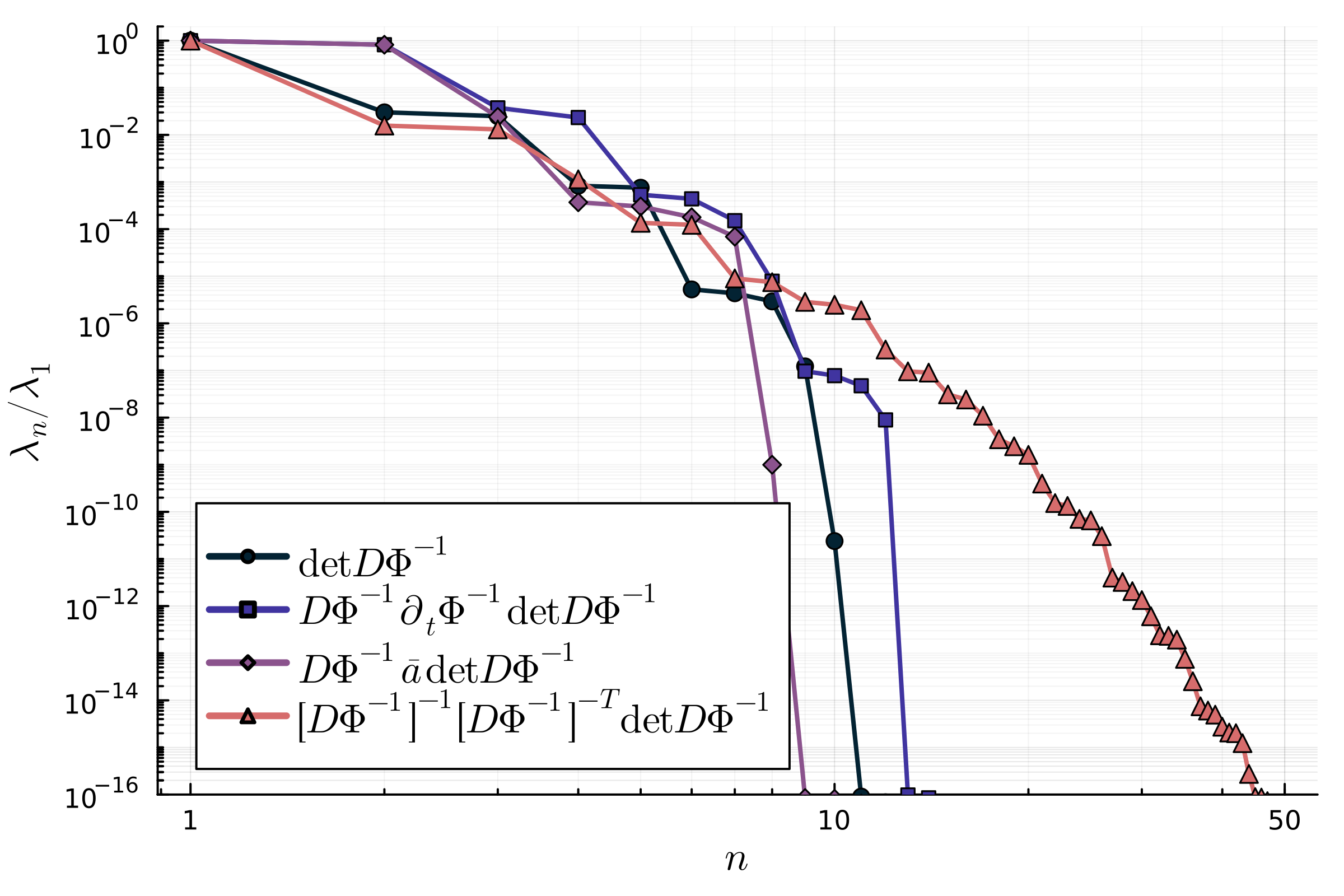}
	\caption{Left: Eigenvalues of the correlation matrices $\bC^u$, $\bC^{\Phi_\ast u}$, and $\bC^\psi$. Right: Eigenvalues of the correlation matrices used in the hyper-reduction. }
	\label{fig:evds_2}
\end{figure}

The selected energy criterion leads to $n = 24$ for the classical RB method and $n_m = 5, m=3$ for the proposed one. For the hyper-reduction, the RB method with $m=0$ requires only an EIM approximation for $\bar a$, which can be done to machine precision with $Q = 2$ modes. The method with transport requires EIM approximations for $\det D\Phi_{\mu}^{-1}$, $D\Phi_{\mu}^{-1} \partial_t \Phi_{\mu}^{-1} \det D\Phi_{\mu}^{-1}$, $D\Phi_{\mu}^{-1} \bar a(\mu) \, \det D\Phi_{\mu}^{-1}$, and $[D\Phi_{\mu}^{-1}]^{-1} [D\Phi_{\mu}^{-1}]^{-T} \det D\Phi_{\mu}^{-1}$. The values of $Q$ for these terms are $4, 4, 3,$ and $4$, respectively. The corresponding eigenvalue decays are shown in \cref{fig:evds_2}.

%The online phase of the proposed method takes approximately 15 seconds per snapshot (plus 2.5 seconds to reconstruct the solution in the physical space at $t = T$.) The RB method with $m=0$ is much faster (approximately $5.0 \times 10^{-2}$ seconds), as the forms are time-independent in this example.

\begin{figure}[h!]
	\centering
	\includegraphics[width=0.45\textwidth]{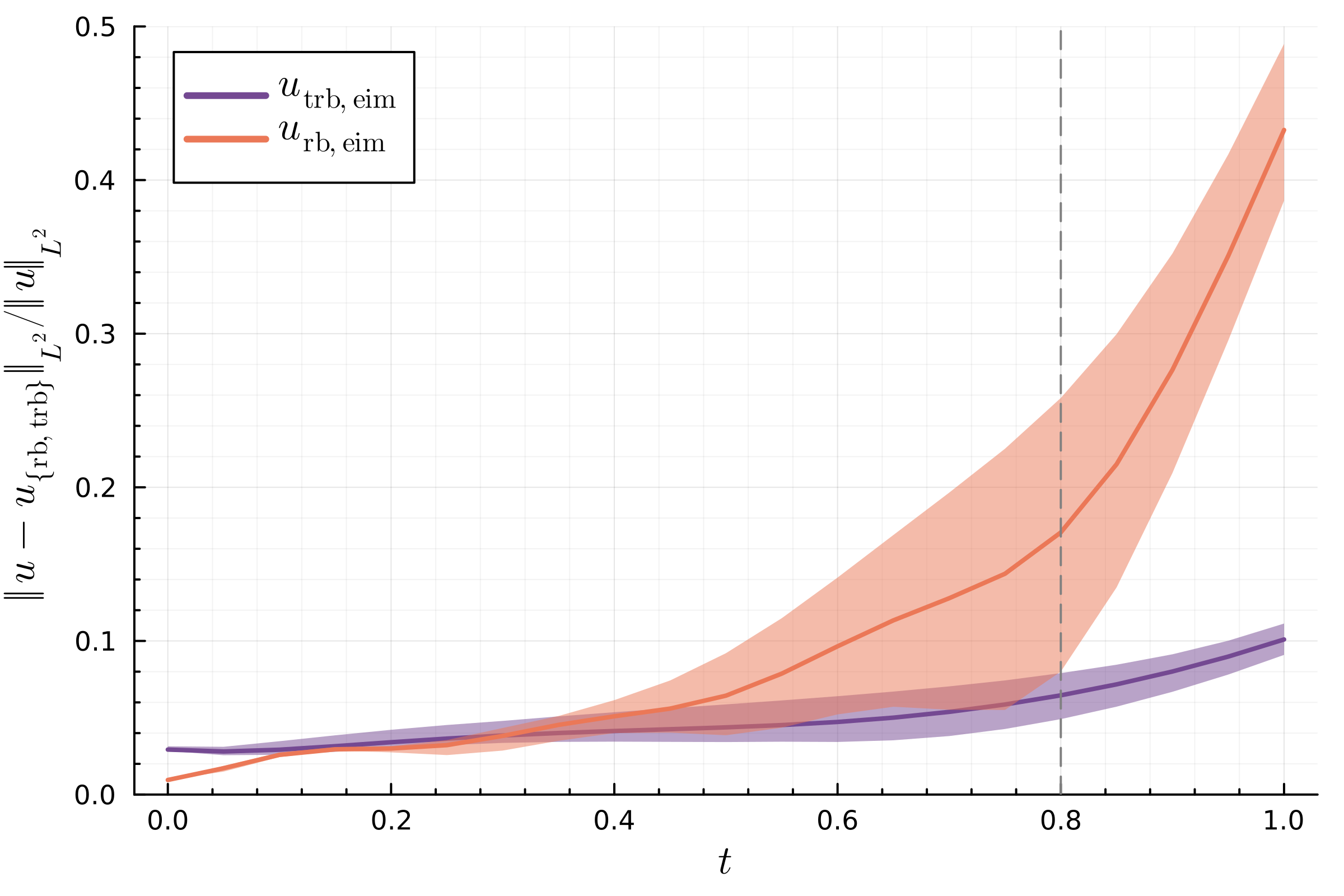}
	\includegraphics[width=0.45\textwidth]{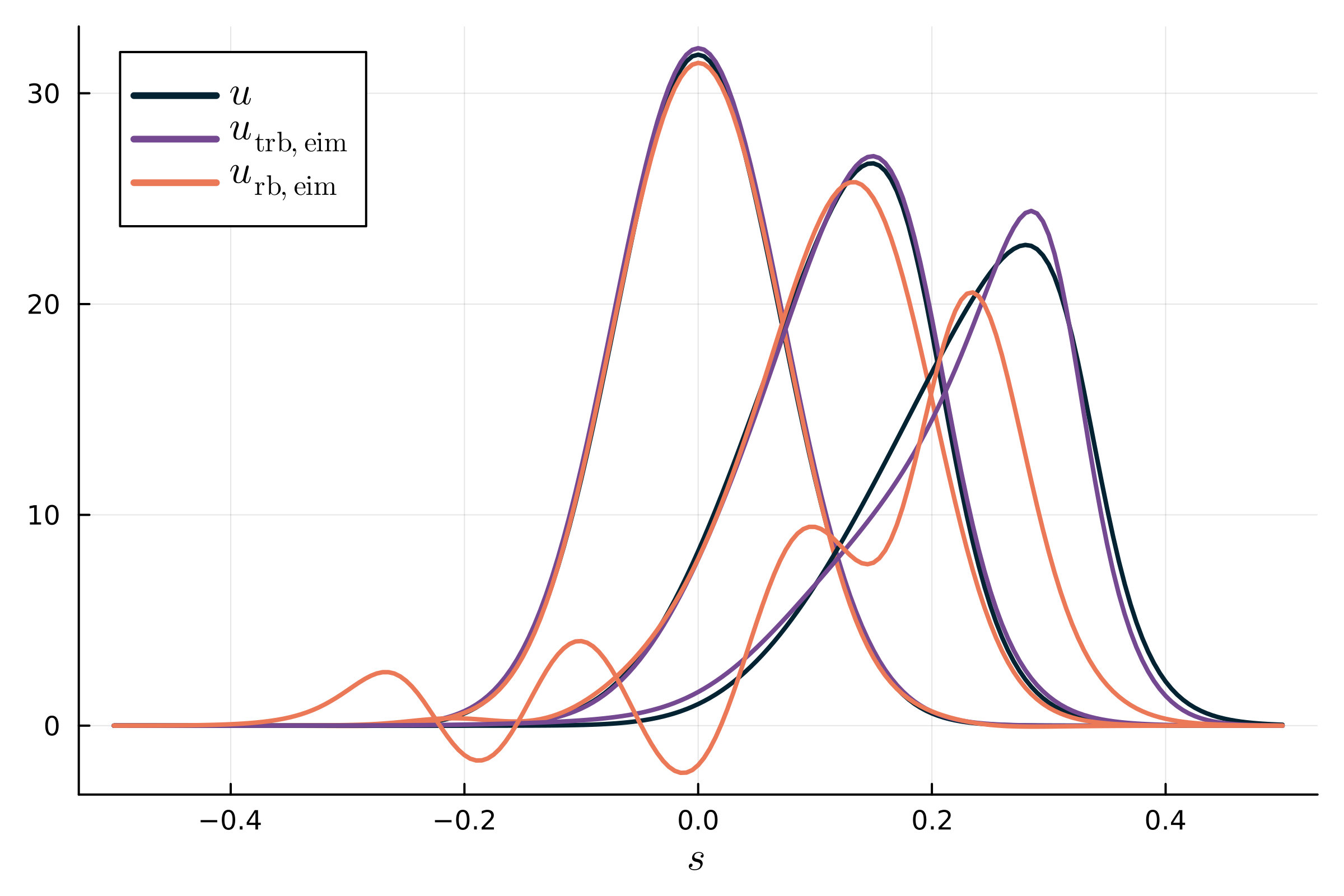}
	\caption{Left: Relative $L^2$ errors for the non-linear advection-diffusion equation as a function of time. Plotted is the average error across all ten values of $\alpha$. The shaded area is bounded above and below by the maximum and minimum error. The dashed line indicates $T^{\mathrm{train}}$ and the beginning of the extrapolation region.  Right: cross-sections in propagation direction (parametrized by $s$) for the value of $\alpha$ where the proposed method performs the worst and for $t \in \{0, \tfrac{1}{2}, 1\}$.}
	\label{fig:errs_ex2}
\end{figure}

\Cref{fig:errs_ex2} shows the relative $L^2$ errors over time for $20$ values of $\alpha$ randomly chosen from $[0, 2\pi]$ and for $t \in [0, T^{\mathrm{test}} = 1]$. The average approximation $L^2$ error of $\{ u(T^{\mathrm{test}}, \alpha_j) \}_{1 \leq j \leq 10}$ in the proposed method is $1.01 \times 10^{-1}$, while the maximum is $1.11 \times 10^{-2}$. We see that the classical RB method is only accurate for solutions close to the initial condition and for some select values of $\alpha$ close to those in the training set. This issue cannot be remedied by adding more reduced basis functions, since the solutions for values of $\alpha$ not in the training set and $t > T^{\mathrm{train}}$ cannot be expressed by any linear combination of training snapshots. In contrast, the proposed method yields qualitatively correct results even with the low number of modes employed. 
\section{Summary}\label{sec:summary}

In this work, we presented a data-driven registration method based on Monge embeddings. In the simple numerical examples we considered, the method proved robust and showed the expected improvements compared to an approach without registration.

Several possible extensions and open questions regarding the method have been mentioned throughout the text. For example, would be desirable to enforce the mapping boundary conditions in the construction itself instead of a later projection step. Furthermore, the OT computations, which rely mostly on convolution operations, could be performed in the same finite element framework as the PPDE solution itself. In the present work, the potentials are represented by $H^1$ conforming finite element functions, which are crucially not of class $C^1$. We would expect the method to profit from a $C^2$ conforming approximation, e.g. using a spectral polynomial basis.

Due to the parameter-dependence of the mapped forms, hyper-reduction methods are crucial to make the method performant in the online phase. Several other hyper-reduction methods exist that have been applied successfully in many applications, such as the empirical quadrature method \cite{yano_lp_2019}. 

The method requires three central inputs from the user: An entropic regularization parameter $\varepsilon$, a energy cut-off criterion $\tau$, and a choice of density $\rho$. In order to stay close to the true OT problem, $\varepsilon$ should be chosen small, keeping in mind that the transport plan will essentially ignore all features smaller than $\sqrt{\varepsilon}$. The choice of $\tau$ depends on the desired accuracy versus cost of the approximation. Selecting $\rho$ is the most intricate issue (of course, if the PPDE solution itself is a probability density, this is a natural choice). For maximum robustness, $\rho$ should be bounded from below on the entire domain by a positive constant. The size of this constant can have a significant effect on $\Vert \Phi_{\mu} \Vert$ as discussed in \cref{sec:regularity}. The numerical experiments indicate that entropic smoothing of the reference density can make the method applicable in some cases where $\rho$ takes very small values in the domain.

\section{Acknowledgements}

The presented work is supported by the Helmholtz Association and the Munich School of Data Science. \\
We would like to thank Michael Kraus and Olga Mula for stimulating discussions about the method and its numerical implementation.

\appendix
\section{Notation} The following table repeats some of the notation that is repeatedly used throughout the article:\\

\begin{longtable}{l | l}
	$\langle \cdot, \cdot \rangle_H$ & inner product in the Hilbert space $H$ \\
	$\dtimes, \ddiv, \dots$ & element-wise operations, e.g. $(\hat a \dtimes \hat b)_i = \hat a_i \ast \hat b_i$;  \\
	$\hat{(\cdot)}$ & collocated quantity $\hat a_i = a(x_i)$ \\
	$\oplus$ & notation for $(\psi_\rho \oplus \psi_\sigma)(x,y) = \psi_\rho(x) + \psi_\sigma(y)$  \\
	$(\cdot)_\sharp$ & push-forward operation for a density, see \cref{eq:push-forward} \\
	$(\cdot)^*$ & Legendre transform $f^*(y) := \sup_x \left ( \langle x, y \rangle - f(x) \right )$ \\
	$\mathbbm{1}$ & indicator function: $\mathbbm{1}_{\Omega'}(x) = 1$ if $x \in \Omega'$ and zero otherwise \\
	$\mathcal A$ & parameter space \\
	%$a, b$ & scaling functions $a: x \mapsto \exp \left ( \tfrac{1}{\varepsilon} \psi^\varepsilon_\rho(x) \right )$ and $b: y \mapsto \exp \left ( \tfrac{1}{\varepsilon} \psi^\varepsilon_\sigma(y) \right )$ \\
	$B$ & interpolation matrix used in empirical interpolation, see \cref{sec:hyperreduction} \\
	$c$ & cost function $c(x,y) = \tfrac{1}{2} \vert x - y \vert^2$ \\
	$(\cdot)^c$ & c-transform, see \cref{def:c-transform} \\
	$\bC^u$ & correlation matrix of $u_1,\dots,u_{n_s}$, see \cref{def:POD} \\
	$\cC_b$ & space of bounded continuous functions \\
	$\mathrm{cdf}$ & cumulative distribution function \\
	$\mathcal{E}(n; \lambda) $ & retained eigenvalue energy, see \cref{eq:eigenvalue_energy} \\
	$h$ & grid width \\
	$\mathrm{id}, \mathrm{Id}$ & identity $x \mapsto x$, identity matrix \\
	$k$ & Gibbs kernel $k(x,y) = \exp \left ( - \tfrac{1}{\varepsilon} c(x,y) \right )$ \\
	$\mathcal{L}(\cdot; \mu)$ & parameter-dependent PDE operator $V_h \rightarrow \bR$ \\
	$m$ & number of transport modes, see \cref{def:transport_modes} \\
	$\mathcal{M}$ & solution manifold, see \cref{eq:solution_manifold} \\
	$\max_{(\cdot) \sim u}^\varepsilon, \min_{(\cdot) \sim u}^\varepsilon$ & softmax and softmin operations, see \cref{def:softmin} \\
	$n, n_m$ & number of reduced basis functions, see \cref{def:POD} and \cref{eq:OTRB_approx} \\
	$n_s, n_t$ & number of snapshots in the training and test set \\
	$N$ & dimension (i.e. number of degrees of freedom) of the high fidelity discretization \\
	$\mathcal N$ & normal distribution \\
	$\mathcal{P}$ & set of probability measures \\
	$\LOT(\rho,\sigma)$ & Linear Optimal Transport distance between $\rho$ and $\sigma$, see \cref{def:lot_distance} \\
	$\mathcal O$ & $f = \mathcal{O}(1 + x^k) \Leftrightarrow$ $f$ is bounded by an order $k$ polynomial with positive coefficients  \\
	$\OT(\rho,\sigma)$ & Optimal Transport or Wasserstein distance between $\rho$ and $\sigma$, see \cref{def:OT_dist} \\
	$\mathrm{OTBar}, \mathrm{LOTBar} $ & (Linear) Optimal Transport barycenter, see \cref{def:OT_barycenter} and \cref{def:lot_barycenter} \\
	$Q$ & number of EIM modes, see \cref{sec:hyperreduction}. \\
	$T$ & Transport or Monge map, see \cref{eq:Transport_map} \\
	$u, v$ & elements of $V_h$ \\
	$\tilde u$ & coefficients used to approximate $u_\mathrm{trb}$, see \cref{eq:OTRB_approx} \\
	$\mathrm{v}$ & matrix eigenvector \\
	$V_h$ & discretized function space \\
	$w$ & transport mapping coefficient, see \cref{sec:embedding_reduction} \\
	$X$ & EIM mode, see \cref{sec:hyperreduction} \\
	$\delta$ & penalization parameter, set to $10^{-9}$ or smaller \\
	$\varepsilon$ & entropic regularization parameter, see \cref{eq:OT_entropic_primal} \\
	$\varepsilon_{\mathrm{fine}}$ & regularization used when inverting $\Phi^{-1}$. Set to the order of $h^2$. \\
	$\zeta$ & POD basis function, see \cref{def:POD} \\
	$\theta$ & EIM coefficients, see \cref{sec:hyperreduction}. \\
	$\kappa$ & $H^1$ projection parameter, see \cref{sec:boundary_cond}. Set to $\varepsilon^{-1/2}$ \\
	$\lambda$ & matrix eigenvalue, non-increasing: $\lambda_1 \geq \lambda_2 \geq \dots$ \\
	$\mu$ & parameter in $\mathcal A$ \\
	$\pi$ & transport plan, see \cref{eq:OT_primal} \\
	$\Pi(\rho, \sigma)$ & set of admissible transport plans between $\rho$ and $\sigma$, see \cref{eq:adm_plans} \\
	$\rho, \sigma$ & probability densities \\
	$\bar \rho$ & reference density, see \cref{def:monge_embedding} and \cref{sec:embedding_reduction} \\
	$\varsigma$ & measure used in the regularisation term of entropic OT, see \cref{eq:OT_entropic_primal} \\
	$\tau, \tau_{\mathrm{eim}}$ & energy criterion for POD: $1 - \mathcal{E}(n; \lambda) < \epsilon$ defines $n$ \\
	$\varphi$ & $x \mapsto \tfrac{|x|^2}{2} - \psi(x)$, see \cref{thm:brenier} \\
	$\phi$ & POD basis function in the reference domain, see \cref{eq:OTRB_approx} \\
	$\Phi_\mu$ & parameter-dependent mapping, see \cref{eq:mapped_manifold} \\
	$\xi^c$ & transport mode, see \cref{def:transport_modes} \\
	$\Xi$ & POD modes used in the EIM construction, see \cref{sec:hyperreduction} \\
	$\psi$ & transport potential, see \cref{eq:OT_dual} \\
	$\omega$ & barycenter weight, see \cref{def:OT_barycenter}. \\
	$\Omega$ & physical domain $\subset \bR^d$ \\
\end{longtable}
\section{Optimal transport computations}\label{sec:computation}

In this section, we will provide some details on how we compute the OT ingredients from Section \ref{sec:method}. The routines are available in the package \texttt{WassersteinDictionaries.jl}\footnote{https://github.com/JuliaRCM/WassersteinDictionaries.jl}.

\subsection{Sinkhorn's algorithm}

We rely on the tools from entropic OT for our computations. In order to not needlessly bloat the notation, we will omit the superscript $^\varepsilon$ from the entropic transport potential $\psi^\varepsilon$ and Monge map $T^\varepsilon$. Furthermore, we will write $N$ (rather than $N_{\mathrm{fine}})$ for the dimension of the discrete problem, which does not necessarily coincide with the dimension of the discrete space used to approximate the PPDE solution. \\
Using the regularized formulation of OT allows us to leverage very fast algorithms developed in recent years, which we repeat for convenience. Suppose that the quantities are discretized by collocation, i.e. $\hat \rho_i = \rho(x_i) : 1 \leq i \leq N$, $\hat \sigma_i = \sigma(x_i) : 1 \leq i \leq N$, and  $\hat C_{ij} = c(x_i, y_j) : 1 \leq i,j \leq N$. 
\begin{algorithm}[h]
	\caption{Sinkhorn's algorithm}\label{alg:sinkhorn}
	\begin{algorithmic}[1]
		\Function{sinkhorn}{$\hat \rho,\hat \sigma, \hat C, \varepsilon, \mathtt{tol}$}
			\State $\hat a, \hat b \dgets \tone$ \Comment{. denotes element-wise operations}
			\State $\hat K \dgets \exp.( - \varepsilon^{-1} \hat C )$
			\While{$\Vert \hat \rho - \hat \rho \dtimes \hat a \dtimes \hat K(\hat b \dtimes \hat \sigma ) \Vert_1  > \mathtt{tol} $} \Comment{$L^1$ error of the marginal condition}
				\State $\hat a \dgets \tone \ddiv \hat K(\hat b \dtimes \hat \sigma )$
				\State $\hat b \dgets \tone \ddiv \hat K (\hat a \dtimes \hat \rho )$
			\EndWhile
			\State \Return $\varepsilon \log. \hat a,  \varepsilon \log. \hat b$ \Comment{The Kantorovich potentials}
		\EndFunction
	\end{algorithmic}
\end{algorithm}
Algorithm \ref{alg:sinkhorn} exclusively relies on element-wise operations and matrix-vector products. The Matrix $\hat C$ can be evaluated lazily, which reduces the memory footprint from $N^2$ to $N$. However, the quantities $\hat K$, $\hat a$, and $\hat b$ become numerically unstable as $\varepsilon \rightarrow 0$. \\
A very similar algorithm can be used to calculate OT barycenters. In the case of the barycenters, even moderate values of $\varepsilon > 0$ cause a blurring of the result, which can be corrected by using a slightly modified algorithm (see also Section \ref{sec:entropic_bias}). Our implementation follows \cite{janati_debiased_2020} and we refer to this work for further details. \\

\subsection{Log-domain c-transform}

For many applications, moderate values of $\varepsilon$ might be sufficient (around $10^{-3}$ in our numerical examples) and Algorithm \ref{alg:sinkhorn} can be run without numerical over- or underflow. However, in our approach, there is one step that relies on $\varepsilon$ to be very small, say of order $10^{-6}$: The approximation of the c-transform to invert the transport mappings. For this, we rely on the softmin function \eqref{eq:softmin}:
\begin{equation}\label{eq:softmin_discrete}
	\hat \psi^{c,\varepsilon}_j := - \varepsilon \log \sum_{i=1}^N \exp \left (\frac{\hat \psi^\varepsilon_i - \hat C_{ij}}{\varepsilon} \right ) \hat \rho_i  \quad 1 \leq j \leq N.
\end{equation}
%\begin{algorithm}[h]
%	\begin{algorithmic}[1]
%		\Function{softmin}{$\tpsi,\tu,\tC,\varepsilon$}
%		\For{j = 1, \dots, N}
%			\State  $\mathtt{res}[j] = - \varepsilon \log \mathtt{sum} ( \exp . (( \tpsiu - \tC[:,j]  ) / \varepsilon ) \dtimes \tu )$  \Comment{":" used as in Julia, NumPy, ...}
%		\EndFor
%		\State \Return $\tpsiv$
%		\EndFunction
%	\end{algorithmic}
%\end{algorithm}
The non-zero weights $\hat \rho_i$ can be absorbed in the exponent as $\exp \log \hat \rho_i$, obtaining a \emph{LogSumExp} operation. This is a standard function implemented in many programming languages and can be calculated very accurately, without round-off errors, using only a single pass over $i$. We use the implementation from the package \texttt{LogExpFunctions.jl}\footnote{https://github.com/JuliaStats/LogExpFunctions.jl}.

\subsection{Seperable kernels}

Naive implementations of both the matrix-vector products from Algorithm \ref{alg:sinkhorn} and the LogSumExp evaluation in \eqref{eq:softmin_discrete} are of $\mathcal{O}(N^2)$ complexity due to nested loops over $i$ and $j$. Note that $N$ scales exponentially with the spatial dimension of the problem when it is discretized on a grid. However, in the special case of $c(x,y) = \tfrac{1}{2} |x - y|^2$ we are working with, we can do better, as pointed out in \cite{solomon_convolutional_2015}. Note that in dimension $d$, the cost is separable in $d$ terms along each dimension: $\tfrac{1}{2} |x - y|^2 = \tfrac{1}{2} |x^1 - y^1|^2 + \dots + \tfrac{1}{2} |x^d - y^d|^2$. Now assume the points $x_i$ are sampled on a regular tensor grid and therefore can be indexed as $x_{ i_1,\dots,i_d} : 1 \leq i_1,\dots,i_d \leq N^{1/d}$. Note that $x^l_{ i_1,\dots,i_d}$ can be denoted with $x^l_{i_l} $, as only the $l$th coordinate changes when varying the indices $i_1,\dots,i_d$. Let $\hat c^l_{ij} := |x_{i_l}^l - y_{j_l}^l|^2$ and $\hat k^l = \exp.(- \hat c^l / \varepsilon)$. In this case,
\begin{align}
	(\hat K \hat a )_j &= \sum_{1 \leq i \leq N} \hat K_{ij} \hat a_i \nonumber \\
	\Leftrightarrow (\hat K \hat a)_{j_1,\dots,j_d} &= \sum_{1 \leq i_1,\dots,i_d \leq N^{1/d}} \hat k_{i_1 j_1} \cdots \hat k_{i_d j_d} \hat a_{ i_1,\dots,i_d} \nonumber \\
	&= \sum_{1 \leq i_1 \leq N^{1/d}} \hat k_{i_1 j_1} \cdots  \sum_{1 \leq i_d \leq N^{1/d}} \hat k_{i_d j_d} \hat a_{ i_1,\dots,i_d}
\end{align}
As a result, instead of computing one large matrix-vector product of complexity $N^2$, we are computing $d$ tensor contractions of complexity $N^{1+1/d}$ each. An analogous trick can be applied in the log-domain as well.

In order to utilize this separable kernel trick, when it comes to the OT calculations, we sample the densities $\rho$ on a fine, regular tensor grid and perform the calculations as outlined above. This yields point-wise approximations to the Kantorovich potentials $\{ \psi(x_i), \psi^c(y_i) \}_{i = 1}^N$. In order to obtain the transport map $T$ and its Jacobian $DT$, we use this data to construct finite element approximations to $\psi$ and $\psi^c$ in $V_h$. Note that we could also use a different approximation space at this point.

\section{EIM algorithm}\label{sec:eim_appendix}

For convenience, we give the full algorithm to construct an empirical interpolation approximation from \cref{sec:hyperreduction}:

\begin{algorithm}[h]
	\caption{Empirical Interpolation Method algorithm}\label{alg:eim}
	\begin{algorithmic}[1]
		\Function{empiricalinterpolation}{$\{ \Xi_q \}_{q=1}^Q$}
		\For{$q = 1, \dots, Q$}
		\If{$q$ is $1$}
		\State $r \gets \Xi_q$
		\Else
		\State $\theta \gets B^{-1} \, [\Xi_q(y_1^{\mathrm{eim}}), \dots,  [\Xi_q(y_{q-1}^{\mathrm{eim}})] $
		\State $r \gets \Xi_q - \theta \cdot X$
		\EndIf
		\State $y_q^{\mathrm{eim}} \gets \argmax_{y \in \Omega} \Vert r(y) \Vert_\infty$
		\State $X_q \gets r_q / r_q(y_q^{\mathrm{eim}}) $
		\State $B_{q,q'} \gets X_q(y_{q'}^{\mathrm{eim}}) \quad \forall q' = 1, \dots, Q$ \Comment{ $B$ is lower-triangular with unit diagonal.}
		\EndFor
		\State \Return $\{ X_q \}_{q=1}^Q, \{ y_q^{\mathrm{eim}} \}_{q=1}^Q, B$ \Comment{The interpolation functions, points, and matrix.}
		\EndFunction
	\end{algorithmic}
\end{algorithm}

\section{An analytical example in one dimension}\label{sec:n-m-supp}

In this section we will show how the method proposed in \cref{sec:method} performs in a one-dimensional example where analytical solutions can be calculated. The example is taken from appendix B of \cite{taddei_registration_2020}.

\begin{prop}
	The solutions to the equation
	\begin{eqnarray}\label{eq:boundary_layer}
		- \partial_{xx}^2 u_\mu + \mu^2 u_\mu = 0
	\end{eqnarray}
	on the domain $\Omega = (0,1)$ with boundary conditions $u_\mu(0) = 1$, $u_\mu(1) = 0$ and $\mu, \bar \mu \in [ \mu_{\mathrm{min}}, \mu_{\mathrm{max}}] =: \mathcal{A}$, $\mu_{\mathrm{max}} = \epsilon^{-2} \mu_{\mathrm{min}}$, $\mu_{\mathrm{min}} > 1$, $\epsilon \in (0,1)$ satisfy
	\begin{equation}
		\Vert u_{\bar \mu} - u_\mu \circ T_{\rho_{\bar \mu} \rightarrow \rho_{\mu}} \Vert_{L^2(\Omega)} \leq \left \vert \frac{1}{\cosh \mu} - \frac{1}{\cosh \bar \mu} \right \vert \leq 2 e^{-\mu_{\mathrm{min}}},
	\end{equation}
	where $\rho(u) = \tfrac{u}{\int u}$.
\end{prop}

\begin{proof}
	The solution manifold is given by $\mathcal{M} = \{ \tfrac{ \cosh(\mu(1-x)) }{\cosh \mu} : \mu \in \mathcal{A} \}$. In this one-dimensional example, the OT maps can be calculated analytically using cumulative density functions. Since $\rho(u) = \tfrac{u}{\int u}$, $\rho(u_\mu) =: \rho_\mu = \mu \tfrac{ \cosh(\mu(1-x))}{\sinh \mu }$. Furthermore, $\mathrm{cdf}(\rho_\mu)(x) = 1 - \tfrac{\sinh(\mu(1-x))}{\sinh \mu}$, $\mathrm{cdf}(\rho_\mu)^{[-1]}(p) = 1 - \tfrac{1}{\mu} \sinh^{-1} ( (1-p) \sinh \mu )$, and
	\begin{equation}\label{eq:layer_map}
		T_{\rho_{\bar \mu} \rightarrow \rho_{\mu}}(y) = \mathrm{cdf}(\rho_\mu)^{[-1]} \circ \mathrm{cdf}(\rho_{\bar \mu}) (y) = 1 - \frac{1}{\mu} \sinh^{-1} \left ( \frac{\sinh \mu}{\sinh \bar \mu} \sinh(\bar \mu(1-y)) \right ).
	\end{equation}
	The map $T_{\rho_{\bar \mu} \rightarrow \rho_{\mu}}$ is a bijection as it is strictly increasing and $T_{\rho_{\bar \mu} \rightarrow \rho_{\mu}}( \partial \Omega) = \partial \Omega$. The former is a consequence of $0 < \rho_{\mu} < +\infty \; \forall \mu$.
	
	Using $\cosh \circ \sinh^{-1}(x) = \sqrt{1 + x^2}$, and letting $z = \bar \mu(1-y)$, we write
	\begin{align}\label{eq:layer_bound_1}
		\Vert u_{\bar \mu} - u_\mu \circ T_{\rho_{\bar \mu} \rightarrow \rho_{\mu}} \Vert_{L^2(\Omega)}^2 
		%&= \frac{1}{\bar{\mu}} \int_0^{\bar \mu} \frac{1 + \tfrac{\sinh \mu^2}{\sinh \bar \mu^2} \sinh z^2}{\cosh \mu^2} \left ( 1 - \frac{\cosh \mu}{\cosh \bar \mu} \frac{\cosh z}{\sqrt{1 + \tfrac{\sinh \mu^2}{\sinh \bar \mu^2} \sinh z^2}} \right )^2 \rd z \nonumber \\
		&= \frac{1}{\bar{\mu}} \int_0^{\bar \mu} \frac{1}{\cosh \mu^2} \left (  \frac{\cosh \mu}{\cosh \bar \mu} \cosh z -  \sqrt{ 1 + \frac{\sinh \mu^2}{\sinh \bar \mu^2} \sinh z^2 } \right )^2 \rd z
	\end{align}
	where $\sinh z^2$ denotes $(\sinh z)^2$. We can check using symbolic numerical software that the integrand has no local maximum, as any stationary condition at $z \neq 0$ requires either $\tfrac{\sinh \mu^2}{\sinh \bar \mu^2} < \tfrac{\cosh \mu}{\cosh \bar \mu} < \tfrac{\sinh \mu}{\sinh \bar \mu} \Leftrightarrow \tfrac{\sinh \mu}{\sinh \bar \mu} < \tfrac{\tanh \bar \mu}{\tanh \mu} < 1$ or the same relation will all inequalities reversed, either of which lead to contradiction for $\mu, \bar \mu > 0$. As the integrand vanishes at $z = \bar \mu$, its maximum value is attained at $z = 0$, i.e. $y = 1$, and we find
	\begin{equation}
		\Vert u_{\bar \mu} - u_\mu \circ T_{\rho_{\bar \mu} \rightarrow \rho_{\mu}} \Vert_{L^2(\Omega)}
		\leq \vert \Omega \vert \vert u_{\bar \mu}(1) - u_\mu \circ T_{\rho_{\bar \mu} \rightarrow \rho_{\mu}}(1) \vert = \left \vert \frac{1}{\cosh \bar \mu} - \frac{1}{\cosh \mu} \right \vert.
	\end{equation}
\end{proof}

\begin{rem}\label{rem:rho_bar}
	 If we chose $\bar \rho = \mathrm{OTBar}\{ \rho_\mu : \mu \in \mathcal{A} \}$ with uniform weights, we find $\bar \rho = \rho_{\bar \mu}$ with $\bar \mu = \tfrac{\mu_{\mathrm{max}} - \mu_{\mathrm{min}}}{\log \mu_{\mathrm{max}} - \log \mu_{\mathrm{min}}}$, the logarithmic mean. The choice made in \cite{taddei_registration_2020} is $\bar \mu = \sqrt{ \mu_{\mathrm{min}} \mu_{\mathrm{max}} }$.
\end{rem}
%
%	To obtain the bound on the $L^2$ error note that, again by construction, $\rho_{\bar \mu}(y) = ( \rho_{\mu} \circ \Phi_\mu^{-1}(y)) \partial_y \Phi_\mu^{-1}(y)$ on $\Omega$. Since $u_{\mu} = \tfrac{\tanh \mu}{\mu} \rho_\mu$, 
%	\begin{equation}
%		\Vert u_{\bar \mu} - u_\mu \circ \Phi_\mu^{-1} \Vert_{L^2(\Omega)}^2 = \left \Vert \frac{\tanh (\bar \mu)}{\bar \mu} ( \rho_\mu \circ \Phi_\mu^{-1} ) \partial_y \Phi_\mu^{-1} - \frac{\tanh \mu}{\mu} \rho_\mu \circ \Phi_\mu^{-1} \right \Vert_{L^2(\Omega)}^2.
%	\end{equation}
%	Substituting the expressions for $\rho_{\mu}$ and $\partial_y \Phi_\mu^{-1}$, using identities of the hyperbolic functions (e.g. $\cosh \circ \sinh^{-1}(x) = \sqrt{1 + x^2}$), and letting $z = \bar \mu(1-y)$, we arrive at
%
\begin{rem}
	 Note that $T_{\rho_{\bar \mu} \rightarrow \rho_{\mu}}(y) \approx T_{\rho_{\bar \mu} \rightarrow \rho_{\mu}}(0) + y \partial_y T_{\rho_{\bar \mu} \rightarrow \rho_{\mu}}(0) = \tfrac{\bar \mu}{\mu} \tfrac{\tanh \mu}{\tanh \bar \mu} y \approx  \tfrac{\bar \mu}{\mu} y$ and this approximation is close until either $y \approx \tfrac{\mu}{\bar \mu}$ (when $\mu < \bar \mu)$ or $y \approx 1$ (when $\mu > \bar \mu)$. 
	 
	The derivative of $T_{\rho_{\bar \mu} \rightarrow \rho_{\mu}}$ can take extreme values: $\partial_yT_{\rho_{\bar \mu} \rightarrow \rho_{\mu}}(1) =  \tfrac{\bar \mu}{\mu} \tfrac{\sinh \mu}{\sinh \bar \mu} \approx \tfrac{\bar \mu}{\mu} e^{\mu - \bar \mu}$. %Unfortunately no analytical expression can be given for the smoothed OT problem even in one dimension. %, but we recall that $\Vert \partial_y T^\varepsilon \Vert_\infty = \mathcal{O}(1 + \tfrac{1}{\epsilon})$ \cite{genevay_sample_2019}.
\end{rem}

We now give the proof of \cref{prop:boundary_layer}:

\begin{proof}
	We want to show that (c.f. \cref{eq:boundary_layer_nm})
		\begin{equation}
		\inf_{\xi^c_{1} \in \mathrm{span} \{ \psi^c_\mu : \mu \in \mathcal{A} \} } \; \sup_{\mu \in \mathcal A} \; \inf_{\substack{ w_1(\mu) : \Phi^{-1}(y) = y - w_1(\mu) \partial_y \xi^c_1(y) \\ \Phi_\mu^{-1}: \Omega \rightarrow \Omega \text{ is a bijection}}} \Vert u_{\bar \mu} - u_\mu \circ \Phi_{\mu}^{-1} \Vert_{L^2(\Omega)} \leq e^{-\mu_{\mathrm{min}}}( 4 + \epsilon). \nonumber
	\end{equation}
	By \cref{rem:rho_bar}, $\bar \rho = \rho_{\bar \mu}$ with $\bar \mu = \tfrac{\mu_{\mathrm{max}} - \mu_{\mathrm{min}}}{\log \mu_{\mathrm{max}} - \log \mu_{\mathrm{min}}}$. Let $c(\mu) = \tfrac{\bar \mu - \mu}{\mu}$ and $w(\mu) = - \tfrac{c(\mu)}{c(\mu_{\mathrm{min}})}$. We will show the bound by evaluating it at the trial function 
	\begin{equation}
		\Phi^{-1}_\mu(y) := y - w(\mu) \left ( T_{\rho_{\bar \mu} \rightarrow \rho_{\mu_{\mathrm{min}}}}(y) - y \right ),
	\end{equation}
	i.e. $\partial_y \xi^c_1(y) := T_{\rho_{\bar \mu} \rightarrow \rho_{\mu_{\mathrm{min}}}}(y) - y$. Note that, by the properties of the logarithmic mean, $\tfrac{\mu_{\mathrm{min}}}{\epsilon} \leq \bar \mu \leq \tfrac{\mu_{\mathrm{min}}}{2 \epsilon^2}$. As a consequence, $-1 \leq w(\mu) \leq \epsilon$. As $T_{\rho_{\bar \mu} \rightarrow \rho_{\mu_{\mathrm{min}}}}$ is concave, $\Phi^{-1}_\mu$ is concave for $\mu < \bar \mu$ and convex otherwise. Consequently, $\Phi^{-1}_\mu$ is strictly increasing, with $\partial_y \Phi^{-1}_\mu(y) \geq \partial_y \Phi^{-1}_{\mu_{\mathrm{min}}}(1) > \tfrac{\bar \mu}{\mu_{\mathrm{min}}} e^{\mu_{\mathrm{min}} - \bar \mu} $ for $\mu < \bar \mu$ and $\partial_y \Phi^{-1}_\mu(y) \geq \Phi^{-1}_{\mu_{\mathrm{max}}}(0) \geq \epsilon$ for $\mu \geq \bar \mu$.
	
	Let $\delta' := \tfrac{\mu_{\mathrm{min}}}{\bar \mu}$ and write
	\begin{multline}
		\Vert u_{\bar \mu} - u_\mu \circ \Phi_{\mu}^{-1} \Vert_{L^2(\Omega)}^2 \leq 
		\int_{0}^{\delta'} \left ( u_{\bar \mu}(y) - u_\mu \circ \frac{\bar \mu}{\mu} y \right )^2 \rd y \\
		+ \int_{0}^{\delta'} \left ( u_\mu \circ \Phi_{\mu}^{-1}(y) - u_\mu \circ \frac{\bar \mu}{\mu} y \right )^2 \rd y
		+ \int_{\delta'}^1 \left ( u_{\bar \mu}(y) - u_\mu \circ \Phi_{\mu}^{-1}(y) \right )^2 \rd y.
	\end{multline}

	For the first term, we can simplify $\vert u_{\bar \mu}(y) - u_\mu \circ \frac{\bar \mu}{\mu} y \vert = \sinh(\bar \mu y) \vert \tanh \bar \mu - \tanh \mu \vert \leq 2 \sinh( \mu_{\mathrm{min}}) \vert e^{-2 \mu} - e^{-2 \bar \mu} \vert \leq e^{-\mu_{\mathrm{min}}} \; \forall y \in [0,\delta']$.
	
	For the second term, since $\Phi^{-1}_\mu$ is either concave or convex, we find that $\vert \frac{\bar \mu}{\mu} y - \Phi^{-1}_\mu(y) \vert \leq \vert \frac{\bar \mu}{\mu} y - \partial_y \Phi^{-1}_\mu(0) y \vert \leq \vert \tfrac{\bar \mu}{\mu} - 1 + w(\mu) \tfrac{\bar \mu}{\mu_{\mathrm{min}}} \tfrac{\tanh \mu_{\mathrm{min}}}{\tanh \bar \mu} - w(\mu) \vert \tfrac{\mu_{\mathrm{min}}}{\bar \mu} = \vert w(\mu) \vert (1 - \tfrac{\tanh \mu_{\mathrm{min}}}{\tanh \bar \mu})  \; \forall y \in [0,\delta']$. Therefore,
	\begin{equation}
		\left \vert u_\mu \circ \Phi_{\mu}^{-1}(y) - u_\mu \circ \frac{\bar \mu}{\mu} y \right \vert_{0 \leq y \leq \delta'}  \leq \mu \underbrace{\max_{0 \leq y \leq \delta'} \frac{\sinh (\mu (1-y))}{\cosh \mu}}_{= \tanh \mu} \vert w(\mu) \vert \left (1 - \frac{\tanh \mu_{\mathrm{min}}}{\tanh \bar \mu} \right ).
	\end{equation}
	Using the bounds on $w(\mu)$, $\bar \mu$, and $1 - \tanh \mu \leq 2 e^{-2\mu}$, this expression is bounded by $\bar \mu \cdot 1 \cdot 1 \cdot 2 e^{-2\mu_{\mathrm{min}}} = 2 \bar \mu e^{-2\mu_{\mathrm{min}}}$ for $\mu \leq \bar \mu$ and $\mu_{\mathrm{max}} \cdot 1 \cdot \epsilon \cdot 2 e^{-2\mu_{\mathrm{min}}} \leq 2 \bar \mu e^{-2\mu_{\mathrm{min}}}$ otherwise.
	
	For the third term, assume $\mu \leq \bar \mu$. Recall that in this case, $\Phi_{\mu}^{-1}(y) \leq \tfrac{\bar \mu}{\mu} y$ and therefore $u_\mu \circ \tfrac{\bar \mu}{\mu} y \leq u_\mu \circ \Phi_{\mu}^{-1}(y)$. Since the mapping is increasing, the maximum of the integrand is reached at $y = \delta'$. We find the following chain of inequalities: $u_{\bar \mu}(\delta') = \cosh \mu_{\mathrm{min}} - \tanh \bar \mu \sinh \mu_{\mathrm{min}} \leq \cosh \mu_{\mathrm{min}} - \tanh \mu \sinh \mu_{\mathrm{min}} = u_{\mu} \circ \tfrac{\bar \mu}{\mu} \delta' \leq u_\mu \circ \Phi_{\mu}^{-1}(\delta')$. Lastly, using $T_{\rho_{\bar \mu} \rightarrow \rho_{\mu_{\mathrm{min}}}}(\delta') = 1 - \tfrac{1}{\mu_{\mathrm{min}}} \sinh^{-1} (  \tfrac{\sinh \mu_{\mathrm{min}}}{\sinh \bar \mu} \sinh (\bar \mu -  \mu_{\mathrm{min}}) ) \geq 1 - \tfrac{1}{\mu_{\mathrm{min}}}$, we arrive, after some simplifications, at $\mu(1 - \Phi_{\mu}^{-1}(\delta')) \leq \mu - \mu_{\mathrm{min}} + \tfrac{\bar \mu - \mu}{\bar \mu - \mu_{\mathrm{min}}} \leq \mu - \mu_{\mathrm{min}} + 1$ and $u_{\mu} \circ \Phi_{\mu}^{-1}(\delta') \leq \cosh(\mu_{\mathrm{min}} - 1 ) - \sinh(\mu_{\mathrm{min}} - 1 ) \tanh \mu = e^{1- \mu_{\mathrm{min}}} + \sinh (\mu_{\mathrm{min}} - 1 ) (1 - \tanh \mu)$ and therefore
	\begin{equation}
		\left \vert u_{\bar \mu} (\delta') - u_\mu \circ \Phi_{\mu}^{-1}(\delta') \right \vert \leq u_{\mu} \circ \Phi_{\mu}^{-1}(\delta') \leq e^{1 - \mu_{\mathrm{min}}} + e^{\mu_{\mathrm{min}} - 1} e^{- 2 \mu} \leq  e^{- \mu_{\mathrm{min}}} (e + e^{-1}).
	\end{equation}
	When $\mu \geq \bar \mu$, the reversed inequalities hold and we obtain $\left \vert u_{\bar \mu} (\delta') - u_\mu \circ \Phi_{\mu}^{-1}(\delta') \right \vert \leq u_{\bar \mu} \circ \Phi_{\bar \mu}^{-1}(\delta') = \cosh \mu_{\mathrm{min}} - \sinh \mu_{\mathrm{min}} \tanh \bar \mu \leq 2 e^{- \mu_{\mathrm{min}}}$.
	
	Collecting all terms, we find
	\begin{equation}
		\Vert u_{\bar \mu} - u_\mu \circ \Phi_{\mu}^{-1} \Vert_{L^2(\Omega)} \leq \delta' e^{-\mu_{\mathrm{min}}} + \delta' 2 \bar \mu e^{-2\mu_{\mathrm{min}}} + (1-\delta') (e + e^{-1}) e^{-\mu_{\mathrm{min}}} < e^{- \mu_{\mathrm{min}}} (4 + \epsilon),
	\end{equation}
	using $0 < \delta' = \tfrac{\mu_{\mathrm{min}}}{\bar \mu} < \epsilon$ and $\mu_{\mathrm{min}} > 1$.
\end{proof}

\bibliographystyle{siamplain}
\bibliography{../OTRB.bib}

\end{document}